\documentclass[12pt,a4paper]{amsart}

\usepackage{amsmath}
\usepackage{amssymb}
\usepackage{rotating}
\usepackage{hyperref,amscd}

\theoremstyle{thm}
\newtheorem{proposition}{\addtocounter{corollary}{1}\addtocounter{lemma}{1}\addtocounter{theorem}{1}Proposition}

\newtheorem{corollary}{\addtocounter{proposition}{1}\addtocounter{lemma}{1}\addtocounter{theorem}{1}Corollary}

\newtheorem{lemma}{\addtocounter{proposition}{1}\addtocounter{corollary}{1}\addtocounter{theorem}{1}Lemma}

\newtheorem{theorem}{\addtocounter{proposition}{1}\addtocounter{corollary}{1}\addtocounter{lemma}{1}Theorem}

\makeatletter
\@addtoreset{equation}{section}
\@addtoreset{proposition}{section}
\@addtoreset{corollary}{section}
\@addtoreset{lemma}{section}
\@addtoreset{theorem}{section}
\makeatother
\newcommand\Aut{\mathop{\mathrm{A\hspace{-0.2ex}u\hspace{-0.1ex}t}}\nolimits}
\newcommand\card{\mathop{\mathrm{card}}\nolimits}
\newcommand\Hom{\mathop{\mathrm{Hom}}\nolimits}
\newcommand\id{\mathop{\mathrm{id}}\nolimits}
\newcommand\IM{\mathop{\mathrm{Im}}\nolimits}
\newcommand\Ker{\mathop{\mathrm{Ker}}\nolimits}
\newcommand\sgn{\mathop{\mathrm{sgn}}\nolimits}
\newcommand\St{\mathop{\mathrm{St}}\nolimits}
\newcommand\supp{\mathop{\mathrm{supp}}\nolimits}
\newcommand\fg{\mathfrak g}
\newcommand\fm{\mathfrak m}
\newcommand\Nee{\mathbb N}
\newcommand\ugol{\mskip\medmuskip\llcorner\mskip\medmuskip}
\newcommand\im{\mathop{\mathrm{Im}}\nolimits}
\newcommand\Ham{\mathop{\mathrm{Ham}}\nolimits}
\newcommand\Cont{\mathop{\mathrm{Cont}}\nolimits}
\newcommand\Prim{\mathop{\mathrm{Prim}}\nolimits}

\newcommand\Acal{\mathcal A}
\newcommand\Bcal{\mathcal B}
\newcommand\Fcal{\mathcal F}
\newcommand\Gcal{{\mathcal G}}
\newcommand\Hcal{{\mathcal H}}
\newcommand\Jcal{{\mathcal J}}
\newcommand\Kcal{{\mathcal K}}
\newcommand\Lcal{{\mathcal L}}
\newcommand\Ocal{\mathcal O}
\newcommand\gf{\mathfrak g}
\newcommand\Dfrak{\mathfrak D}
\newcommand\Efrak{\mathfrak E}
\newcommand\Ffrak{\mathfrak F}
\newcommand\Gfrak{\mathfrak G}
\newcommand\Pfrak{\mathfrak P}
\newcommand\Qfrak{\mathfrak Q}
\newcommand\Qfrakov{\overline{\Qfrak}}
\newcommand\Rfrakov{\overline{\mathfrak R}}
\newcommand\elp{\vphantom{e}'\mskip-2mu e}
\newcommand\mOm{\fm\mskip1mu\Omega}
\newcommand\om{\omega}
\newcommand\Om{\Omega}
\newcommand\ufm{u\mskip1mu\fm}
\newcommand\uOcal{u\mskip1mu\Ocal}
\newcommand\uOm{u\mskip1mu\Omega}
\newcommand\umOm{u\mskip1mu\fm\mskip1mu\Omega}
\newcommand\KtmE{K^\times\mskip-5mu\times\mskip-2mu E}
\newcommand\KtmVd{K^\times\mskip-5mu\times\mskip-2mu V^*}
\newcommand\OrtmG{\Ocal^\times\mskip-5mu\rtimes G}
\newcommand\OFrtmGF{\Ocal(\Fcal)^\times\mskip-4mu\rtimes G(\Fcal)}
\newcommand\LaE{{\textstyle\bigwedge^2\!E}}
\newcommand\LaEnd{{\textstyle\bigl(\LaE\bigr){}_\#}}
\newcommand\LaVd{{\textstyle\bigwedge^2}\mskip2mu V^*}
\newcommand\LaVdnd{{\textstyle\bigl(\LaVd\bigr){}_\#}}
\newcommand\Ma{\setbox0\hbox{$\bigwedge$}\hbox{\copy0\hskip-.65\wd0\box0}}
\newcommand\MaE{\Ma E}
\newcommand\MaVd{\Ma V^*}
\newcommand\FNup{\Ffrak\strut_\Nee^{\mskip2mu\uparrow}}
\newcommand\FNinftyup{\Ffrak\strut_{\Nee,\infty}^{\mskip2mu\uparrow}}
\newcommand\Finftyup{\Ffrak\strut_\infty^{\mskip2mu\uparrow}}
\newcommand\FNdown{\Ffrak\strut_\Nee^{\mskip2mu\downarrow}}
\newcommand\FNinftydown{\Ffrak\strut_{\Nee,\infty}^{\mskip2mu\downarrow}}
\newcommand\Finftydown{\Ffrak\strut_\infty^{\mskip2mu\downarrow}}
\newcommand\iipair{\vphantom{I}_{\mskip1mu\iov}^{\mskip2mu\iov}}
\newcommand\ijpair{\vphantom{I}_{\mskip1mu\jov}^{\mskip2mu\iov}}
\newcommand\jipair{\vphantom{I}_{\mskip1mu\iov}^{\mskip2mu\jov}}
\newcommand\iov{\overline i}
\newcommand\jov{\overline j}
\newcommand\ijov{\iov,\mskip1mu\jov}
\newcommand\kkpair{\vphantom{I}_{\mskip1mu\kov}^{\mskip1mu\kov}}
\newcommand\klpair{\vphantom{I}_{\mskip1mu\lov}^{\mskip1mu\kov}}
\newcommand\lkpair{\vphantom{I}_{\mskip1mu\kov}^{\mskip1mu\lov}}
\newcommand\kov{\overline k}
\newcommand\lov{\overline\ell}
\newcommand\klov{\kov,\mskip1mu\lov}
\newcommand\tausiP{\tau\sigma\mskip-1mu P}
\newcommand\siP{\sigma\mskip-1mu P}
\newcommand\tauP{\tau\mskip-1mu P}
\newcommand\repl[1]{{\csname xx#1\endcsname}}

\expandafter\newcommand\csname xxNormal shape\endcsname{Normal shape}
\expandafter\newcommand\csname xxThe normal shape\endcsname{The normal shape}
\expandafter\newcommand\csname xxTHE NORMAL SHAPES\endcsname{THE NORMAL SHAPES}
\expandafter\newcommand\csname xxNORMAL SHAPES\endcsname{NORMAL SHAPES}

\begin{document}
 
\thispagestyle{empty}
\title[\repl{Normal shapes} of \repl{symplectic} and contact forms]{\repl{The normal shapes} of the \repl{symplectic} and contact forms over algebras of divided powers}
\author{S.M. Skryabin}

\begin{abstract}
The Lie algebras of Cartan type over the fields of positive characteristic are para\-met\-rized by differential forms with coefficients in the divided power algebras. The normal shapes of the volume forms have been classified by Tyurin 
and Wilson. 
In this paper, a similar problem is solved for the contact and \repl{symplectic} forms.
\end{abstract}

\maketitle

\markboth{\itshape Serge Skryabin}{{\itshape Normal shapes of symplectic and contact forms}}

\thispagestyle{empty}

\setcounter{tocdepth}{2}
\tableofcontents

\section*{Preface to the Translation}

This text is a translation into English of a manuscript written in Russian and deposited in 1986 in VINITI archive. The manuscript gives a classification of the differential forms that parametrize the finite-dimensional Lie algebras of alternating hamiltonian and contact types over fields of characteristic $p>0$. Part of the material covered by this manuscript was published in a journal article \cite{Sk}. 
It was shown there how the classification of 
symplectic forms reduces to certain questions in the linear, or rather semilinear, algebra. For the 
symplectic forms of the 2nd type the classification is then completed quickly. An exposition of these results can be found also in the book of H.~Strade \cite{St}. 

As to the 
symplectic forms of the 1st type one eventually has to classify objects consisting of a pair of finite-dimensional vector spaces, alternating bilinear forms and flags of subspaces on the two spaces, and a collection of semilinear isomorphisms between the factors of the two flags. This problem is solved in \ref{6} of the manuscript, and the details have not been published elsewhere.  

The manuscript explains also the classification of contact forms.

It has been the initiative of Dimitry Leites to bring back to life this almost 
forgotten manuscript. The result has come as a pleasant surprise to me. I am 
very grateful to Dimitry Leites for translating, 
and to Dimitry Vorobiev for \TeX-coding the manuscript. 
Several minor errors in the original text have been corrected by myself.

The translator has made several changes to the terminology used in the 
manuscript and other papers on modular Lie algebras. The term ``\emph{hamiltonian 
form}'' has been replaced by ``\emph{symplectic form}'' speaking of the non-degenerate closed differential 2-form, 
and the expression ``\emph{canonical form}\/ of something'' by ``\emph{normal shape}\/ of something''.

\medskip
The class of Cartan type Lie algebras in the setting of the characteristic $p>0$ theory was introduced in a 1969 paper of A.I.~Kostrikin and I.R.~Shafarevich \cite{KSh}. Such a Lie algebra consists of special derivations of the free divided power algebra $\Ocal(E)$ on a finite-dimensional vector space $E$, and is defined with respect to a decreasing flag $\Fcal$ of subspaces of $E$ and a certain differential form $\omega$. Actually the paper \cite{KSh} dealt only with the Lie algebras of Cartan type that admit a standard grading. For their generalization one has to consider a larger class of differential forms with coefficients in the completion $\widehat\Ocal(E)$ of $\Ocal(E)$, but still related in some way to flags of subspaces of $E$.

The algebra $\Ocal(E)$ is presented by generators $\{x^{(r)}\mid x\in E,\ r=0,1,\ldots\}$ and relations
\begin{gather*}
x^{(0)}=1,\qquad
(x+y)^{(r)}=\sum_{i=0}^rx^{(i)}y^{(r-i)},\qquad
(\lambda x)^{(r)}=\lambda^rx^{(r)},\\
\hspace{1cm}x^{(r)}x^{(s)}=\begin{pmatrix}r+s\\r\end{pmatrix}x^{(r+s)},
\end{gather*}
where $\,x,y\in E$, $\,\lambda\in K$ (the ground field), and $r,s\in\Nee$. The subalgebra $\Ocal(\Fcal)$ associated with the flag $\Fcal$ (see \ref{1}) is generated by
$$
\{x^{(r)}\mid x\in E_i,\ r<p^{i+1},\ i=0,1,\ldots\}.
$$
The ideals $\fm(E)^{(k)}$ of $\Ocal(E)$ are defined inductively by the relation
$$
\fm(E)^{(k)}=\sum_{x\in E,\ r\geqslant k}\Ocal(E)\,x^{(r)}+\sum_{0<j<k}\fm(E)^{(j)}\,\fm(E)^{(k-j)}\quad\textrm{for }\,k>0
$$
starting with $\fm(E)^{(0)}=\Ocal(E)$. The ideals $\fm(\Fcal)^{(k)}$ of $\Ocal(\Fcal)$ are similarly defined with the use of only those divided powers $x^{(r)}$ that lie in $\Ocal(\Fcal)$.

A derivation $D$ of $\Ocal(E)$ is called \textit{special} \cite{KSh} (or \textit{distinguished} by Leites and co-authors) if
$$
D(x^{(r)})=x^{(r-1)}D(x)\quad\text{for all }\,x\in E,\ r>0.
$$
A derivation $D$ of $\Ocal(\Fcal)$ is called \textit{special} if the equalities above hold whenever $x^{(r)}\in\Ocal(\Fcal)$. In this case $D$ extends in a unique way to a special derivation of $\Ocal(E)$, and uniquely extends to a continuous derivation of the completed algebra $\widehat\Ocal(E)$ defined below.

Other relevant notions will be defined in the main body of the paper.

\section{Divided power algebras, differential forms, cohomology,\\{ }automorphism groups}\label{1}

Let $K$ be a field of characteristic $p>0$. Let $E$ be a finite-dimensional vector space over $K$ and
$$
\Fcal:E=E_0\supseteq E_1\supseteq\ldots\supseteq 0=E_m=E_{m+1}=\ldots
$$
a flag of subspaces. Consider the following objects defined in \cite{KSh}: the divided power algebra $\Ocal(E)$, its subalgebra $\Ocal(\Fcal)$, their standard filtrations $\fm(E)^{(k)}$ and $\fm(\Fcal)^{(k)}$, the completion $\widehat{\Ocal}(E)$ of $\Ocal(E)$ with respect to this filtration and the closure $\widehat\fm(E)^{(k)}$ of the ideal $\fm(E)^{(k)}$ in $\widehat{\Ocal}(E)$. The ideal $\widehat\fm(E)=\widehat\fm(E)^{(1)}$ (resp. $\fm(\Fcal)=\fm(\Fcal)^{(1)}$) of the algebra $\widehat{\Ocal}(E)$ (resp. $\Ocal(\Fcal)$) is maximal and $f^p=0$ for any $f\in\widehat\fm(E)$.

Let $W(\Fcal)$ (resp. $\widehat W(E)$) be the Lie algebra of \textit{\repl{distinguished}} (resp. \textit{continuous \repl{distinguished}}) derivations of the algebra $\Ocal(\Fcal)$ (resp. $\widehat{\Ocal}(E)$). Let $\widehat\Om(E)$ and $\Om(\Fcal)$ be the respective complexes of differential forms. The terms of the standard filtration in $\widehat W(E)$ and $W(\Fcal)$ are as follows:
$$
\widehat W(E)_k=\widehat\fm(E)^{(k+1)}\,\widehat W(E),\qquad W(\Fcal)_k=\fm(\Fcal)^{(k+1)}\,W(\Fcal).
$$

Let $x_1,\ldots,x_n$ be a basis of $E$ compatible with the flag, and let $\partial_i\in\widehat W(E)$ be such that $\partial_ix_j=\delta_{ij}$ (the Kronecker symbol)\footnote{In other words, $\partial_i$ is the \textit{partial derivative} whose action is given by the formula $\partial_ix_j^{(k)}=\delta_{ij}x_j^{(k-1)}$ for any $k$.}. For $i=1,2,\ldots,n$ we set $m_i=\min\limits\{k\mid x_i\notin E_k\}$. We will consider the totality of $k$-forms as the set of $\widehat{\Ocal}(E)$-multilinear (resp. $\Ocal(\Fcal)$-multilinear) \repl{antisymmetric}\footnote{In the original manuscript the word ``skewsymmetric" was used with the understanding that these are mappings which take zero value whenever any pair of arguments are equal. However, ``alternating" is the standard term for such mappings, while ``skewsymmetric" should have the same meaning as ``symmetric" when $p=2$. For this reason ``alternating" replaces ``skewsymmetric" throughout the whole paper.} mappings
$$
\underbrace{\widehat W(E)\times\ldots\times\widehat W(E)}_k\to\widehat{\Ocal}(E)\quad\text{(resp. }\underbrace{W(\Fcal)\times\ldots\times W(\Fcal )}_k\to\Ocal(\Fcal)\text{)}.
$$
We identify $\Om^k(\Fcal)$ with the subspace of the forms $\om\in\widehat\Om^k(E)$ such that 
$$
\om(\delta_1,\ldots,\delta_k)\in\Ocal(\Fcal)\quad\text{for all}\quad\delta_1,\ldots,\delta_k\in W(\Fcal).
$$
The differential $d$ and the multiplication of forms are defined by the relations 
\begin{gather}
\begin{split}
\label{(1.1)}
d\om(\delta_1,\ldots,\delta_{k+1})
&=\sum_{i=1}^{k+1}(-1)^{i-1}\delta_i\,\om(\delta_1,\ldots,\widehat\delta_i,\ldots,\delta_{k+1})+{}\\
&\hskip3cm+\sum_{i<j}(-1)^{i+j}\om([\delta_i,\delta_j],\ldots,\widehat\delta_i,\ldots,\widehat\delta_j,\ldots),
\end{split}\\
\label{(1.2)}
(\om\wedge\om')(\delta_1,\ldots,\delta_{k+\ell})=\sum\nolimits'\sgn(\pi)\,\om(\delta_{\pi(1)},\ldots,\delta_{\pi(k)})\,\om'(\delta_{\pi(k+1)},\ldots,\delta_{\pi(k+\ell)}),
\end{gather}
where $\om\in\widehat\Om^k(E)$, $\om'\in\widehat\Om^\ell(E)$ and the sum in \eqref{(1.2)} runs over all the permutations $\pi$ of the set $\{1,\ldots,k+\ell\}$ such that 
$\,\pi(1)<\ldots<\pi(k)$ and $\,\pi(k+1)<\ldots<\pi(k+\ell)$.

For every $\delta\in\widehat W(E)$ and $\om\in\widehat\Om^k(E)$, define $\delta\om\in\widehat\Om^k(E)$ and $\delta\ugol\om\in\widehat\Om^{k-1}(E)$ by the formulas 
\begin{gather}
\label{(1.3)}
(\delta\om)(\delta_1,\ldots,\delta_k)=\delta\,\om(\delta_1,\ldots,\delta_k)-\sum_{k=1}^k\om(\delta_1,\ldots,[\delta,\delta_i],\ldots,\delta_k),\\
\label{(1.4)}
(\delta\ugol\om)(\delta_1,\ldots,\delta_{k-1})=\om(\delta,\delta_1,\ldots,\delta_{k-1}).
\end{gather}
The following identities are well-known
\begin{gather}
\label{(1.5)}
d\delta\om=\delta d\om,\\
\label{(1.6)}
\delta\om=\delta\ugol d\om+d(\delta\ugol\om),\\
\label{(1.7)}
d(\om\wedge\om')=d\om\wedge\om'+(-1)^{\deg\om}\om\wedge d\om',\\
\label{(1.8)}
\delta(\om\wedge\om')=\delta\om\wedge\om'+\om\wedge\delta\om',\\
\label{(1.9)}
\delta\ugol(\om\wedge\om')=(\delta\ugol\om)\wedge\om'+(-1)^{\deg\om}\om\wedge(\delta\ugol\om').
\end{gather}

Let $Z^k(\widehat\Om(E))$, resp. $B^k(\widehat\Om(E))$, resp. $H^k(\widehat\Om(E))$ denote the spaces of \textit{closed}, resp. \textit{exact} $k$-forms, resp. \textit{$k$-cohomology} of the complex $\widehat\Om(E)$ 
The notation $Z^k(\Om(\Fcal))$, $B^k(\Om(\Fcal))$, $H^k(\Om(\Fcal))$ has a similar meaning. 

In \cite{10} continuous mappings $\widehat\fm\to\widehat{\Ocal}(E)$, $f\mapsto f^{(r)}$ for $r=0,1,\ldots$ have been constructed such that
\begin{gather}
\label{(1.10)}
(x_1^{(\alpha_1)}\ldots x_n^{(\alpha_n)})^{(r)}=\prod_{i=1}^n\frac{(r\alpha_i)!}{r!\,(\alpha_i!)^r}\,x_1^{(r\alpha_1)}\ldots x_n^{(r\alpha_n)},\\
\label{(1.11)}
f^{(0)}=1,\\
\label{(1.12)}
(f+g)^{(r)}=\sum_{i=0}^rf^{(i)}g^{(r-i)},\\
\label{(1.13)}
(\lambda f)^{(r)}=\lambda^rf^{(r)},\\
\label{(1.14)}
f^{(r)}f^{(s)}=\begin{pmatrix}r+s\\r\end{pmatrix}f^{(r+s)},\\
\label{(1.15)}
(fg)^{(r)}=0\quad\text{for}\quad r\geqslant p,\\
\label{(1.16)}
(f^{(r)})^{(s)}=\frac{(rs)!}{(r!)^ss!}\,f^{(rs)},\\
\label{(1.17)}
f^{(r)}=\widehat\fm(E)^{(rj)},\quad\text{if}\ f\in\widehat\fm(E)^{(j)},
\end{gather}
for any $f,g\in\widehat\fm(E)$, $\lambda\in K$, and $r,s,\alpha_1,\ldots,\alpha_n\geqslant 0$. Note a particular case of formula \eqref{(1.16)}:
\begin{equation}
\label{(1.18)}
(f^{(p^k)})^{(p^\ell)}=f^{(p^{k+\ell})}.
\end{equation}
From \eqref{(1.10)} and continuity of raising to the divided powers it follows that
\begin{equation}
\label{(1.19)}
df^{(r)}=f^{(r-1)}df.
\end{equation}
If $f\in\fm(\Fcal)$, then, generally speaking, $f^{(r)}\notin\Ocal(\Fcal)$. We set
\begin{equation}
\label{(1.20)}
C_k(\Fcal)=\{f\in\fm(\Fcal)\mid f^{(p^k)}\in\Ocal(\Fcal)\}.
\end{equation}
In particular, $C_k(\Fcal)\cap E=E_k$.

\begin{proposition}[\cite{10}]\label{1.1}The space $C_k(\Fcal)$ has a basis consisting of monomials $x_1^{(\alpha_1)}\ldots x_n^{(\alpha_n)}$ with $\alpha_i<p^{m_i}$ which are not of the form $x_i^{(p^\ell)}$, where $\ell+k\geqslant m_i$.
\end{proposition}

\begin{corollary}\label{1.2}$\fm(\Fcal)^2\subseteq C_k(\Fcal)$ for all $k$.
\end{corollary}

Indeed, the monomials $x_1^{(\alpha_1)}x_2^{(\alpha_2)}\ldots x_n^{(\alpha_n)}$ with $\alpha_i<p^{m_i}$, other than $x_i^{(p^\ell)}$, form a basis of $\fm(\Fcal)^2$.

Let $\widehat{\Ocal}(E)^\times$ (resp. $\Ocal(\Fcal)^\times$) be the group of invertible elements of the algebra $\widehat{\Ocal}(E)$ (resp. $\Ocal(\Fcal)$). For each $f\in\widehat\fm(E)$ we set
\begin{equation}
\label{(1.21)}
\exp f=\sum_{i=0}^\infty f^{(i)}\in 1+\widehat\fm(E)\subseteq\widehat{\Ocal}(E)^\times.
\end{equation}
We have\quad$\displaystyle d(\exp f)=\sum df^{(i)}=\sum f^{(i-1)}df=(\exp f)\,df$,\quad i.e., 
\begin{equation}\label{(1.22)}
(\exp f)^{-1}d(\exp f)=df.
\end{equation}

\begin{proposition}\label{1.3}
The mapping $u\mapsto u^{-1}du$ gives a surjective homomorphism of the multiplicative group $\widehat {\Ocal }(E)^{\times }$ onto the additive group $Z^1(\widehat\Om(E)),$ and the group $K^\times$ is its kernel.
\end{proposition}

\noindent
\textbf{Proof.}
The fact that the mapping in question is a homomorphism follows from the formula
\begin{equation}
\label{(1.23)}
(uv)^{-1}d(uv)=u^{-1}du+v^{-1}dv.
\end{equation}
Since $d(u^{-1})=-u^{-2}du$, it follows that $u^{-1}du\in Z^1(\widehat\Om(E))$. Surjectivity follows from (1.22) and the fact that $H^1(\widehat\Om (E))=0$. Finally, if $u^{-1}du=0$, then $du=0$, implying $u\in K$.\qed

\begin{proposition}\label{1.4}
The mapping $\exp$ gives an isomorphism of the additive group\/ $\widehat\fm(E)$ onto the multiplicative group $1+\widehat\fm(E)$. We have $\exp(\fm(\Fcal)^i)=1+\fm(\Fcal)^i$ for any $i\geqslant 2$.
\end{proposition}

\noindent
\textbf{Proof.} The restriction $\varphi$ of the homomorphism $u\mapsto u^{-1}du$ to the subgroup $1+\widehat\fm(E)$ is an isomorphism onto $Z^1(\widehat\Om (E))$. Due to \eqref{(1.22)} we have $\exp(f)=\varphi^{-1}(df)$. Since $d$ maps $\widehat\fm(E)$ onto $Z^1(\widehat\Om (E))$ isomorphically, it follows that $\exp$ maps $\widehat\fm(E)$ onto $1+\widehat\fm(E)$ isomorphically.

If $f\in\fm(\Fcal)^i$ is such that $f=gh$,  where $g,h\in\fm(\Fcal)$, then taking into account formula \eqref{(1.15)} we have
$$
\exp(f)=\sum_{r=0}^{p-1}f^{(r)}=\sum_{r=0}^{p-1}\frac 1{r!}\,f^r\in 1+\fm(\Fcal)^i,
$$
and, moreover,
\begin{equation}
\label{(1.24)}
\exp(f)\equiv 1+f\mod{(1+\fm(\Fcal)^{i+1})}.
\end{equation}
Since $\exp$ is a homomorphism, it follows that $\exp(f)\in 1+\fm(\Fcal)^i$ and congruence \eqref{(1.24)} holds for any $f\in\fm(\Fcal)^i$. For $N$ sufficiently large we have $\fm(\Fcal)^N=0$. Now, by the downward induction on $i$ we prove with the aid of \eqref{(1.24)} that $\,1+\fm(\Fcal)^i\subseteq\exp(\fm(\Fcal)^i)$.\qed

\medskip
Denote by $\widehat G(E)$ the group of continuous automorphisms $\sigma$ of the algebra $\widehat{\Ocal}(E)$ such that 
$$
(\sigma f)^{(r)}=\sigma f^{(r)}\quad\text{for all}\quad f\in\widehat\fm(E),
$$
by $G(\Fcal)$ the group of automorphisms $\sigma$ of the algebra $\Ocal(\Fcal)$ such that $f^{(r)}$ and $(\sigma f)^{(r)}$ for any $f\in\fm(\Fcal)$ simultaneously either belong or do not belong to $\Ocal(\Fcal)$, and if $f^{(r)},(\sigma f)^{(r)}\in\Ocal(\Fcal)$, then $(\sigma f)^{(r)}=\sigma f^{(r)}$. 

Finally, let $G'(\Fcal)$ be the subgroup of automorphisms $\sigma\in G(\Fcal)$ such that $\sigma f-f\in\fm(\Fcal)^2$ for all $f\in\Ocal(\Fcal)$.

The group $\widehat G(E)$ leaves the spaces $\widehat\fm(E)^{(k)}$ invariant, the group $G(\Fcal)$ leaves the spaces $\fm(\Fcal)^{(k)}$ and $C_k(\Fcal)$  invariant. Define filtrations of subgroups setting for $j\geqslant 0$

\begin{gather*}
\widehat G(E)_j=\{\sigma\in\widehat G(E)\mid\sigma f-f\in\widehat\fm(E)^{(j+\ell)}\quad\text{for all}\quad f\in\widehat\fm(E)^{(\ell)},\ \ell\geqslant 0\},\\
G(\Fcal)_j=\{\sigma\in G(\Fcal)\mid\sigma f-f\in\fm(\Fcal)^{(j+\ell)}\quad\text{for all}\quad f\in\fm(\Fcal)^{(\ell)},\ \ell\geqslant 0\},\\
G'(\Fcal)_j=G(\Fcal)_j\cap G'(\Fcal).
\end{gather*}
The groups $G(\Fcal)_j$ and $G'(\Fcal)_j$ are normal in $G(\Fcal)$.

\begin{proposition}[\cite{10}]\label{1.5}
The group $\widehat G(E)$ consists of all continuous automorphisms $\sigma$ of the algebra $\widehat{\Ocal}(E)$ such that $\,\sigma\,\widehat W(E)\,\sigma^{-1}=\widehat W(E)$. Similarly,
$$
G(\Fcal)=\{\sigma=\Aut\Ocal(\Fcal)\mid\sigma\,W(\Fcal)\,\sigma^{-1}=W(\Fcal)\}.
$$
\end{proposition}

Therefore, we can define an action of $\widehat G(E)$ in the complex $\widehat\Om (E)$ and an action of $G(\Fcal)$ in the complex $\Om(\Fcal)$ by setting 
\begin{equation}
\label{(1.25)}
(\sigma\om)(\delta_1,\ldots,\delta_k)=\sigma\om(\sigma^{-1}\delta_1\sigma,\ldots,\sigma^{-1}\delta_k\sigma).
\end{equation}
Directly from the definition we derive the identities
\begin{gather}
\label{(1.26)}
\sigma d\om=d\sigma\om,\\
\label{(1.27)}
\sigma(\om\wedge\om')=\sigma\om\wedge\sigma\om'.
\end{gather}

\begin{proposition}[\cite{10}]\label{1.6}
The mapping $\sigma\mapsto(\sigma x_1,\sigma x_2,\ldots,\sigma x_n)$ establishes a bijection of $\widehat G(E)$ onto the set of tuples  $(y_1,y_2,\ldots,y_n)$ such that $y_i\in\widehat\fm(E)$ and\/ $\det(\partial_iy_j)_{1\leqslant i,j\leqslant n}$ is invertible in $\widehat{\Ocal}(E)$.
The mapping $\sigma\mapsto(\sigma x_1,\sigma x_2,\ldots,\sigma x_n)$ establishes a bijection of $G(\Fcal)$ onto the set of tuples $(y_1,y_2,\ldots,y_n)$ such that $y_i\in C_{m_i-1}(\Fcal)$ and\/ $\det(\partial_iy_j)_{1\leqslant i,j\leqslant n}$ is invertible in $\Ocal(\Fcal)$.
\end{proposition}

\begin{corollary}[\cite{10}]\label{1.7}
The group $G(\Fcal)$ is identified with the subgroup of automorphisms in $\widehat G(E)$ leaving the algebra $\Ocal(\Fcal)$ invariant.
\end{corollary}

The embedding $G(\Fcal)\subseteq\widehat G(E)$ allows us to consider the action of $G(\Fcal)$ on $\widehat\Om(E)$.

\begin{proposition}\label{1.8}
The mapping $\sigma\mapsto(\sigma x_1,\sigma x_2,\ldots,\sigma x_n)$ establishes a bijection of $G'(\Fcal)$ with the set of tuples $(y_1,y_2,\ldots,y_n)$ such that $y_i-x_i\in\fm(\Fcal)^2$. Under this bijection the subgroup $G'(\Fcal)_j$ corresponds to the tuples such that $y_i-x_i\in\fm(\Fcal)^2\cap\fm(\Fcal)^{(j+1)}$.
\end{proposition}

\noindent
\textbf{Proof.}
If $\sigma\in G'(\Fcal)_j$, then $\sigma x_i-x_i\in\fm(\Fcal)^2$ by definition of $G'(\Fcal)$, and $\sigma x_i-x_i\in\fm(\Fcal)^{(j+1)}$ by definition of $G(\Fcal)_j$. 

Conversely, let there be given a tuple $(y_1,y_2,\ldots,y_n)$ such that 
$y_i-x_i\in\fm(\Fcal)^2\cap\fm(\Fcal)^{(j+1)}$. Proposition~\ref{1.1} and Corollary~\ref{1.2} show that $y_i\in C_{m_i-1}(\Fcal)$.

Further, $\partial_iy_j\equiv\delta_{ij}\bmod{\fm(\Fcal)}$, whence $\det(\partial_iy_j)\equiv 1\bmod{\fm(\Fcal)}$, i.e., $\det(\partial_iy_j)$ is invertible. Therefore, there exists a unique $\sigma\in G(\Fcal)$ which sends $x_i$ to $y_i$, and we have to prove that $\sigma\in G'(\Fcal)_j$. For any $r<p^{m_i}$ we have 
$$
\sigma(x_i^{(r)})-x_i^{(r)}=y_i^{(r)}-x_i^{(r)}=\sum_{s=1}^r(y_i-x_i)^{(s)}x_i^{(r-s)}.
$$
Since $y_i-x_i\in\fm(\Fcal)^2$, it follows that $(y_i-x_i)^{(s)}\in\fm(\Fcal)^2$ for $s\geqslant 1$, and so $\sigma x_i^{(r)}-x_i^{(r)}\in\fm(\Fcal)^2$. 

On the other hand, $y_i-x_i\in\fm(\Fcal)^{(j+1)}$ implies 
$$
(y_i-x_i)^{(s)}x_i^{(r-s)}\in\fm(\Fcal)^{(s(j+1)+r-s)}\subseteq\fm(\Fcal)^{(r+j)}\quad\text{when}\quad s\geqslant 1,
$$
whence $\sigma x_i^{(r)}-x_i^{(r)}\in\fm(\Fcal)^{(j+r)}$. Set 
$$
A_k=\{f\in\fm(\Fcal)^{(k)}\mid\sigma f-f\in\fm(\Fcal)^2\cap\fm(\Fcal)^{(j+k)}\}.
$$
It has been proved that $x_i^{(r)}\in A_r$. We claim that $A_kA_\ell\subseteq A_{k+\ell}$. Indeed, if $f\in A_k$ and $g\in A_\ell$, then 
$$
\sigma(fg)-fg=(\sigma f-f)\sigma g+f(\sigma g-g)\in\fm(\Fcal)^2\cap\fm(\Fcal)^{(j+k+\ell)}.
$$
Now it is clear that $x_1^{(\alpha_1)}x_2^{(\alpha_2)}\ldots x_n^{(\alpha_n)}\in A_{\alpha_1+\alpha_2+\ldots +\alpha_n}$, i.e., 
$$
\sigma(x_1^{(\alpha_1)}\ldots x_n^{(\alpha_n)})-x_1^{(\alpha_1)}\ldots x_n^{(\alpha_n)}\in\fm(\Fcal)^2\cap\fm(\Fcal)^{(j+\alpha_1+\ldots +\alpha_n)},
$$
which means that $\sigma\in G'(\Fcal)_j$.\qed 

\medskip
The space $E$ is not stable under the $G(\Fcal)$-action on $\Ocal(\Fcal)$. Let us transfer the $G(\Fcal)$-module structure to $E$ by means of the isomorphism $E\cong\fm(\Fcal)/\fm(\Fcal)^{(2)}$. Thus, if $\sigma*x$, where $\sigma\in G(\Fcal)$ and $x\in E$, denotes the respective action, then 
\begin{equation}
\label{(1.28)}
\sigma x\equiv\sigma*x\mod{\fm(\Fcal)^{(2)}}\quad\text{for any}\quad\sigma\in G(\Fcal),\quad x\in E.
\end{equation}
The differential $d$ induces a $G(\Fcal)$-equivariant isomorphism 
$$
\fm(\Fcal)/\fm(\Fcal)^{(2)}\to\Om^1(\Fcal)/\fm(\Fcal)\Om^1(\Fcal), 
$$
whence we get an isomorphism of $G(\Fcal)$-modules $E\cong\Om^1(\Fcal)/\fm(\Fcal)\Om^1(\Fcal)$.  Proposition~\ref{1.6} shows that $G(\Fcal)$ acts on $E$ as the group of all linear transformations preserving the flag $\Fcal$.

Let us study the action of $G(\Fcal)$ in the cohomology groups. Note that the multiplication of differential forms induces an algebra structure in the space $H(\Om(\Fcal))$.

\begin{proposition}[\cite{kry}]\label{1.9}
The algebra $H(\Om(\Fcal))$ is isomorphic to the exterior algebra of the space $H^1(\Om(\Fcal))$. The classes of cocycles $x_i^{(p^{m_i}-1)}dx_i,$ where $i=1,2,\ldots,n,$ form a basis of $H^1(\Om(\Fcal))$.
\end{proposition}

We will define a $G(\Fcal)$-invariant decreasing filtration in $H^1(\Om(\Fcal))$. For $k\geqslant 0$, consider the mappings $\beta_k:C_k(\Fcal)\to Z^1(\Om(\Fcal))$:
\begin{equation}
\label{(1.29)}
\beta_k(f)=f^{(p^{k+1}-1)}df=df^{(p^{k+1})}.
\end{equation}
Let $\gamma_k:C_k(\Fcal)\to H^1(\Om(\Fcal))$ be the composite of $\beta_k$ and the projection $Z^1(\Om(\Fcal))\to H^1(\Om(\Fcal))$. The mappings $\beta_k$ are, clearly, $G(\Fcal)$-equivariant, and therefore so are the mappings $\gamma_k$. The mappings $\gamma_k$ are, moreover, $p^{k+1}$-semilinear. Indeed, for any $\lambda\in K$ and $f,g\in C_k(\Fcal)$, we have 
\begin{gather}
\notag
\beta_k(\lambda f)=\lambda^{p^{k+1}}\beta_k(f),\\
\label{(1.30)}
\beta_k(f+g)-\beta_k(f)-\beta_k(g)=d\left(\sum_{i=1}^{p^{k+1}-1}f^{(i)}g^{(p^{k+1}-i)}\right)\in B^1(\Om(\Fcal)).
\end{gather}
Further, $\Ker\gamma_k=C_{k+1}(\Fcal)$ since $\beta_k(f)\in B^1(\Om(\Fcal))$ if and only if $f^{(p^{k+1})}\in\Ocal(\Fcal)$.

Let $H^1(\Om(\Fcal))_k$ be the linear span of the image of $\gamma_k$. If $K$ is perfect, then $H^1(\Om(\Fcal))_k=\im\gamma_k$. The $G(\Fcal)$-equivariance of $\gamma_k$ implies the $G(\Fcal)$-invariance of $H^1(\Om(\Fcal))_k$. If $f\in C_k(\Fcal)$, then $f^{(p)}\in C_{k-1}(\Fcal)$ and thanks to \eqref{(1.18)} we have 
$$
\beta_k(f)=df^{(p^{k+1})}=\beta_{k-1}(f^{(p)}),
$$
implying $H^1(\Om(\Fcal))_k\subseteq H^1(\Om(\Fcal))_{k-1}$. 

\begin{lemma}\label{1.10}
\textit{For a basis of $H^1(\Om(\Fcal))_k$ one can take the cohomology classes of cocycles $\{x_i^{(p^{m_i}-1)} dx_i\mid m_i>k\}$. In particular, $H^1(\Om(\Fcal))=H^1(\Om(\Fcal))_0$}.
\end{lemma}

\noindent
\textbf{Proof.}
It follows from the explicit description of the spaces $C_\ell(\Fcal)$ that 
$$
C_k(\Fcal)=C_{k+1}(\Fcal)\oplus\langle x_i^{(p^{m_i-1-k})}\mid m_i>k\rangle_K.
$$
Therefore, $\IM\gamma_k$ is a vector space over $K^{p^{k+1}}$ whose basis consists of the classes of cocycles
\begin{equation}
\label{(1.31)}
\beta_k(x_i^{(p^{m_i-1-k})})=d({(x_i^{(p^{m_i-1-k})})}^{(p^{k+1})})=dx_i^{(p^{m_i})}=x_i^{(p^{m_i}-1)}dx_i\,,\qquad m_i>k.\qed
\end{equation}

\begin{proposition}\label{1.11}
There exist $G(\Fcal)$-equivariant $p^k$-semilinear embeddings 
$$
E_{k-1}/E_k\to H^1(\Om(\Fcal))_{k-1}/H^1(\Om(\Fcal))_k
$$
which are bijective if $K$ is perfect.
\end{proposition}

\noindent
\textbf{Proof.}
Since $E_{k-1}\subseteq C_{k-1}(\Fcal)$ and $E_k\subseteq C_k(\Fcal)$, we have $\gamma_{k-1}(E_{k-1})\subseteq H^1(\Om(\Fcal))_{k-1}$ and $\gamma_{k-1}(E_k)=0$, i.e., $\gamma_{k-1}$ induces a $p^k$-semilinear map
\begin{equation}
\label{(1.32)}
E_{k-1}/E_k\to H^1(\Om(\Fcal))_{k-1}\to H^1(\Om(\Fcal))_{k-1}/H^1(\Om(\Fcal))_k.
\end{equation}
If $x\in E_{k-1}$, then it follows from \eqref{(1.28)} that
$$
\sigma*x-\sigma x\in\fm(\Fcal)^{(2)}\cap C_{k-1}(\Fcal)\subseteq C_k(\Fcal)\oplus\langle x_i^{(p^{m_i-k})}\mid m_i>k\rangle_k,
$$
whence, by \eqref{(1.31)},
\begin{gather*}
\gamma_{k-1}(\sigma*x-\sigma x)\in\langle\gamma_{k-1}(x_i^{(p^{m_i-k})})\mid m_i>k\rangle_K\subseteq H^1(\Om(\Fcal))_k\,;\\
\gamma_{k-1}(\sigma*x)\equiv\gamma_{k-1}(\sigma x)=\sigma\gamma_{k-1}(x)\mod{H^1(\Om(\Fcal))_k}\,.
\end{gather*}
This proves $G(\Fcal)$-equivariance of the mappings \eqref{(1.32)}. Since 
$\gamma_{k-1}$ takes the basis $\{x_i\mid m_i=k\}$ of a complement of $E_k$ 
in $E_{k-1}$ to the basis of a complement of $H^1(\Om(\Fcal))_k$ in 
$H^1(\Om(\Fcal))_{k-1}$ consisting of the classes of cocycles 
$\{x_i^{(p^{m_i}-1)} dx_i\mid m_i=k\}$, the mapping \eqref{(1.32)} 
is injective (and even bijective if the field $K$ is perfect).\qed 

\begin{proposition}\label{1.12}
If $K$ is perfect, then $G(\Fcal)$ acts in $H^1(\Om(\Fcal))$ as the group of all linear transformations preserving all subspaces $H^1(\Om(\Fcal))_k$ for $k\geqslant 0$. Moreover, the group $G(\Fcal)_1$ acts as the group of those transformations which induce the identity transformations in each quotient space $H^1(\Om(\Fcal))_{k-1}/H^1(\Om(\Fcal))_k$.
\end{proposition}

\noindent
\textbf{Proof.}
We already know that the group $G(\Fcal)$ leaves the spaces $H^1(\Om(\Fcal))_k$ invariant, and since $G(\Fcal)_1$ acts on $E$ identically, by the isomorphisms of Proposition \ref{1.11} the group $G(\Fcal)_1$ acts identically on $H^1(\Om(\Fcal))_{k-1}/H^1(\Om(\Fcal))_k$.

Conversely, let $A$ be an invertible endomorphism of the space $H^1(\Om(\Fcal))$ with respect to which the subspaces $H^1(\Om(\Fcal))_k$ are invariant. Let $z_i$ be the class of the cocycle $x_i^{(p^{m_i}-1)}dx_i$, i.e., $z_i=\gamma_{m_i-1}(x_i)$. Let us write $Az_i=\sum\limits a_{ij}z_j$, where $a_{ij}\in K$ and $a_{ij}=0$ when $m_j<m_i$. Set 
$$
y_i=\sum_{m_j\geqslant m_i}b_{ij}x_j^{(p^{m_j-m_i})},
$$
where the $b_{ij}\in K$ are found from the condition $b_{ij}^{p^{m_i}}=a_{ij}$. Since $A$ is invertible, it follows that $\det(a_{ij})_{i,j\in I_k}\ne 0$, where $I_k=\{i\mid m_i=k\}$ for each $k$, whence $\det(b_{ij})_{i,j\in I_k}\ne 0$, and since $y_i\equiv\sum\limits_{m_j=m_i}b_{ij}x_j\mod{\fm(\Fcal)^{(2)}}$, the determinant $\det(\partial_iy_j)_{1\leqslant i,j\leqslant n}$ is also invertible. Moreover, $y_i\in C_{m_i-1}(\Fcal)$. Therefore, there exists $\sigma\in G(\Fcal)$ (Proposition~\ref{1.6}) such that $\sigma x_i=y_i$ for all $i$. Properties of the mappings $\gamma_k$ imply that 
$$
\sigma z_i=\sigma\gamma_{m_i-1}(x_i)=\gamma_{m_i-1}(\sigma x_i)=\gamma_{m_i-1}\left(\sum\limits b_{ij}x_j^{(p^{m_j-m_i})}\right)=\sum b_{ij}^{p^{m_i}}z_j=Az_i.
$$
Finally, if $A$ acts identically on the quotients $H^1(\Om(\Fcal))_{k-1}/H^1(\Om(\Fcal))_k$, then $a_{ij}=\delta_{ij}$ whenever $m_j=m_i$, and the same holds for the $b_{ij}$; hence, $y_i-x_i\in\fm(\Fcal)^{(2)}$, i.e., $\sigma\in G(\Fcal)_1$.\qed 

\begin{proposition}\label{1.13}
The group $G'(\Fcal)$ acts trivially on $H(\Om(\Fcal))$.
\end{proposition}

\noindent
\textbf{Proof.}
Since $G(\Fcal)$ acts on $H(\Om(\Fcal))$ by automorphisms, it suffices to prove that $G'(\Fcal)$ fixes\linebreak
$H^1(\Om(\Fcal))$ elementwise. From the properties of the mappings $\gamma_k$ it follows that
$$
\sigma\gamma_k(f)=\gamma_k(\sigma f)=\gamma_k(f)+\gamma_k(\sigma f-f)=\gamma_k(f)\quad\text{for}\ \ \sigma\in G'(\Fcal),\ \,f\in C_k(\Fcal)
$$
since
$$
\sigma f-f\in\fm(\Fcal)^2\subseteq C_{k+1}(\Fcal)=\Ker\gamma_k.
$$
It remains to observe that
$$
H^1(\Om(\Fcal))=H^1(\Om(\Fcal))_0=K\cdot\IM\gamma_0.\eqno\qed
$$

\medskip
Define the following groups:
\begin{gather*}
U(\Fcal):=\{u\in\widehat{\Ocal}(E)^\times\mid u^{-1} du\in\Om^1(\Fcal)\},\\
\Tilde B^1(\Fcal):=\{u^{-1} du\mid u\in\Ocal(\Fcal)^\times\}\subseteq Z^1(\Om(\Fcal)),\\
\Tilde H^1(\Fcal):=Z^1(\Om(\Fcal))/\Tilde B^1(\Fcal).
\end{gather*}

\begin{proposition}\label{1.14}
The mapping $u\mapsto u^{-1} du$ induces a $G(\Fcal)$-equivariant isomorphism of the multiplicative group $U(\Fcal)/\Ocal(\Fcal)^\times$ onto the additive group $\Tilde H^1(\Fcal)$.
\end{proposition}

\noindent
\textbf{Proof.}
The mapping $u\mapsto u^{-1}du$ is obviously $G(\Fcal)$-equivariant, and by Proposition~\ref{1.3} it maps $U(\Fcal)$ surjectively onto $Z^1(\Om(\Fcal))\subseteq Z^1(\widehat\Om(E))$. The pre-image of $\Tilde B^1(\Fcal)$ under this mapping coincides with $\Ocal(\Fcal)^\times$.\qed 

\medskip
By Proposition~\ref{1.4} and formula \eqref{(1.22)} we have 
\begin{equation}
\label{(1.33)}
\Dfrak:=\{df\mid f\in\fm(\Fcal)^2\}=\{u^{-1}du\mid u\in 1+\fm(\Fcal)^2\}\subseteq B^1(\Om(\Fcal))\cap\Tilde B^1(\Fcal).
\end{equation}
For $f\in\fm(\Fcal)$ we have
\begin{gather}
\label{(1.34)}
\begin{aligned}
(1+f)^{-1}d(1+f)=\sum_{i=0}^{p-1}(-1)^if^idf=df+f^{p-1}df+d\left(\sum_{i=1}^{p-2}\frac{(-1)^i}{i+1}\,f^{i+1}\right)&\\
\equiv df-f^{(p-1)}df & \mod{\Dfrak};
\end{aligned}\\
\label{(1.35)}
df\equiv f^{(p-1)}df\mod{\Tilde B^1(\Fcal)}.
\end{gather}

\begin{proposition}\label{1.15}We have
$$
Z^1(\Om(\Fcal))=\Tilde B^1(\Fcal)\oplus\left\{\sum c_ix_i^{(p^{m_i}-1)}dx_i\mid c_i\in K\right\}.
$$
If $K$ is perfect, then $Z^1(\Om(\Fcal))=\Tilde B{}^1(\Fcal)\oplus dE$.
\end{proposition}

\noindent
\textbf{Proof.}
Since every element of $1+\fm(\Fcal)$ is uniquely expressed  in the form $u\prod\limits_{k<m_i}(1+a_{ik}x_i^{(p^k)})$, where $a_{ik}\in K$, $u\in 1+\fm(\Fcal)^2$, and since the mapping $v\mapsto v^{-1}dv$ isomorphically maps $1+\fm(\Fcal)$ onto $\Tilde B^1(\Fcal)$, we see thanks to \eqref{(1.34)} that
$$
\Tilde B^1(\Fcal)=\Dfrak\oplus\left\{\sum_{k<m_i}(a_{ik}dx_i^{(p^k)}-a_{ik}^p dx_i^{(p^{k+1})})\mid a_{ik}\in K\right\}.
$$
The proof is completed by comparing with the equality
$$
Z^1(\Om(\Fcal))=\Dfrak\oplus\langle dx_i^{(p^\ell)}\mid\ell\leqslant m_i\rangle .
\eqno\qed
$$

\begin{corollary}\label{1.16}
Every element of $U(\Fcal)$ is uniquely represented in the form\/ $\exp\bigl(\sum c_ix_i^{(p^{m_i})}\bigr)\cdot f$, where $c_i\in K,$ $f\in\Ocal(\Fcal)^\times$; and if $K$ is perfect, then also in the form $\exp(x)\cdot f$, where $x\in E$ and $f\in\Ocal(\Fcal)^\times$.
\end{corollary}

Note that $u\,\Om(\Fcal)$ is a subcomplex of $\widehat\Om(E)$ for any $u\in U(\Fcal)$ since for $\om\in\Om(\Fcal)$ we have
$$
d(u\om)=u(d\om+u^{-1}du\wedge\om)\in u\,\Om(\Fcal).
$$

\begin{proposition}\label{1.17}
For $u\in U(\Fcal)\diagdown\Ocal(\Fcal)$ the complex $u\,\Om(\Fcal)$ is acyclic.
\end{proposition}

\noindent
\textbf{Proof.} By Corollary~\ref{1.16} we may assume that $u=\exp\bigl(\sum\limits c_ix_i^{(p^{m_i})}\bigr)$. 
Multiplication by $u$ defines an isomorphism of spaces $\Om(\Fcal)\simeq u\,\Om(\Fcal)$. Let us transfer the differential of the complex $u\,\Om(\Fcal)$ by means of this isomorphism. The new differential $d'$ on $\Om(\Fcal)$ is  as follows:
$$
d'\om=d\om+\sum c_ix_i^{(p^{m_i}-1)}dx_i\wedge\om.
$$
Direct calculation of the cohomology of the complex $(\Om(\Fcal),d')$ has been performed in~\cite{kry}.\qed

\begin{proposition}\label{1.18}
The group $G'(\Fcal)$ acts trivially in $\Tilde H^1(\Fcal),$ and therefore also in $U(\Fcal)/\Ocal(\Fcal)^\times\!$.
\end{proposition}

\noindent
\textbf{Proof.} If $\sigma\in G'(\Fcal)$, then $\sigma x_i-x_i\in\fm(\Fcal)^2$ and 
$$
\sigma x_i^{(p^{m_i})}-x_i^{(p^{m_i})}=\sum_{s=1}^{p^{m_i}}(\sigma x_i-x_i)^{(s)}x_i^{(p^{m_i}-s)}\in\fm(\Fcal)^2.
$$
By \eqref{(1.33)} $\,\sigma\varphi-\varphi\in\Dfrak\subset\Tilde B^1(\Fcal)\,$ for each cocycle $\varphi=\sum c_i\,dx_i^{(p^{m_i})}$ with $c_i\in K$. Therefore, by Proposition 1.15, $\,\sigma\varphi\equiv\varphi\bmod{\Tilde B^1(\Fcal)}\,$ for all $\varphi\in Z^1(\Om(\Fcal))$.\qed

\begin{corollary}\label{1.19}
The group $G'(\Fcal)$ leaves invariant every complex $u\,\Om(\Fcal),$ where $\,u\in U(\Fcal)$.
\end{corollary}

Let us construct a $G(\Fcal)$-invariant decreasing filtration in $\Tilde H^1(\Fcal)$. Denote by
$$
\Tilde d:\Ocal(\Fcal)\to\Tilde H^1(\Fcal)
$$
the composition of $d:\Ocal(\Fcal)\to Z^1(\Ocal(\Fcal))$ and the projection $Z^1(\Om(\Fcal))\to\Tilde H^1(\Fcal)$. Set 
$$
\Tilde H^1(\Fcal)_k=\Tilde dC_k(\Fcal).
$$

\begin{lemma}\label{1.20}If $K$ is perfect, then $\Tilde H^1(\Fcal)_k=\Tilde dE_k$.
\end{lemma}

\noindent
\textbf{Proof.}
Since $d(\fm(\Fcal)^2)=\Dfrak\subseteq\Tilde B^1(\Fcal)$, it follows from Proposition~\ref{1.1} that
$$
\Tilde H^1(\Fcal)_k=\left\{\,\sum a_{i\ell}\,\Tilde dx_i^{(p^\ell)}\mid k+\ell<m_i,\ a_{i\ell}\in K\right\}.
$$
In view of \eqref{(1.35)} we have $a_{i\ell}\,\Tilde dx_i^{(p^\ell)}=b_{i\ell}\,\Tilde dx_i$, where the elements $b_{i\ell}$ are determined from the equation $b_{i\ell}^{p^\ell}=a_{i\ell}$. Thus,
$$
\Tilde H^1(\Fcal)_k=\Bigl\{\,\sum_{m_i>k}c_i\,\Tilde dx_i\mid c_i\in K\Bigr\}=\Tilde dE_k.
\eqno\qed
$$

\begin{proposition}\label{1.21}If $K$ is perfect, then the $G(\Fcal)$-modules $E_{k-1}/E_k$ and $\Tilde H^1(\Fcal)_{k-1}/\Tilde H^1(\Fcal)_k$ are isomorphic.
\end{proposition}

\noindent
\textbf{Proof.} By Lemma~\ref{1.20} and Proposition~\ref{1.15} $\Tilde d$ induces isomorphisms 
$$
E_{k-1}/E_k\cong\Tilde H^1(\Fcal)_{k-1}/\Tilde H^1(\Fcal)_k.
$$
Let us check that they are $G(\Fcal)$-equivariant. If $x\in E_{k-1}$ and $\sigma\in G(\Fcal)$, then
$$
\sigma*x-\sigma x\in\fm(\Fcal)^{(2)}\cap C_{k-1}(\Fcal)=\fm(\Fcal)^2\oplus\langle x_i^{(p^\ell)}\mid\ell\geqslant 1,\ \ell+k-1<m_i\rangle_K.
$$
Since $d(\fm(\Fcal)^2)=\Dfrak\subseteq\Tilde B^1(\Fcal)$ and $d(ax_i^{(p^\ell)})\equiv a^{p^{-\ell}}dx_i\bmod{\Tilde B^1(\Fcal)}$ by \eqref{(1.35)}, it follows that 
\begin{gather*}
\Tilde d(\sigma*x-\sigma x)\in\Tilde dE_k=\Tilde H^1(\Fcal)_k;\\
\Tilde d(\sigma*x)\equiv\Tilde d(\sigma x)\equiv\sigma\Tilde d(x)\mod{\Tilde H^1(\Fcal)_k}.
\lefteqn{\hskip115pt\qed}
\end{gather*}

\begin{lemma}\label{1.22}{\normalfont ($K$ is perfect).}
Two elements of $\Tilde H^1(\Fcal)_{k-1}\diagdown\Tilde H^1(\Fcal)_k$ are conjugate under $G(\Fcal)_1$ if and only if their images in $\Tilde H^1(\Fcal)_{k-1}/\Tilde H^1(\Fcal)_k$ coincide.
\end{lemma}

\noindent
\textbf{Proof.} By Proposition~\ref{1.21} $G(\Fcal)_1$ acts trivially on $\Tilde H^1(\Fcal)_{k-1}/\Tilde H^1(\Fcal)_k$. Therefore, if two elements of $\Tilde H^1(\Fcal)_{k-1}$ are conjugate under $G(\Fcal)_1$, then their images in $\Tilde H^1(\Fcal)_{k-1}/\Tilde H^1(\Fcal)_k$ must coincide. Further, any element of $\Tilde H^1(\Fcal)_{k-1}$ is represented uniquely in the form $\Tilde dx$, where $x\in E_{k-1}$ and $\Tilde dx\notin\Tilde H^1(\Fcal)_k$ if and only if $x\notin E_k$.

Let $x,y\in E_{k-1}\diagdown E_k$. Assume that $\Tilde dx-\Tilde dy\in\Tilde H^1(\Fcal)_k$. Then $x-y\in E_k$. Let us show that $\Tilde dx$ and $\Tilde dy$ are conjugate under $G(\Fcal)_1$. Proposition~\ref{1.6} implies the existence of $\sigma\in G(\Fcal)_1$ such that $\sigma x=x+(y-x)^{(p)}$. Taking into account \eqref{(1.35)} we get
$$ 
\sigma\Tilde dx=\Tilde d\,\sigma x=\Tilde dx+\Tilde d\,(y-x)^{(p)}=\Tilde dx+\Tilde d(y-x)=\Tilde dy.
\eqno\qed 
$$

\medskip
Consider the following Lie subalgebras in $W(\Fcal)$:
\begin{align*}
\gf(\Fcal)&=\left\{\delta\in W(\Fcal)\mid\delta C_k(\Fcal)\subseteq C_k(\Fcal)\quad\text{for all}\ \,k\geqslant 0\right\},\\
\gf'(\Fcal)&=\left\{\delta\in W(\Fcal)\mid\delta\Ocal(\Fcal)\subseteq\fm(\Fcal)^2\right\}\subseteq\gf(\Fcal),\\
\gf(\Fcal)_j&=\gf(\Fcal)\cap W(\Fcal)_j,\\
\gf'(\Fcal)_j&=\gf'(\Fcal)\cap W(\Fcal)_j.
\end{align*}

\begin{lemma}\label{1.23}
\quad$\gf(\Fcal)=\sum C_{m_i-1}(\Fcal)\,\partial_i\,;\quad\gf'(\Fcal)=\fm(\Fcal)^2\,W(\Fcal)$.
\end{lemma}

\noindent
\textbf{Proof.} Let $\delta=\sum f_i\partial_i\in W(\Fcal)$. If $\delta\in\gf(\Fcal)$ (resp. $\delta\in\gf'(\Fcal)$), then $f_i=\delta x_i\in C_{m_i-1}(\Fcal)$ (resp. $f_i\in\fm(\Fcal)^2$) since $x_i\in C_{m_i-1}(\Fcal)$.

Conversely, let $f_i\in C_{m_i-1}(\Fcal)$ for all $i $. In particular, since $f_i\in\fm(\Fcal)$ for all $i$, it follows that $\delta(\fm(\Fcal)^{(2)})\subseteq\fm(\Fcal)^2\subseteq C_\ell(\Fcal)$ for any $\ell $. Since $C_k(\Fcal)\subseteq\langle x_j\mid m_j>k\rangle_K+\fm(\Fcal)^{(2)}$ and $\delta x_j=f_j\in C_k(\Fcal)$ for any $k<m_j$, it follows that $\delta C_k(\Fcal)\subseteq C_k(\Fcal)$, i.e., $\delta\in\gf(\Fcal)$. If, moreover, $f_i\in\fm(\Fcal)^2$ for all $i$, then $\delta\Ocal(\Fcal)\subseteq\fm(\Fcal)^2$ and $\delta\in\gf'(\Fcal)$.\qed 

\begin{proposition}\label{1.24}{\normalfont (cf.~\cite{11})}
For $j\geqslant 1,$ given $\sigma\in G'(\Fcal)_j$ we can select $\delta\in\gf'(\Fcal)_j,$ and the other way round, given $\delta\in\gf'(\Fcal)_j$ we can select $\sigma\in G'(\Fcal)_j$ so that 
\begin{gather}
\notag
\sigma x=x+\delta x\quad\text{\textit{for all}}\quad x\in E\quad\text{and}\\
\label{(1.36)}
(\sigma-\id-\delta)(\fm(\Fcal)^{(\ell)})\subseteq\fm(\Fcal)^{(j+\ell+1)}\quad\text{\textit{for all}}\quad\ell\geqslant 0.
\end{gather}
\end{proposition}

\noindent
\textbf{Proof.} Let $\sigma\in G'(\Fcal)_j$. Then $f_i=\sigma x_i-x_i\in\fm(\Fcal)^2\cap\fm(\Fcal)^{(j+1)}$, whence $\delta=\sum f_i\partial_i\in\gf'(\Fcal)_j$. By construction $\sigma x=x+\delta x$ for every $x\in E$. If $r<m_i$, then
$$
\sigma x_i^{(r)}-x_i^{(r)}-\delta x_i^{(r)}=(x_i+f_i)^{(r)}-x_i^{(r)}-x_i^{(r-1)}f_i=\sum\limits_{s=2}^rx_i^{(r-s)}f_i^{(s)}.
$$
Since 
\[
\,x_i^{(r-s)}f_i^{(s)}\in\fm(\Fcal)^{(r-s)}\fm(\Fcal)^{((j+1)s)}\subseteq\fm(\Fcal)^{(r+js)}\subseteq\fm(\Fcal)^{(r+j+1)}\,\text{~~ when $s\geqslant 2$,}
\]
 it follows that
$$
\sigma x_i^{(r)}-x_i^{(r)}-\delta x_i^{(r)}\in\fm(\Fcal)^{(r+j+1)}.
$$

Denote
$$
A_k=\left\{h\in\fm(\Fcal)^{(k)}\mid\sigma h-h-\delta h\in\fm(\Fcal)^{(j+k+1)}\right\}.
$$
Thus, $x_i^{(r)}\in A_r$. We have $A_kA_\ell\subseteq A_{k+l}$ since
\begin{multline*}
\sigma(gh)-gh-\delta(gh)
=(\sigma g-g)(\sigma h-h)+g(\sigma h-h-\delta h)+(\sigma g-g-\delta g)h\in\\
\in\fm(\Fcal)^{(2j+k+\ell)}+\fm(\Fcal)^{(j+k+\ell+1)}+\fm(\Fcal)^{(j+k+\ell+1)}\subseteq\fm(\Fcal)^{(j+k+\ell+1)}
\end{multline*}
for $g\in A_k$ and $h\in A_l$. Hence,
$$
(\sigma-\id-\delta)(x_1^{(\alpha_1)}x_2^{(\alpha_2)}\ldots x_n^{(\alpha_n)})\subseteq\fm(\Fcal)^{(j+\alpha_1+\alpha_2+\ldots+\alpha_n+1)},
$$
implying \eqref{(1.36)}.

Conversely, if $\delta\in\gf'(\Fcal)_j$, then $\delta x_i\in\fm(\Fcal)^2\cap\fm(\Fcal)^{(j+1)}$, and by Proposition~\ref{1.8} there exists a $\sigma\in G'(\Fcal)_j$ that sends $x_i$ to $x_i+\delta x_i$. Again, $\sigma x=x+\delta x$ for any $x\in E$. Repeating the above arguments we get \eqref{(1.36)}.\qed 

\begin{lemma}\label{1.25}
Let $\,\psi\in u\,\fm(\Fcal)^{(r)}\,\Om^k(\Fcal)$ and $\sigma\in G'(\Fcal)_j$, where $u\in U(\Fcal)$ and $j\geqslant 1$. Then
$$
\sigma\psi-\psi\in u\,\fm(\Fcal)^{(j+r)}\,\Om^k(\Fcal).
$$
Moreover, if $\delta\in\gf'(\Fcal)_j$ is connected with $\sigma$ by relation 
\textrm{\eqref{(1.36)}}, then 
$$
\sigma\psi-\psi-\delta\psi\in u\,\fm(\Fcal)^{(j+r+1)}\,\Om^k(\Fcal).
$$
\end{lemma}

The proof is standard.

\section{The $G'(\Fcal)$-orbits of \repl{symplectic} forms}\label{2}

For each form $\om\in u\,\Om^2(\Fcal)$ define a homomorphism of $\Ocal(\Fcal)$-modules $i_\om:W(\Fcal)\to u\,\Om^1(\Fcal)$ by the rule $i_\om(\delta)=\delta\ugol\om$. A \textit{\repl{symplectic} form \repl{corresponding} to the flag $\Fcal$} is a closed form $\om\in u\,\Om^2(\Fcal)$ for any $u\in U(\Fcal)$ such that $i_\om$ is an isomorphism. Note that if
$$
\om=u\sum_{i<j}h_{ij}dx_i\wedge dx_j
$$
where $h_{ij}\in\Ocal(\Fcal)$, then $i_\om$ is bijective if and only if $\det(h_{ij})$ is invertible in $\Ocal(\Fcal)$.

Denote by $\Ham(\Fcal)$ the set of all \repl{symplectic} forms \repl{corresponding} to the flag $\Fcal$. The forms belonging, respectively, to
$$
\Ham^1(\Fcal)=\Ham(\Fcal)\cap\Om(\Fcal)\quad\text{and}\quad\Ham^2(\Fcal)=\Ham(\Fcal)\diagdown\Ham^1(\Fcal)
$$
will be called \repl{symplectic} forms of the \textit{first} (resp. \textit{second}) type.

The action of $G(\Fcal)$ in $\widehat\Om(E)$ leaves invariant the sets $\Ham(\Fcal)$, $\Ham^1(\Fcal)$, and $\Ham^2(\Fcal)$. The group $G'(\Fcal)$ preserves the sets of \repl{symplectic} forms in each complex $u\,\Om(\Fcal)$, where $u\in U(\Fcal)$ (Corollary~\ref{1.19}).

In what follows we will write briefly $\Ocal=\Ocal(\Fcal)$, $\fm=\fm(\Fcal)$, $U=U(\Fcal)$, $\Om=\Om(\Fcal)$, and so on. Note that if $\om,\om'\in\uOm^2$ are closed forms such that $\om$ is \repl{symplectic} and $\om'-\om\in\umOm^2$, then $\om'$ is also \repl{symplectic}.

\begin{lemma}\label{2.1}
Any form $\psi\in\uOm^1$, where $u\in U$, can be expressed as $df+\Tilde\psi$,  where $f\in\uOcal$ and $\Tilde\psi\in\umOm^1$.
\end{lemma}

\noindent
\textbf{Proof.} Let us write $\psi=u(dx+\varphi)$,  where $x\in E$ and $\varphi\in\mOm^1$. We have
$$
\psi=d(ux)+u(\varphi-xu^{-1}du).\eqno\qed 
$$

\begin{lemma}\label{2.2}
Any form $\psi\in\ufm^{(j+1)}\Om^1$, where $u\in U,$ $j\geqslant 1$, can be expressed as $df+\Tilde\psi$,  where $f\in u(\fm^2\cap\fm^{(j+2)})$ and $\Tilde\psi\in u(\fm^2\cap\fm^{(j+1)})\,\Om^1$.
\end{lemma}

\noindent
\textbf{Proof.} Let $\psi=u\sum f_i dx_i$, where $f_i\in\fm^{(j+1)}$. Then
$$
\psi=d\Bigl(u\sum f_ix_i\Bigr)-u\Bigl(\,\sum x_i\,df_i+\sum x_if_iu^{-1}du\Bigr).
$$
Obviously, $x_if_i\in\fm^2\cap\fm^{(j+2)}$ and, since $df_i\in\fm^{(j)}\Om^1$, we have $x_i\,df_i\in(\fm^2\cap\fm^{(j+1)})\,\Om^1$.\qed 

\begin{lemma}\label{2.3}
Let $\om\in\uOm^k$, where $u\in U$, be a closed form. Let $\sigma\in G'_j$ and $\delta\in\gf'_j$, where $j\geqslant 1$, be related by condition \eqref{(1.36)}. Then there exists $\varphi\in\uOm^{k-1}$ such that
$$
\sigma\om-\om=d\varphi\quad\textit{and}\quad\varphi-\delta\ugol\om\in\ufm^{(j+2)}\Om^{k-1}.
$$
\end{lemma}

\noindent
\textbf{Proof.}
Suppose first that $\om$ is an exact form, i.e., $\om=d\psi$, where 
$\psi\in\uOm^{k-1}$. Adjusting $\psi$ by a coboundary we may assume 
thanks to Lemma~\ref{2.1} that $\psi\in\umOm^{k-1}$. Setting 
$$
\varphi=\sigma\psi-\psi-d(\delta\ugol\psi),
$$
we get $\sigma\om-\om=d(\sigma\psi-\psi)=d\varphi$. By \eqref{(1.6)} $\delta\ugol\om=\delta\ugol d\psi=\delta\psi-d(\delta\ugol\psi)$, whence, in view of Lemma~\ref{1.25},
$$
\varphi-\delta\ugol\om=\sigma\psi-\psi-\delta\psi\in\ufm^{(j+2)}\Om^{k-1}.
$$

If $u\notin\Ocal$, then all closed forms are exact (by Proposition~\ref{1.17}), and so the lemma is proved in this case. 

Let $u=1$. For every integer $\ell\geqslant 0$, denote by $A^\ell$ the subspace of all forms $\om\in Z^\ell(\Om)$ for which the lemma holds. Set $A=\bigoplus\limits_{\ell=0}^nA^\ell $. We have proved already that $B(\Om)\subseteq A$. We claim that $A$ is a subalgebra of $\Om $. Indeed, let $\om\in A^k$ and $\om'\in A^\ell$, so that
$$ 
\sigma\om-\om=d\varphi,\quad\sigma\om'-\om'=d\varphi'\quad\text{for some}\ \,\varphi\in\Om^{k-1},\ \,\varphi'\in\Om^{\ell-1}
$$
such that
\begin{equation}
\label{(2.1)}
\varphi-\delta\ugol\om\in\fm^{(j+2)}\Om^{k-1}\quad\text{and}\quad\varphi'-\delta\ugol\om'\in\fm^{(j+2)}\Om^{\ell-1}.
\end{equation}
Then 
$$
\sigma(\om\wedge\om')-\om\wedge\om'=(\sigma\om-\om)\wedge\sigma\om'+\om\wedge (\sigma\om'-\om')=d(\varphi\wedge\sigma\om'+(-1)^k\om\wedge\varphi'),
$$
and also
\begin{multline}
\label{(2.2)}
\varphi\wedge\sigma\om'+(-1)^k\om\wedge\varphi'-\delta\ugol(\om\wedge\om')=\\
=(\varphi-\delta\ugol\om)\wedge\sigma\om'+(-1)^k\om\wedge(\varphi'-\delta\ugol\om')+(\delta\ugol\om)\wedge(\sigma\om'-\om').
\end{multline}
Since $\delta\in\gf_j'\subseteq\fm^{(j+1)}W\!$, we have $\delta\ugol\om\in\fm^{(j+1)}\,\Om^{k-1}\!$. By Lemma~\ref{1.25} $\sigma\om'-\om'\in\fm^{(j)}\,\Om^\ell\!$. Making use of \eqref{(2.1)} as well, we see that the element given by the expression \eqref{(2.2)} lies in $\fm^{(j+2)}\,\Om^{k+\ell-1}$. Thus, $\,\om\wedge\om'\in A^{k+\ell}$.

Let us show that $Z^1(\Om)\subseteq A$. For $f\in\fm$, we have
$$
\sigma\,df^{(p)}-df^{(p)}=d\bigl((\sigma f)^{(p)}-f^{(p)}\bigr)=dh,
$$
where $h=\sum\limits_{r=0}^{p-1}f^{(r)}(\sigma f-f)^{(p-r)}\in\Ocal(\Fcal)$ since $\sigma f-f\in\fm^2$. Moreover,
\begin{multline}
\label{(2.3)}
\hskip1cm
h-\delta\ugol df^{(p)}=h-\delta f^{(p)}=h-f^{(p-1)}\delta f=\\
=f^{(p-1)}(\sigma f-f-\delta f)+\sum_{r=0}^{p-2}f^{(r)}(\sigma f-f)^{(p-r)}.
\end{multline}
Thanks to \eqref{(1.36)} we see that 
$$
f^{(p-1)}(\sigma f-f-\delta f)\in\fm^{(p+j+1)}\quad\text{and}\quad 
f^{(r)}(\sigma f-f)^{(p-r)}\in\fm^{{\textstyle(}r+(p-r)(j+1){\textstyle)}}.
$$
So expression \eqref{(2.3)} is definitely in $\fm^{(j+2)}$, i.e., $df^{(p)}\in A$.

But every 1-cocycle is congruent modulo the space of coboundaries to an appropriate cocycle $df^{(p)}$. This proves the inclusion $Z^1(\Om(\Fcal))\subseteq A$.

It follows from Proposition~\ref{1.9} that the algebra $Z(\Om)$ is generated by $B(\Om)$ and $Z^1(\Om)$. It is clear now that $A=Z(\Om)$.\qed 

\begin{corollary}\label{2.4}
If $\om\in\uOm^k$ is a closed form, and if $\sigma\in G'_j$,  where $j\geqslant 1$, then $\sigma\om-\om=d\varphi$ for a suitable form $\varphi\in\ufm^{(j+1)}\,\Om^{k-1}$.
\end{corollary}

\begin{lemma}\label{2.5}
Let $\om,\om'\in\uOm^2$,  where $u\in U$, be two \repl{symplectic} forms such that $\om'-\om=d\varphi$,  where $\varphi\in\ufm^{(j+1)}\,\Om^1$ and $j\geqslant 1$. Then there exists $\sigma\in G_j'$ such that $\om'-\sigma\om=d\Tilde\varphi$ for some $\Tilde\varphi\in\ufm^{(j+2)}\,\Om^1$.
\end{lemma}

\noindent
\textbf{Proof.}
Adjusting $\varphi$ by a coboundary we may assume thanks to Lemma~\ref{2.2} that 
$$
\varphi\in u(\fm^2\cap\fm^{(j+1)})\,\Om^1.
$$
Since $i_\om$ is an isomorphism of $\Ocal$-modules, we can find $\delta\in(\fm^2\cap\fm^{(j+1)})W=\gf_j'$ such that $\delta\ugol\om=\varphi$. Proposition~\ref{1.24} ensures the existence of a $\sigma\in G_j'$ such that \eqref{(1.36)} holds. By Lemma~\ref{2.3} $\,\sigma\om-\om=d\psi\,$ for some $\psi\in\uOm^1$ such that 
$$
\psi-\delta\ugol\om=\psi-\varphi\in\ufm^{(j+2)}\,\Om^1.
$$
Besides, $\,\om'-\sigma\om=\om'-\om-(\sigma\om-\om)=d(\varphi-\psi)$.\qed 

\begin{proposition}\footnote{It should be added to the hypothesis that here $j\ge1$. In fact, $G'_0=G'=G'_1$.}
\label{2.6}
The $G'_j$-orbit of a \repl{symplectic} form $\om\in\uOm^2$,  where $u\in U$, consists of the forms $\om+d\varphi$ such that $\varphi\in\ufm^{(j+1)}\,\Om^1$.
\end{proposition}

\noindent
\textbf{Proof.}
By Corollary~\ref{2.4} the orbit $G'_j\,\om$ is contained in the indicated set of forms. Conversely, let us prove that for any \repl{symplectic} form $\om\in\uOm^2$ and an integer $k\geqslant 1$ we have 
\begin{equation}
\label{(2.4)}
\{\om+d\varphi\mid\varphi\in\ufm^{(k+1)}\,\Om^1\}\subseteq G'_k\,\om 
\end{equation}
by the downward induction on $k$. For $k$ sufficiently large we have $\fm^{(k+1)}=0$ and the set in the left-hand side of inclusion \eqref{(2.4)} consists of a single form $\om$, so that everything is obvious.

Suppose that \eqref{(2.4)} has been proved for $k=j+1$,  where $j\geqslant 1$. Let $\varphi\in\ufm^{(j+1)}\,\Om^1$. By Lemma~\ref{2.5} applied to $\om$ and $\om'=\om+d\varphi$ we can find $\sigma\in G'_j$ and $\Tilde\varphi\in\ufm^{(j+2)}\,\Om^1$ such that $\om'-\sigma\om=d\Tilde\varphi$. Hence, $\,\om'\in G_{j+1}'\sigma\om\subseteq G_j'\om\,$ by the induction hypothesis.\qed 

\begin{proposition}\label{2.7}
For $p>2$ \repl{symplectic} forms $\om,\om'\in\uOm^2$,  where $u\in U$, are conjugate under $G'$ if and only if their images in $\uOm^2/\umOm^2$ and in $H^2(\uOm)$ coincide.
\end{proposition} 

\noindent
\textbf{Proof.}
Since the mappings $\uOm^2\to\uOm^2/\umOm^2$ and $Z^2(\uOm)\to H^2(\uOm)$ are $G'$-equivariant and $G'$ acts trivially in the space $\uOm^2/\umOm^2$ because $G'\subseteq G_1$, and in the space $H^2(\uOm)$ by Propositions~\ref{1.13} and~\ref{1.17}, it follows that the images of conjugate forms in these spaces coincide. 

Conversely, let $\om$ and $\om'$ be cohomologous and their \textit{initial} terms coincide. In particular, $\om'-\om=d\psi$,  where by Lemma~\ref{2.1} we may assume that $\psi\in\umOm^1$. Let $\psi=u\sum b_{ij}x_i dx_j+\Tilde\psi$, where $b_{ij}\in K$ and $\Tilde\psi\in\ufm^{(2)}\Om^1$. Then 
$$
d\psi\equiv u\sum b_{ij}dx_i\wedge dx_j\pmod{\umOm^2}.
$$
Since the images of $\om$ and $\om'$ in $\uOm^2/\umOm^2$ coincide, we have $b_{ij}=b_{ji}$. Set
$$
\varphi:=\psi-d\Bigl(u\sum_{i<j}b_{ij}x_ix_j+u\sum b_{ii}x_i^{(2)}\Bigr)=\Tilde\psi-u\Bigl(\sum_{i<j}b_{ij}x_ix_j+\sum b_{ii}x_i^{(2)}\Bigr)u^{-1}du.
$$
Then $\om'-\om=d\varphi$, but now $\varphi\in\ufm^{(2)}\Om^1$. It remains to use Proposition~\ref{2.6} with $j=1$.\qed 

\begin{corollary}\label{2.8}{\normalfont($p>2$).}
Any \repl{symplectic} form $\om\in\uOm^2$ is conjugate under $G'$ to a unique form of the following expression, where $a_{ij}$, $b_{ij}\in K$\/{\rm:}
$$
\begin{array}{cl}
\displaystyle\sum_{i<j}\,\bigl(a_{ij}+b_{ij}x_i^{(p^{m_i}-1)}x_j^{(p^{m_j}-1)}\bigr)\,dx_i\wedge dx_j & \text{when}\ \,u=1,\quad\text{or}\\
\displaystyle d\Bigl(u\sum_{i<j}a_{ij}x_idx_j\Bigr) & \text{when}\ \,u\in U\diagdown\Ocal.
\end{array}
$$
\end{corollary}

We will use the results just obtained for a preliminary description of the orbits of \repl{symplectic} forms with respect to the whole group $G$.

Let $K$ be perfect of characteristic $p>2$. Since $G'$ is normal in $G$, the group $G$ acts on the set $\Ham(\Fcal)/G'$, and there is a bijection
$$
\Ham(\Fcal)/G\cong\bigl(\Ham(\Fcal)/G'\bigr)/G.
$$
Denote by $\LaEnd$ the subset of bivectors in $\LaE$ corresponding to the non-degenerate bilinear forms on the dual space $E^*$. There is a mapping
$$
\Ham(\Fcal)\subseteq\widehat\Om^2(E)\to\widehat\Om^2(E)/\,\widehat\fm(E)\widehat\Om^2(E)\cong\Om^2/\,\mOm^2\cong{\textstyle\bigwedge^2}\bigl(\Om^1/\,\mOm^1\bigr)\cong\LaE
$$
which is $G$-equivariant and whose image coincides with $\LaEnd$. Since $G'$ acts trivially in $E$, we get a mapping
$$
\Ham(\Fcal)/G'\to\LaEnd\,.
$$
The projection $Z^2(\Om)\to H^2(\Om)$ together with the isomorphism $H^2(\Om)\cong\bigwedge^2\bigl(H^1(\Om)\bigr)$ (Proposition~\ref{1.9}) give rise to a $G$-equivariant mapping
$$
\Ham^1(\Fcal)\subseteq Z^2(\Om)\to{\textstyle\bigwedge^2}H^1(\Om).
$$
Since $G'$ acts trivially on $H^1(\Om)$ (Proposition~\ref{1.13}), we get a mapping
$$
\Ham^1(\Fcal)/G'\to{\textstyle\bigwedge^2}H^1(\Om).
$$
Consider the $G$-equivariant mapping
\begin{equation}\label{(2.5)}
\Ham^1(\Fcal)/G'\to\LaEnd\times{\textstyle\bigwedge^2}H^1(\Om).
\end{equation}
Its injectivity is equivalent to the conclusion of Proposition~\ref{2.7}, while surjectivity follows from consideration of the forms explicitly written out in Corollary~\ref{2.8}. Thus, the mapping \eqref{(2.5)} is a bijection.

Recall (see \ref{1}) that there is a $G$-invariant flag of subspaces in $H^1(\Om)$; we will denote it by $\Hcal$. Let $\St_\Fcal$ (resp. $\St_\Hcal$) be the subgroup of linear transformations of the space $E$ (resp. $H^1(\Om)$) preserving the flag $\Fcal$ (resp. $\Hcal$). The image $\overline G$ of $G$ in $GL(E)\times GL(H^1(\Om))$ is contained in the subgroup $T$ of all pairs $\,(A,B)\in\St_\Fcal\times\St_\Hcal\,$ rendering commutative the diagrams
\begin{equation}\label{(2.6)} 
\begin{CD}
E_{k-1}/E_k@>>>H^1(\Om)_{k-1}/H^1(\Om)_k\\
@VAVV@VVBV\\
E_{k-1}/E_k@>>>H^1(\Om)_{k-1}/H^1(\Om)_k
\end{CD}
\end{equation}
where horizontal arrows are the isomorphisms of Proposition~\ref{1.11}, while vertical arrows are induced by $A$ and $B$, respectively. Let $\rho:T\to\St_\Fcal$ be the projection onto the first factor. Then $\Ker\rho$ consists of all pairs $(\id_E,B)$,  where $B\in\St_\Hcal$ induces the identity transformation in each factor of the flag $\Hcal$. By Proposition~\ref{1.12} $\Ker\rho\subseteq\overline G$. Also, $\rho(\overline G)=\St_\Fcal=\rho(T)$, whence $\overline G=T$. As a result we get a bijection 
\begin{equation}
\label{(2.7)}
\Ham^1(\Fcal)/G\cong\bigl(\Ham^1(\Fcal)/G'\bigr)/G\cong\bigl(\LaEnd\times{\textstyle\bigwedge^2}H^1(\Om)\bigr)/T.
\end{equation}

If $\om\in\Ham(\Fcal)$ is such that $\om\in\uOm^2\cap v\Om^2$,  where $u,v\in U$, then $u^{-1}v\in\Ocal^\times$. Therefore, to each \repl{symplectic} form there corresponds a well-defined class in $U\!/\Ocal^\times$. The resulting mapping $\Ham(\Fcal)\to U\!/O^\times$ is $G$-equivariant. Composing it with the isomorphism of Proposition~\ref{1.14} we get a mapping $\Ham(\Fcal)\to\Tilde H{}^1$. Under this mapping non-zero elements of the group $\Tilde H^1$ correspond precisely to the \repl{symplectic} forms of the second type. Since $G'$ acts trivially on $\Tilde H^1$ (Proposition~\ref{1.18}), we obtain a map
$$
\Ham^2(\Fcal)/G'\to\Tilde H^1\diagdown\{0\}.
$$
Consider the $G$-equivariant map
\begin{equation}\label{(2.8)}
\Ham^2(\Fcal)/G'\to\LaEnd\times\bigl(\Tilde H^1\diagdown\{0\}\bigr).
\end{equation}
Its injectivity follows from Proposition~\ref{2.7} since $H^2(\uOm)=0$ when $u\in U\diagdown\Ocal$, while surjectivity follows from Corollary~\ref{2.8}. Thus,  \eqref{(2.8)} is a bijection. Since $G_1$ acts trivially in $E$, we derive from \eqref{(2.8)} a $G$-equivariant bijection
\begin{equation}
\label{(2.9)}
\Ham^2(\Fcal)/G_1\to\LaEnd\times\bigl((\Tilde H^1\diagdown \{0\})/G_1\bigr).
\end{equation}
Let us bring the filtration in the space $\Tilde H^1$ (see \ref{1}) into play. By Lemma~\ref{1.22} the projections $\Tilde H^1_{k-1}\to\Tilde H^1_{k-1}/\Tilde H{}_k^1$ induce bijections
$$
(\Tilde H^1_{k-1}\diagdown\Tilde H^1_k)/G_1\to (\Tilde H{}_{k-1}^1/\Tilde H{}_k^1)\diagdown \{0\}.
$$
Proposition~\ref{1.21} allows us to re-express the target set in the form $(E_{k-1}/E_k)\diagdown \{0\}$. Now the set $(\Tilde H^1\diagdown \{0\})/G_1$ is represented as a disjoint union
$$
\coprod_{k\geqslant 1}\,\bigl((E_{k-1}/E_k)\diagdown\{0\}\bigr).
$$
Substituting this into \eqref{(2.9)} we get a $G$-equivariant bijection
\begin{equation}
\label{(2.10)}
\Ham^2(\Fcal)/G_1\to\coprod_{k\geqslant 1}\,\Bigl(\LaEnd\times\bigl((E_{k-1}/E_k)\diagdown\{0\}\bigr)\Bigr).
\end{equation}
Since the image of $G$ in $GL(E)$ coincides with the group $\St_\Fcal$, it follows that \eqref{(2.10)} yields a bijection
\begin{equation}
\label{(2.11)}
\Ham^2(\Fcal)/G\to\coprod_{k\geqslant 1}\,\Bigl(\LaEnd\times\bigl((E_{k-1}/E_k)\diagdown\{0\}\bigr)\Bigr)/\St_\Fcal.
\end{equation}
\textbf{Remarks.}
1) The condition $u\in U$ in the definition of a \repl{symplectic} form was proposed in a somewhat more general form in~\cite{KfiD}.

2) The results of this section are valid also for $p=2$ provided that it is additionally\footnote{This additional condition concerns only the symplectic forms of the 2nd type.} assumed that all $m_i>1$ whenever $u\notin\Ocal(\Fcal)$.

\section{The $G'(\Fcal)$-orbits of contact forms}\label{3}

For every form $\om\in\Om^1(\Fcal)$, set
$$
P_\om=\{\delta\in W(\Fcal)\mid\delta\ugol d\om=0\},\qquad 
Q_\om=\{\delta\in W(\Fcal)\mid\om(\delta)=0\}.
$$
A form $\om\in\Om^1(\Fcal)$ will be called a \textit{contact form \repl{corresponding} to the flag} $\Fcal\,$ if $\,W(\Fcal)=P_\om\oplus Q_\om$. The set of all contact forms \repl{corresponding} to $\Fcal$ will be denoted by $\Cont (\Fcal)$. Set
$$ 
P_\om^\perp=\{\psi\in\Om^1(\Fcal)\mid\psi(P_\om)=0\},\qquad Q_\om^\perp=\{\psi\in\Om^1(\Fcal)\mid\psi(Q_\om)=0\}.
$$

Recall that $\,i_{d\om}(\delta)=\delta\ugol d\om$.

\begin{proposition}\label{3.1}
Let $\om\in\Cont(\Fcal)$. Then\/{\rm:} \ $P_\om$ and $Q_\om$ are free $\Ocal$-modules{\rm;}

\smallskip
$\om$ isomorphically maps $P_\om$ onto $\Ocal,\,$ while $i_{d\om}$ isomorphically maps $Q_\om$ onto $P_\om^\perp ${\rm;}

\smallskip
$Q_\om^\perp=\Ocal\om,$\qquad$\Om^1(\Fcal)=\Ocal\om\oplus P_\om^\perp,$\qquad$f\om\in\Cont(\Fcal)\,$ for every $f\in\Ocal^\times$.
\end{proposition}

\noindent
\textbf{Proof.}
From the decomposition $W=P_\om\oplus Q_\om$ it follows that $P_\om$ and $Q_\om$ are projective finitely generated modules over a local ring $\Ocal$. Hence, they are free. Since $\Ker\om=Q_\om$, the mapping $\om|_{P_\om}$ is injective, and, since $P_\om$ is a free $\Ocal$-module, comparison of dimensions yields $\om(P_\om)=\Ocal$. Thus,  $\om|_{P_\om}$ is an isomorphism, implying $Q_\om^\perp=\Ocal\om$. It follows now from the decomposition $W=P_\om\oplus Q_\om$ that
$$
\Om^1=\Hom_{\Ocal}(W,\Ocal)=P_\om^\perp\oplus Q_\om^\perp=P_\om^\perp\oplus\Ocal\om .
$$
Since $\Ker i_{d\om}=P_\om$, the mapping $i_{d\om}\mskip1mu|_{Q_\om}$ is injective and $i_{d\om}(Q_\om)\subseteq P_\om^\perp$. Since $\dim\Om^1=\dim W$, we have $\dim Q_\om=\dim P_\om^\perp$, whence $i_{d\om}(Q_\om)=P_\om^\perp$.

Obviously, $Q_{f\om}=Q_\om$. Since $d(f\om)=df\wedge\om+fd\om$, the restrictions of the forms $d(f\om)$ and $fd\om$ to $Q_\om\times Q_\om$ coincide. Therefore, $i_{d(f\om)}$ induces an isomorphism of $Q_\om$ onto $\Hom_{\Ocal}(Q_\om,\Ocal)$, implying $W=P_{f\om}\oplus Q_\om$. Thus,  $f\om\in\Cont(\Fcal)$.\qed 

The semidirect product $\OFrtmGF$ acts on $\Cont(\Fcal)$. Note that in coordinates the condition on a form $\sum f_idx_i$ to be contact means precisely that the determinant
$$
\det\left(\begin{array}{c|c}
\partial_if_j-\partial_jf_i&f_i\\
\hline
-f_j&0
\end{array}
\right)
$$
is invertible.
If $\om,\om'\in\Om^1$ are such that $\om\in\Cont(\Fcal)$ and $\om'-\om\in\fm^{(2)}\Om^1$, then $d\om'-d\om\in\mOm^2$ and, as is easy to see, $\om'\in\Cont(\Fcal)$. 

\begin{lemma}\label{3.2}
Let $\om,\om'\in\Cont(\Fcal)$ and $\om'-\om=df$,  where $f\in\fm^2\cap\fm^{(j+1)}$ with $j\geqslant 1$. There exists $\sigma\in G_j'$ such that $\om'-\sigma\om\in\fm^{(j+1)}\Om^1$.
\end{lemma}

\noindent
\textbf{Proof.}
Since $\om|_{P_\om}$ is an isomorphism, there is a unique 
$\partial\in P_\om$ such that $\om(\partial)=1$. 
Find some $\sigma\in G_j'$ such that \eqref{(1.36)} is satisfied with 
$\delta=f\partial\in\fg_j'$. Since $\delta\ugol d\om=0$ and 
$\om(\delta)=f$, we have
$$
\delta\om=\delta\ugol d\om+d(\delta\ugol\om)=df,
$$
whence (by Lemma 1.25 )
\[
\quad\om'-\sigma\om=df+\om-\sigma\om=\delta\om+\om-\sigma\om\in\fm^{(j+1)}\Om^1.\qed
\]

\begin{lemma}\label{3.3}
Let $\om,\om'\in\Cont(\Fcal)$ be such that $\,\om'-\om\in\fm^{(j+1)}\,\Om^1\,$ for some $j\geqslant 1$. Then there exist $f\in 1+\fm^{(j+1)}$ and $\sigma\in G'_j$ such that $\,f\sigma\om'-\om\in\fm^{(j+2)}\,\Om^1$.
\end{lemma}

\noindent
\textbf{Proof.} Using Lemma~\ref{2.2} write $\,\om'-\om=dt+\chi\,$,  where
$$
t\in\fm^2\cap\fm^{(j+2)}\quad\text{and}\quad\chi\in(\fm^2\cap\fm^{(j+1)})\,\Om^1.
$$
By Lemma~\ref{3.2} applied to the forms $\om'$ and $\om'-dt$ we have $\om'-dt-\tau\om'\in\fm^{(j+2)}\,\Om^1$ for some $\tau\in G_{j+1}'$. Hence
\begin{equation}
\tau\om'-\om-\chi\in\fm^{(j+2)}\Om^1.
\label{(3.1)}
\end{equation}

Let $\chi=h\om+\psi$, where $h\in\Ocal$ and $\psi\in P_\om^\perp$ (Proposition~\ref{3.1}). Since the projections ${\Om^1\to P_\om^\perp}$ and ${\Om^1\to Q_\om^\perp}$ are $\Ocal$-module homomorphisms, we have
$$
h\in\fm^2\cap\fm^{(j+1)}\quad\text{and}\quad\psi\in(\fm^2\cap\fm^{(j+1)})P_\om^\perp.
$$
Setting $\,g=(1+h)^{-1}\in 1+\fm^{(j+1)}\,$ and $\,\varphi=(1+h)^{-1}\psi\in(\fm^2\cap\fm^{(j+1)})P_\om^\perp$, we get by \eqref{(3.1)}
\begin{equation}
g\tau\om'-\om-\varphi\in\fm^{(j+2)}\Om^1.
\label{(3.2)}
\end{equation}

Since $i_{d\om}:Q_\om\to P_\om^\perp$ is an isomorphism of $\Ocal$-modules, there is a unique derivation
$$
\delta\in(\fm^2\cap\fm^{(j+1)})\,Q_\om\subseteq\fg_j'
$$
such that $\varphi=\delta\ugol d\om$. Take $\sigma\in G_j'$ for which \eqref{(1.36)} is satisfied and $\sigma x=x+\delta x$ for all $x\in E$ (Proposition~\ref{1.24}). Write $\,\om=dy+\Tilde\om\,$,  where $y\in E$, $\,\Tilde\om\in\mOm^1$. Then
$$
\sigma\om-\om-\delta\om=\sigma\Tilde\om-\Tilde\om-\delta\Tilde\om\in\fm^{(j+2)}\Om^1.
$$
Since $\om(\delta)=0$, we have 
\[
\delta\om=\delta\ugol d\om+d(\delta\ugol\om)=\varphi,
\]
i.e., $\,\sigma\om-\om-\varphi\in\fm^{(j+2)}\Om^1$. Together with 
\eqref{(3.2)}, this yields 
$$
g\tau\om'-\sigma\om\in\fm^{(j+2)}\Om^1.
$$
Hence
$$
\sigma^{-1}(g)\,\sigma^{-1}\tau\om'-\om=\sigma^{-1}(g\tau\om'-\sigma\om)\in\fm^{(j+2)}\Om^1.
$$
Here $\,\sigma^{-1}(g)\in 1+\fm^{(j+1)}\,$ and $\,\sigma^{-1}\tau\in G_j'\,$.\qed 

\begin{corollary}\label{3.4}
Two forms $\,\om,\om'\in\Cont(\Fcal)$ such that $\om-\om'\in\fm^{(j+1)}\Om^1$,  where $j\geqslant 1$, are conjugate under $\,(1+\fm^{(j+1)})\rtimes G_j'\,$.
\end{corollary}

The proof is obtained from Lemma~\ref{3.3} by the downward induction on $j$, taking into account that $\fm^{(j+1)}=0$ for $j$  sufficiently large.

\begin{proposition}\label{3.5}
For $p>2$ two forms $\om,\om'\in\Cont(\Fcal)$ are conjugate with respect 
to the group $(1+\fm^{(2)})\rtimes G'$ if and only if the images of $\om$ 
and $\om'$ in $\Om^1/\mOm^1$ coincide, and the images of 
$d\om$ and $d\om'$ coincide in $\Om^2/\mOm^2$.  
\end{proposition}

\noindent
\textbf{Proof.}
For each free $\Ocal$-module $M$ denote by $\pi_M$ the canonical mapping $M\to M/\fm M$, and set $\pi=\pi_{\Ocal}$. Define an action of the group $\Ocal^\times$ on $M/\fm M$ by assigning to each element $f\in\Ocal^\times$ the multiplication by $\pi(f)\in\Ocal/\fm\cong K$. Оbviously, the mappings $\pi_{\Om^1}$ and $\pi_{\Om^2}\circ d$ are $G$-equivariant. They are also $(1+\fm^{(2)})$-equivariant since 
\begin{gather}
\notag
\pi_{\Om^1}(f\om)=\pi(f)\,\pi_{\Om^1}(\om),\\
\begin{aligned}[t]
\pi_{\Om^2}\,d(f\om)&=\pi_{\Om^2}(df\wedge\om+fd\om)=\pi_{\Om^1}(df)\wedge\pi_{\Om^1}(\om)+\pi(f)\,\pi_{\Om^2}(d\om)=\\
&=\pi(f)\bigl(\pi_{\Om^1}(f^{-1}df)\wedge\pi_{\Om^1}(\om)+\pi_{\Om^2}(d\om)\bigr),
\end{aligned}
\label{(3.3)}
\end{gather}
and if $f\in 1+\fm^{(2)}$, then $df\in\mOm^1$, implying $\pi_{\Om^2}\bigl(d(f\om)\bigr)=\pi_{\Om^2}\mskip1mu d(\om)$. The group $(1+\fm^{(2)})\rtimes G'$ acts trivially on $\Om^1/\mOm^1$ and on $\Om^2/\mOm^2$. Therefore, if $\om,\om'\in\Om^1$ are conjugate under this group, then $\pi_{\Om^1}(\om)=\pi_{\Om^1}(\om')$ and $\pi_{\Om^2}\mskip1mu d(\om)=\pi_{\Om^2}\mskip1mu d(\om')$.

Conversely, let $\om,\om'\in\Cont(\Fcal)$ be such that $\pi_{\Om^1}(\om)=\pi_{\Om^1}(\om')$ and $\pi_{\Om^2}\mskip1mu d(\om)=\pi_{\Om^2}\mskip1mu d(\om')$. Then
$$
\om'-\om=\sum b_{ij}x_i dx_j+\psi,\quad\text{where}\quad b_{ij}\in K,\quad b_{ij}=b_{ji},\quad\psi\in\fm^{(2)}\Om^1.
$$
Set
$$
f=\sum_{i<j}b_{ij}x_ix_j+\sum b_{ii}x_i^{(2)}\in\fm^2\subset\fm^{(2)}.
$$
By Lemma~\ref{3.2} applied to $\om$ and $\om+df$ there exists $\sigma\in G'$ 
such that $\om+df-\sigma\om\in\fm^{(2)}\Om^1$. Then
$$
\om'-\sigma\om=\om+df-\sigma\om+\psi\in\fm^{(2)}\Om^1.
$$
By Corollary~\ref{3.4} the forms $\om'$ and $\sigma\om$ are conjugate under the group $(1+\fm^{(2)})\rtimes G'$. Hence, so are $\om'$ and $\om$,  too.

\begin{corollary}\label{3.6}{\normalfont ($p>2$).}
Any contact form is $G'$-conjugate to a form $\,f\cdot
\sum\limits_i\bigl(a_i+\sum\limits_jb_{ij}x_j\bigr)dx_i$,  where $a_i,b_{ij}\in K,$ $f\in\Ocal^\times$ and the matrix $(b_{ij})$ can be assumed to be either triangular, or 
skewsymmetric.
\end{corollary}

\begin{corollary}\label{3.7}{\normalfont ($p>2$).}
Any contact Lie algebra is isomorphic to a contact Lie algebra from~\cite{KSh}\footnote{The contact forms considered in \cite{KSh} are expressed as $\sum a_idx_i+\sum b_{ij}x_idx_j$ with antisymmetric matrix $(b_{ij})$. Note, however, that these forms do not satisfy our definition of a contact form when $p=2$ and $n>1$, while unitriangular matrices $(b_{ij})$ give contact forms in any characteristic, including $p=2$.}.
\end{corollary}

Let us identify $E$ (resp. $\LaE$) with the space of linear (resp. bilinear 
\repl{antisymmetric}) forms on $E^*$. Denote by $\MaE$ the set of pairs 
$(x,b)\in E\times\LaE$ such that
$$
x\ne 0\quad\text{and}\quad E^*=\Ker x\oplus\Ker b.
$$
The direct product of groups $\KtmE$ acts on $\MaE$ by the formula
$$ 
(\lambda,x)\cdot(y,b)=(\lambda y,\lambda b+\lambda x\wedge y).
$$
The spaces $\Om^1/\mOm^1$ and $\Om^2/\mOm^2$ are identified, respectively, with $E$ and $\LaE$ by means of $G$-equivariant isomorphisms. The mapping which associates to $\om\in\Om^1$ the pair $(x,b)$ in which $x$ corresponds to $\pi_{\Om^1}(\om)$ and $b\in\LaE$ corresponds to $\pi_{\Om^2}(d\om)$ sends $\Cont(\Fcal)$ to $\MaE$. The $G$-equivariant mapping $\Cont(\Fcal)\to\MaE$ thus obtained induces a map
\begin{equation}
\label{(3.4)}
\Cont(\Fcal)/\bigl((1+\fm^{(2)})\rtimes G'\bigr)\to\MaE.
\end{equation}
Its injectivity is equivalent to the conclusion of Proposition \ref{3.5}, while surjectivity follows from consideration of the forms written out in Corollary~\ref{3.6}. The group $\Ocal^\times$ acts on 
$$
\Cont(\Fcal)/\bigl((1+\fm^{(2)})\rtimes G'\bigr).
$$
In view of \eqref{(3.3)} the corresponding action on $\MaE$ is performed by means of a surjective homomorphism $\Ocal^\times\to \KtmE$. Thus,  bijection \eqref{(3.4)} yields a bijection
\begin{equation}
\label{(3.5)}
\Cont(\Fcal)/(\OrtmG')\to\MaE/(\KtmE).
\end{equation}
The group $G$ acts on $E$ as $\St_\Fcal$. This leads to a bijection
\begin{equation}
\label{(3.6)}
\Cont(\Fcal)/(\OrtmG)\to\MaE/\bigl((\KtmE)\rtimes\St_\Fcal\bigr).
\end{equation}

\section{\repl{Normal shape} of the \repl{antisymmetric} bilinear form with respect to a flag of subspaces}\label{Section4}\label{4}

Let $V$ be a finite-dimensional vector space over $K$, and let
$$
\mathcal G:0=V_0\subset V_1\subset\ldots\subset V_r=V 
$$
be a flag\footnote{The flag $\Gcal$ will be fixed. By a \textit{flag}\/ we mean a sequence of subspaces which form either increasing or decreasing chain. Equal terms  are allowed in this sequence.}. Denote by $\St_\Gcal$ the subgroup of all linear transformations of $V$ preserving $\Gcal$, and by $N(\Gcal)$ the subgroup of transformations in $\St_\Gcal$ inducing the identity transformations in each factor of the flag $\Gcal$ (we will use this notation for arbitrary flags).

For a vector space $W$ denote by $\Ffrak(W)$ the totality of flags of subspaces in $W$ of the form 
$$
W=W_0\supseteq W_1\supseteq\ldots\supseteq W_r
$$
where $r$ is the same as for $\Gcal$. Define a mapping 
$$
\Ffrak(V)\to\prod_{i=1}^r\,\Ffrak(V_i/V_{i-1})
$$
assigning to a flag $\,V=V_0'\supseteq V_1'\supseteq\ldots\supseteq V_r'\,$ the collection of flags
$$
(V_0'\cap V_i+V_{i-1})/V_{i-1}\supseteq(V_1'\cap V_i+V_{i-1})/V_{i-1}\supseteq\ldots\supseteq(V_r'\cap V_i+V_{i-1})/V_{i-1}
$$ 
with $i=1,\ldots,r$. This mapping is $\St_\Gcal$-equivariant. Since the group $N(\Gcal)$ acts trivially on each set $\Ffrak(V_i/V_{i-1})$, we get a mapping 
\begin{equation}
\label{(4.1)}
\Ffrak(V)/N(\Gcal)\to\prod_{i=1}^r\,\Ffrak(V_i/V_{i-1}).
\end{equation}

\begin{lemma}\label{4.1}
The mapping \eqref{(4.1)} is a $\St_\Gcal$-equivariant bijection. Flags 
$$
V=V_0'\supseteq V_1'\supseteq\ldots\supseteq V_r'\quad\text{and}\quad V=V_0''\supseteq V_1''\supseteq\ldots\supseteq V_r''
$$
are conjugate under\/ $\St_\Gcal$ if and only if 
$$
\dim(V_i\cap V_j')=\dim(V_i\cap V_j'')\quad\text{for all }\,1\leqslant i,j\leqslant r.
$$
\end{lemma}

This is an easy fact from linear algebra.

\medskip
For a bilinear \repl{antisymmetric} form $b\in\LaVd$ denote by $M^{\perp_b}$ the orthogonal of a subspace $M\subseteq V$ with respect to $b$. For subspaces $M,L\subseteq V$ and an arbitrary \repl{antisymmetric} form the following known identities hold:
$$
M^{\perp\perp}=M+V^\perp;\quad(M+L)^\perp=M^\perp\cap L^\perp;\quad(M\cap L^\perp)^\perp=M^\perp+L.
$$

Let $\Gcal^{\perp_b}$ denote the flag
$$
V=V_0^{\perp_b}\supseteq V_1^{\perp_b}\supseteq\ldots\supseteq V_r^{\perp_b}.
$$ 
Thus,  $\Gcal^{\perp_b}\in\Ffrak(V)$. The mapping $\LaVd\to\Ffrak(V)$ such that $b\mapsto\Gcal^{\perp_b}$ is $\St_\Gcal$-equivariant.

The form $b\in\LaVd$ induces non-degenerate pairings 
$$
V_i\cap V_{j-1}^{\perp_b}/(V_{i-1}\cap V_{j-1}^{\perp_b}+V_i\cap V_j^{\perp_b})\times V_j\cap V_{i-1}^{\perp_b}/(V_{j-1}\cap V_{i-1}^{\perp_b}+V_j\cap V_i^{\perp_b})\to K.
$$
Indeed,
$$
(V_i\cap V_{j-1}^{\perp_b})^{\perp_b}\cap V_j\cap V_{i-1}^{\perp_b}=(V_i^{\perp_b}+V_{j-1})\cap V_j\cap V_{i-1}^{\perp_b}=(V_i^{\perp_b}\cap V_j+V_{j-1})\cap V_{i-1}^{\perp_b}=V_i^{\perp_b}\cap V_j+V_{j-1}\cap V_{i-1}^{\perp_b}.
$$
With the help of isomorphisms
\begin{gather*}
(V_i\cap V_{j-1}^{\perp_b}+V_{i-1})/(V_i\cap V_j^{\perp_b}+V_{i-1})\cong V_i\cap V_{j-1}^{\perp_b}/(V_i\cap V_j^{\perp_b}+V_{i-1}\cap V_{j-1}^{\perp_b}),\\
(V_j\cap V_{i-1}^{\perp_b}+V_{j-1})/(V_j\cap V_i^{\perp_b}+V_{j-1})\cong V_j\cap V_{i-1}^{\perp_b}/(V_j\cap V_i^{\perp_b}+V_{j-1}\cap V_{i-1}^{\perp_b})
\end{gather*}
we get non-degenerate pairings
\begin{equation}
\label{(4.2)}
(V_i\cap V_{j-1}^{\perp_b}+V_{i-1})/(V_i\cap V_j^{\perp_b}+V_{i-1})\times(V_j\cap V_{i-1}^{\perp_b}+V_{j-1})/(V_j\cap V_i^{\perp_b}+V_{j-1})\to K.
\end{equation}

\medskip
Denote by $X$ the set of tuples $\bigl((\Gcal_i)_{1\leqslant i\leqslant r},(b_{ij})_{1\leqslant i,j\leqslant r}\bigr)$,  where $\Gcal_i\in\Ffrak(V_i/V_{i-1})$, say $\Gcal_i$ is a flag
$$
V_i/V_{i-1}=U_{i0}\supseteq U_{i1}\supseteq\ldots\supseteq U_{ir}\,,\\
$$
and
$$
b_{ij}:(U_{i,j-1}/U_{ij})\times(U_{j,i-1}/U_{ji})\to K
$$
are non-degenerate pairings such that $b_{ij}(x,y)=-b_{ji}(y,x)$ and 
$b_{ii}(x,x)=0$. The group
$$
\prod_{i=1}^r\mskip1mu GL(V_i/V_{i-1})\cong\St_\Gcal/N(\Gcal)
$$
acts on $X$ in a natural way.

To every form $b\in\LaVd$ we assign the collection of flags 
$$
(V_i\cap V_0^{\perp_b}+V_{i-1})/V_{i-1}\supseteq (V_i\cap V_1^{\perp_b}+V_{i-1})/V_{i-1}\supseteq\ldots\supseteq (V_i\cap V_r^{\perp_b}+V_{i-1})/V_{i-1}
$$
and pairings \eqref{(4.2)}. We have thus constructed a $\St_\Gcal$-equivariant mapping $\LaVd\to X$. 
Since $N(\Gcal)$ acts on $X$ trivially, we get a mapping 
\begin{equation}
\label{(4.3)}
\bigl(\LaVd\bigr)/N(\Gcal)\to X.
\end{equation}

\begin{proposition}\label{4.2}The mapping \eqref{(4.3)} is a bijection.
\end{proposition}

\noindent
\textbf{Proof.} Let $(\Gcal_i,b_{ij})\in X$,  where $\Gcal_i$ is a flag 
$$
V_i/V_{i-1}=U_{i0}\supseteq U_{i1}\supseteq\ldots\supseteq U_{ir}\,.
$$
For every pair $i,j$, where $1\leqslant i\leqslant r$ and $1\leqslant j\leqslant r+1$, pick a basis $\{\overline {\overline e}\,_{ij}^k\mid 1\leqslant k\leqslant n_{ij}\}$ of the space $U_{i,j-1}/U_{ij}$, assuming that $U_{i,\mskip1mu r+1}=0$. Choose pre-images $\overline e\,_{ij}^k$ of $\overline{\overline e}\,_{ij}^k$ in $U_{i,j-1}$ and pre-images $e_{ij}^k$ of $\overline e\,_{ij}^k$ in $V_i$. Thus,  for each $1\leqslant i\leqslant r$ the elements $\{\overline e\,_{ij}^k\mid 1\leqslant j\leqslant r+1,\ 1\leqslant k\leqslant n_{ij}\}$ constitute a basis of $V_i/V_{i-1}$ compatible with the flag $\Gcal_i$, and the elements 
$$
\{e_{ij}^k\mid 1\leqslant i\leqslant r,\ 1\leqslant j\leqslant r+1,\ 1\leqslant k\leqslant n_{ij}\}
$$
constitute a basis of $V$ compatible with $\Gcal$.

Define a bilinear \repl{antisymmetric} form $b:V\times V\to K$, specifying its values on the basis elements by the rule
$$
b(e_{ij}^k,e_{\ell m}^n)=\left\{\begin{array}{cl}
b_{ij}(\overline{\overline e}\,_{ij}^k,\mskip1mu\overline{\overline e}\,_{\ell m}^n)&\text{if}\ \,i=m,\ j=\ell,\\
\noalign{\smallskip}
0&\text{otherwise}.
\end{array}
\right.
$$
Since the pairings $b_{ij}$ are non-degenerate, it follows that the space $V_{s-1}^{\perp_b}$ for $1\leqslant s\leqslant r+1$ is spanned by $\{e_{ij}^k\mid j\geqslant s\}$. It is clear now that mapping \eqref{(4.3)} sends $b$ to the collection $(\Gcal_i,b_{ij})$. Thus,  mapping \eqref{(4.3)} is surjective.

Let $b''\in\LaVd$ be another form which is sent by mapping \eqref{(4.3)} to the same collection $(\Gcal_i,b_{ij})$. Since the diagram
$$
\begin{CD}
\LaVd@>(4.3)>>X\\
@VVV@VVV\\
\Ffrak(V)@>(4.1)>>\prod\limits_{i=1}^r\Ffrak(V_i/V_{i-1})
\end{CD}
$$
is commutative, we see that mapping \eqref{(4.1)} sends the flag 
$\Gcal^{\perp_{b''}}$ to $(\Gcal_i)_{1\leqslant i\leqslant r}$. By 
Lemma~\ref{4.1} there exists $\sigma\in N(\Gcal)$ such that 
$\sigma\Gcal^{\perp_{b''}}=\Gcal^{\perp_b}$. Then for $b'=\sigma b''$ we have 
$\Gcal^{\perp_{b'}}=\Gcal^{\perp_b}$. Now, set 
$\Gcal^\perp=\Gcal^{\perp_b}=\Gcal^{\perp_{b'}}$ and 
$V_i^\perp=V_i^{\perp_b}=V_i^{\perp_{b'}}$ for short. Thus,
$$
(V_i\cap V_j^\perp+V_{i-1})/V_{i-1}=U_{ij}.
$$
It remains to show that the forms $b$ and $b'$ are conjugate under $N(\Gcal)$. Let
$$
I=\{(i,j,k)\mid 1\leqslant i\leqslant r,\ 1\leqslant j\leqslant r+1,\ 1\leqslant k\leqslant n_{ij}\}.
$$
It suffices to find elements $\elp_{ij}^k\in V_i\cap V_{j-1}^\perp$ for all $(i,j,k)\in I$ such that $\elp_{ij}^k-e_{ij}^k\in V_{i-1}$ and
\begin{equation}
\label{(4.4)}
b'(\mskip1mu\elp_{ij}^k,\elp_{\ell m}^n)=\left\{\begin{array}{cl}
b_{ij}(\overline{\overline e}\,_{ij}^k,\mskip1mu\overline{\overline e}\,_{\ell m}^n)&\text{if}\ \,i=m,\ j=\ell,\\
\noalign{\smallskip}
0&\text{otherwise}.
\end{array}
\right.
\end{equation}
Indeed, such elements $\elp_{ij}^k$ form a basis of $V$ compatible with the flag $\Gcal$. The linear operator $\tau:V\to V$ sending $e_{ij}^k$ to $\elp_{ij}^k$ belongs to $N(\Gcal)$ and takes $b$ to $b'$.

Introduce a linear order on $I$ by setting $(i,j,k)<(\ell,m,n)$ if and only if
$$
\begin{array}{rl}
\text{either}&\min(i,j)<\min(\ell,m),\\
\noalign{\smallskip}
\text{or}&i=\min(i,j)=\min(\ell,m)<\ell,\\
\noalign{\smallskip}
\text{or}&i=\min(i,j)=\min(\ell,m)=\ell,\quad j<m,\\
\noalign{\smallskip}
\text{or}&\min(i,j)=\min(\ell,m)<i<\ell,\\
\noalign{\smallskip}
\text{or}&i=\ell,\quad j=m,\quad k<n.
\end{array}
$$

Suppose that $(s,t,q)\in I$, the elements $\elp_{ij}^k\in V_i\cap V_{j-1}^\perp$ such that $\elp_{ij}^k\equiv e_{ij}^k\pmod{V_{i-1}}$ have been found for all $(i,j,k)<(s,t,q)$, and condition \eqref{(4.4)} holds for all $(i,j,k),(\ell,m,n)<(s,t,q)$. 

We will construct $\elp_{st}^q\in V_s\cap V_{t-1}^\perp$ such that $\elp_{st}^q\equiv e_{st}^q\pmod{V_{s-1}}$ and \eqref{(4.4)} holds for any ${(i,j,k)<(s,t,q)}$ and $(\ell,m,n)=(s,t,q)$. In this way all elements $\elp_{ij}^k$ with required properties will be constructed by induction.

Set
\begin{align*}
J &=\{(i,j,k)\in I\mid(i,j,k)<(s,t,q),\quad i\geqslant t,\quad j<s\},\\
J'&=\{(i,j,k)\in I\mid(i,j,k)<(s,t,q),\quad i<s,\quad j\geqslant t\}.
\end{align*}
We claim that if $(i,j,k)\in J$, then $(j,i,\ell)\in J'$ for all $1\leqslant\ell\leqslant n_{ji}$. We have to verify that $(j,i,\ell)<(s,t,q)$.

If $j<t$, then $\min(i,j)<\min(s,t)$, and therefore $(j,i,\ell)<(s,t,q)$.

Let $j\geqslant t$. Since $j<s$, we have $t<s$ and $\min(s,t)=t$. Since $(i,j,k)<(s,t,q)$, but $\min(i,j)\geqslant t$, we have either $i=t$, or $j=t$ and $i\leqslant s$. 

Now, $(j,i,\ell)<(s,t,q)$ in the first case because $\min(j,i)=\min(s,t)\leqslant j<s$, and in the second case because $j=\min(j,i)=\min(s,t)<s$.

Let
$$
P=\sum_{(i,j,k)\in J}K\,\elp_{ij}^k,\qquad Q=\sum_{(i,j,k)\in J'}K\,\elp_{ij}^k.
$$
We claim that the pairing $b'|_{P\times Q}$ has no left kernel, and therefore the mapping $Q\to P^*$ induced by $b'$ is surjective. Let us write $P=\bigoplus P_{ij}$ and $Q=\bigoplus Q_{ij}$,  where $P_{ij}$ (resp. $Q_{ij}$) is spanned by the vectors $\elp_{ij}^k$ with $(i,j,k)\in J$ (resp. $J'$) for the given $i,j$. By relation \eqref{(4.4)} $\,b'(P_{ij},Q_{\ell m})=0$ for $i\ne m$ or $j\ne\ell$. At the same time the left kernel of the pairing $b'|_{P_{ij}\times Q_{ji}}$ is zero, as follows from \eqref{(4.4)}, non-degeneracy of the pairings $b_{ij}$, and the fact that $\elp_{ji}^\ell\in Q_{ji}$ for all $\ell$ such that $1\leqslant\ell\leqslant n_{ji}$ whenever $P_{ij}\ne 0$.

Therefore, there exists $x\in Q$ such that $b'(e_{st}^q,y)=b'(x,y)$ for all $y\in P$. Set $\elp_{st}^q=e_{st}^q-x$. If $(i,j,k)\in J'$, then $\elp_{ij}^k\in V_i\cap V_{j-1}^\perp\subseteq V_{s-1}\cap V_{t-1}^\perp$, whence $x\in V_{s-1}\cap V_{t-1}^\perp$. Since $e_{st}^q\in V_s\cap V_{t-1}^\perp$, we have $\elp_{st}^q\in V_s\cap V_{t-1}^\perp$ and $\elp_{st}^q\equiv e_{st}^q\pmod{V_{s-1}}$.

The definition of $x$ implies that $b'(\elp_{st}^q,\elp_{ij}^k)=0$ if $(i,j,k)\in J $. Let $(i,j,k)<(s,t,q)$, but $(i,j,k)\notin J$, i.e., either $i<t$ or $j\geqslant s$.

If $i<t$, then $b'(\elp_{st}^q,\elp_{ij}^k)=0$ because $\elp_{st}^q\in V_{t-1}^\perp$ and $\elp_{ij}^k\in V_i\subseteq V_{t-1}$. If $j>s$, then $b'(\elp_{st}^q,\elp_{ij}^k)=0$ because $\elp_{st}^q\in V_s$ and $\elp_{ij}^k\in V_{j-1}^\perp\subseteq V_s^\perp$.

Let $j=s$, but $i\geqslant t$. By condition $(i,j,k)<(s,t,q)$ we must have $i=t$. Note that for all $u\in V_s\cap V_{t-1}^\perp$ and $v\in V_t\cap V_{s-1}^\perp$ it follows from the definition of $b_{st}$ that\footnote{
More correctly, the equality here follows from the fact that $b'$ is $N(\Gcal)$-conjugate to $b''$. Mapping (4.3) takes $b'$ to the given collection $(\Gcal_i,b_{ij})$.} $\,b'(u,v)=b_{st}(\overline u,\overline v)$,  where $\overline u$ (resp. $\overline v$) is the class of $u$ (resp. $v$) in
$$
(V_s\cap V_{t-1}^\perp+V_{s-1})/(V_s\cap V_t^\perp+V_{s-1})\cong U_{s,t-1}/U_{st}\quad\text{ (resp. in }\,U_{t,s-1}/U_{ts}\text{)}.
$$
Therefore, $\,b'(\elp_{st}^q,\elp_{ts}^k)=b_{st}(\overline{\overline e}\,_{st}^q,\overline{\overline e}\,_{ts}^k)$.

Thus, formula \eqref{(4.4)} has been verified for $(i,j,k)<(s,t,q)$ and $(\ell,m,n)=(s,t,q)$. The proof is complete.\qed 

\begin{proposition}\label{4.3}
Forms $b,b'\in\LaVd$ are conjugate under $\St_\Gcal$ if and only if
$$
\dim(V_i\cap V_j^{\perp_b})=\dim(V_i\cap V_j^{\perp_{b'}})
$$
for all $i,j$ such that $1\leqslant i,j\leqslant r$.
\end{proposition}

\noindent
\textbf{Proof.}
Let us verify injectivity of the mapping $(\LaVd)/\St_\Gcal\to\Ffrak(V)/\St_\Gcal$ induced by the mapping
$$
\LaVd\to\Ffrak(V),\quad b\mapsto\Gcal^{\perp_b}.
$$
The horizontal arrows in the commutative diagram 
$$
\begin{CD}
\LaVd/\St_\Gcal @>>>(\LaVd/N(\Gcal))/(\St_\Gcal /N(\Gcal))@>>> X \Big /\prod\limits_{i=1}^r GL (V_i/V_{i-1})\\
@VVV@VVV@VVV\\
\Ffrak(V)/\St_\Gcal @>>>(\Ffrak(V)/N(\Gcal))/(\St_\Gcal /N(\Gcal))@>>>\prod\limits_{i=1}^r\Ffrak(V_i/V_{i-1}) \Big /\prod\limits_{i=1}^r GL (V_i/V_{i-1})
\end{CD}
$$
are bijections. Therefore, it suffices to prove that the rightmost downward arrow is an embedding. In other words, we have to verify that any two points of the set $X$ lying over the same point of the set $\prod\limits_{i=1}^r\Ffrak(V_i/V_{i-1})$ are conjugate under $\prod\limits_{i=1}^r GL (V_i/V_{i-1})$.

Let $\Gcal_i:V_i/V_{i-1}=U_{i0}\supseteq U_{i1}\supseteq\ldots\supseteq U_{ir}$ be flags in $\Ffrak(V_i/V_{i-1})$ for $i=1,\ldots,r$. Let $b_{ij}$ and $b'_{ij}$ be non-degenerate pairings $U_{i,j-1}/U_{ij}\times U_{j,i-1}/U_{ji}\to K$ satisfying the identities
$$
b_{ij}(x,y)=-b_{ji}(y,x),\qquad b'_{ij}(x,y)=-b'_{ji}(y,x).
$$
There exist nonsingular endomorphisms $A_{ij}$ of the spaces $U_{i,j-1}/U_{ij}$ for $1\leqslant i,j\leqslant r$ such that $b_{ij}(x,y)=b'_{ij}(A_{ij}x,A_{ji}y)$. Pick endomorphisms $A_i\in\St_{\Gcal_i}$ which induce $A_{ij}$ in the quotients $U_{i,j-1}/U_{ij}$. Then the element $(A_1,\ldots,A_r)$ of the group $\prod\limits_{i=1}^rGL (V_i/V_{i-1})$ sends the point $\bigl((\Gcal_i),(b_{ij})\bigr)\in X$ to the point $\bigl((\Gcal_i),(b'_{ij})\bigr)$.\qed

\medskip
Note that the classification of bilinear forms relative to a flag all whose factors are of dimension equal to 1 is given in~\cite[Ch. IX, Section 3, Th. I]{2} . 

The numbers 
$$
n_{ij}=\dim V_i\cap V_{j-1}^\perp /(V_i\cap V_j^\perp+V_{i-1}\cap V_{j-1}^\perp)
$$
are more convenient for us than the dimensions of spaces $V_i\cap V_j^\perp$.

\begin{corollary}\label{4.4}
The $\St_\Gcal$-orbits of bilinear \repl{antisymmetric} forms on $V$ are in a one-to-one correspondence with the tuples of non-negative integers $(n_{qt})_{1\leqslant q,t\leqslant r}$ in which $n_{qq}$ are even{\rm,} $n_{qt}=\nobreak n_{tq},$ and $\,\sum\limits_{t=1}^rn_{qt}\leqslant\dim V_q/V_{q-1}\,$. The tuples with $\,\sum\limits_{t=1}^rn_{qt}=\dim(V_q/V_{q-1})$ correspond to the orbits of non-degenerate forms.
\end{corollary}

\begin{lemma}\label{4.5}Let $1\leqslant k\leqslant r$. The orbits of the $\St_\Gcal$-action on the set 
$$
\LaVdnd\times\bigl((V_k/V_{k-1})^*\,\diagdown\{0\}\bigr)
$$
are in a one-to-one correspondence with the collections of integers $\,\ell,\ (n_{qt})_{1\leqslant q,t\leqslant r}\,$ satisfying the following conditions\/{\rm:}
\begin{equation}
\label{(4.5)}
\begin{array}{c}
n_{qt}\geqslant 0;\quad n_{k\ell}\ne 0;\quad n_{qt}=n_{tq};\quad 1\leqslant\ell\leqslant r;\\
\noalign{\medskip}
n_{qq}\text{ are even};\quad\sum\limits_tn_{qt}=\dim(V_q/V_{q-1}).
\end{array}
\end{equation}
\end{lemma}

\noindent
\textbf{Proof.}
To $(b,f)\in\LaVdnd\times\bigl((V_k/V_{k-1})^*\,\diagdown\{0\}\bigr)$ 
we assign the collection $n_{qt}$ of invariants of the form $b$ (see Corollary~\ref{4.4}) and the number 
$$
\ell=1+\max\{t\mid f\bigl((V_k\cap V_{t-1}^{\perp_b}+V_{k-1})/V_{k-1}\bigr)\ne 0\}.
$$
Then the numbers $n_{qt}$ and $\ell$ depend only on the $\St_(\Gcal)$-orbit, and \eqref{(4.5)} holds. 

Suppose $(b',f')\in\LaVdnd\times\bigl((V_k/V_{k-1})^*\,\diagdown\{0\}\bigr)$ has the same invariants $\ell$ and $n_{qt}$ as $(b,f)$. Let us prove that $(b,f)$ and $(b',f')$ are conjugate under $\St_\Gcal$.

By Corollary~\ref{4.4} the forms $b$, $b'$ are conjugate under means of $\St_\Gcal$. Therefore, we may assume at once that $b'=b$. Denote by $X_\#$ the subset of tuples $(\Gcal_i,b_{ij})\in X$ in which the flags 
$$
\Gcal_i:U_{i0}\supseteq\ldots\supseteq U_{ir}
$$
end with $U_{ir}=0$. Since the form $b$ is non-degenerate, the mapping \eqref{(4.3)} sends it into some tuple
$$ 
\kappa:=(\Gcal_i,b_{ij})\in X_\#. 
$$
By Proposition~\ref{4.2} we have bijections 
\begin{gather*}
\bigl(\LaVdnd\times\bigl((V_k/V_{k-1})^*\,\diagdown\{0\}\bigr)\bigr)/N(\Gcal)\to X_\#\times\bigl((V_k/V_{k-1})^*\diagdown\{0\}\bigl),\\ 
\bigl(\LaVdnd\times\bigl((V_k/V_{k-1})^*\,\diagdown\{0\}\bigr)\bigr)/\St_\Gcal\to\bigl(X_\#\times\bigl((V_k/V_{k-1})^*\diagdown\{0\}\bigr)\bigr)\big/{\textstyle\prod\limits_{q=1}^r}GL(V_q/V_{q-1}).
\end{gather*}
It remains to show that $(\kappa,f)$ and $(\kappa,f')$ are conjugate under $\prod\limits  GL (V_q/V_{q-1})$. By definition of $\ell$ the linear forms $f$ and $f'$ induce non-zero forms $\overline f$ and $\overline{f'}$ on $U_{k,\ell -1}/U_{k\ell}$. There exists an invertible endomorphism $A_{k\ell}$ of the space $U_{k,\ell -1}/U_{k\ell}$ sending $\overline f$ to $\overline{f'}$. If $k=\ell$, we may assume that $A_{kk}$ preserves the \repl{antisymmetric} form $b_{kk}$ since the symplectic group acts transitively on the set of non-zero linear forms. It follows from the definition of $\ell$ that the mappings 
\begin{align*}
(\Ker f\cap U_{k,t-1})/(\Ker f\cap U_{kt})&\to U_{k,t-1}/U_{kt}\,,\\
(\Ker f'\cap U_{k,t-1})/(\Ker f'\cap U_{kt})&\to U_{k,t-1}/U_{kt}
\end{align*}
are isomorphisms for $t\ne\ell$. Using them we find, for $t<\ell$, bases $\{e_t^h\mid h=1,\ldots,n_{kt}\}$ and $\{\elp_t^h\mid h=1,\ldots,n_{kt}\}$ of some complements of $U_{kt}$ to $U_{k,t-1}$ such that 
$$
f(e_t^h)=0,\quad f'(\elp_t^h)=0,\quad\elp_t^h-e_t^h\in U_{kt}\,.
$$
There exists an endomorphism $\Tilde A_k\in\St_{\Gcal_k}$ which acts as the identity on $U_{k\ell}$, induces $A_{k\ell}$ on the quotient $U_{k,\ell -1}/U_{k\ell}$, and sends $e_t^h$ to $\elp_t^h$. Then $\Tilde A_kf=f'$, and in the quotients $U_{k,t-1}/U_{kt}$, for $t\ne\ell$, the endomorphism $\Tilde A_k$ induces the identity transformations. 

Define endomorphisms $A_{qt}$ of the spaces $U_{q,t-1}/U_{qt}$ to be the identity transformations if $(q,t)\ne(k,\ell)$ and $(q,t)\ne(\ell,k)$, while $A_{\ell k}$ in the case where $\ell\ne k$ is determined from the identity
$$
b_{k\ell}(A_{k\ell}^{-1}u,v)=b_{k\ell}(u,A_{\ell k}v)\qquad
u\in U_{k,\ell -1}/U_{k\ell},\ \,v\in U_{\ell,k-1}/U_{\ell k}.
$$
For $q\ne k$ let $\Tilde A_q$ be an arbitrary endomorphism in $\St_{\Gcal_q}$ which induces $A_{qt}$ in the quotients $U_{q,t-1}/U_{qt}\,$. Now the element $(\Tilde A_1,\ldots,\Tilde A_r)\in\prod GL(V_q/V_{q-1})$ sends $f$ to $f'$ and leaves $k$ unchanged.\qed

\medskip
We will need an analogue of Proposition~\ref{4.2} for the case of pairings of two spaces. Let $V$ and $W$ be finite-dimensional vector spaces,
$$
\Gcal:0=V_0\subseteq V_1\subseteq\ldots\subseteq V_r=V\quad\text{and}\quad\Hcal:0=W_0\subseteq W_1\subseteq\ldots\subseteq W_r=W
$$
flags of subspaces, $B(V,W)$ the space of bilinear mappings $V\times W\to K$. As earlier, $b\in B(V,W)$ induces non-degenerate pairings 
\begin{equation}
\label{(4.6)}
\begin{array}{@{}c@{}}
V_i\cap W_{j-1}^{\perp_b}/(V_{i-1}\cap W_{j-1}^{\perp_b}+V_i\cap W_j^{\perp_b})\times W_j\cap V_{i-1}^{\perp_b}/(W_{j-1}\cap V_{i-1}^{\perp_b}+W_j\cap V_i^{\perp_b})\to K;\\
\noalign{\medskip}
(V_i\cap W_{j-1}^{\perp_b}+V_{i-1})/(V_i\cap W_j^{\perp_b}+V_{i-1})\times(W_j\cap V_{i-1}^{\perp_b}+W_{j-1})/(W_j\cap V_i^{\perp_b}+W_{j-1})\to K.
\end{array}
\end{equation}

Denote by $Y$ the set of tuples $\bigl((\Gcal_i,\Hcal_i)_{1\leqslant i\leqslant r},
(b_{ij})_{1\leqslant i,j\leqslant r}\bigr)$,  where $\Gcal_i\in\Ffrak(V_i/V_{i-1})$ 
and\hfil\break 
${\Hcal_i\in\Ffrak(W_i/W_{i-1})}$; say $\Gcal_i$ and $\Hcal_i$ are, 
respectively, flags 
$$
V_i/V_{i-1}=T_{i0}\supseteq T_{i1}\supseteq\ldots\supseteq T_{ir},\qquad W_i/W_{i-1}=U_{i0}\supseteq U_{i1}\supseteq\ldots\supseteq U_{ir},
$$
and $b_{ij}:T_{i,j-1}/T_{ij}\times U_{j,i-1}/U_{ji}\to K$ are non-degenerate pairings.

The group $\,\prod\limits_{i=1}^r\bigl(GL(V_i/V_{i-1})\times GL(W_i/W_{i-1})\bigr)\,$ acts on $Y$.

To a form $b\in B(V,W)$ we assign the collections of flags 
\begin{gather*}
V_i\cap W_0^{\perp_b}+V_{i-1}/V_{i-1}\supseteq V_i\cap W_1^{\perp_b}+V_{i-1}/V_{i-1}\supseteq\ldots\supseteq V_i\cap W_r^{\perp_b}+V_{i-1}/V_{i-1},\\
W_i\cap V_0^{\perp_b}+W_{i-1}/W_{i-1}\supseteq W_i\cap V_1^{\perp_b}+W_{i-1}/W_{i-1}\supseteq\ldots\supseteq W_i\cap V_r^{\perp_b}+W_{i-1}/W_{i-1},
\end{gather*}
and pairings \eqref{(4.6)}. We have thus constructed a $\St_{\Gcal}\!\times\St_{\Hcal}$-equivariant mapping $B(V,W)\to Y$. Since $N(\Gcal)\times N(\Hcal)$ acts trivially on $Y$, we get a mapping
\begin{equation} 
\label{(4.7)}
B(V,W)/\bigl(N(\Gcal)\times N(\Hcal)\bigr)\to Y.
\end{equation}

\begin{proposition}\label{4.6}The mapping \eqref{(4.7)} is a bijection.
\end{proposition}

The proof of Proposition~\ref{4.2} is repeated almost verbatim. Later on we will apply versions of Propositions~\ref{4.2} and~\ref{4.6} for decreasing flags $\Gcal$, $\Hcal$ as well.

\medskip
Consider the set $\MaVd$ and the action of the group $\KtmVd$ on it defined in Section~{3}. The group $GL(V)$, and therefore the semi-direct product $(\KtmVd)\rtimes GL(V)$, also act on $\MaVd$.

\begin{lemma}\label{4.7}
The orbits under the action of $\,(\KtmVd)\rtimes\St_\Gcal\,$ on $\MaVd$ are in a one-to-one correspondence with the collections of integers $\bigl(k,(n_{qt})_{1\leqslant q,t\leqslant r}\bigr)$ such that
\begin{equation}
\label{(4.8)}
\begin{array}{l}
1\leqslant k\leqslant r;\quad n_{qt}\geqslant 0,\\
\noalign{\smallskip}
n_{qq}\text{ are even},\quad n_{qt}=n_{tq},
\end{array}
\quad\sum_{t=1}^rn_{qt}=\begin{cases}
\dim(V_q/V_{q-1})-1&\text{for }\,q=k,\\
\noalign{\smallskip}
\dim(V_q/V_{q-1})&\text{for }\,q\ne k.
\end{cases}
\end{equation}
\end{lemma}

\noindent
\textbf{Proof.} Set $Q=\Ker f$. Assign to $(f,b)\in\MaVd$ the collection of numbers
\begin{gather}
\label{(4.9)}
k=\min\{q\mid f(V_q)\ne 0\},\\
\label{(4.10)}
n_{qt}=\dim\,\bigl(Q\cap V_q\cap(Q\cap V_{t-1})^{\perp_b}\bigr)/\bigl(Q\cap V_q\cap(Q\cap V_t)^{\perp_b}+Q\cap V_{q-1}\cap(Q\cap V_{t-1})^{\perp_b}\bigr).
\end{gather}
They are invariant under the action of the group $(\KtmVd)\rtimes\St_\Gcal$.

Suppose that $(f',b')\in\MaVd$ has the same set of invariants $k,n_{qt}$ as $(f,b)$. Let us prove that $(f,b)$ and $(f',b')$ are conjugate. Set $Q'=\Ker f'$. The definition of $k$ implies that
$$
\dim Q\cap V_q\,/Q\cap V_{q-1}=\dim Q'\cap V_q\,/Q'\cap V_{q-1}=
\begin{cases}
(\dim V_q/V_{q-1})-1&\text{if }q=k\\
\dim V_q/V_{q-1}&\text{if }q\ne k.
\end{cases}
$$
Hence, $\dim(Q\cap V_q)=\dim(Q'\cap V_q)$ for all $q$. Therefore, there exists $A_1\in\St_\Gcal$ such that $A_1Q=Q'$. Then $\mu A_1f=f'$ for some $\mu\in K^\times$.

Thus, we may assume that $f'=f$ and $Q'=Q$. By Corollary~\ref{4.4} there exists $C\in GL(Q)$ which preserves the subspaces $V_q\cap Q$ and sends $b|_Q$ to $b'|_Q$. Let us extend $C$ to some $A_2\in\St_\Gcal$ so that $A_2$ induce the identity transformation in $V/Q$. Then $A_2f=f$ and $(A_2b)|_Q=b'|_Q$.

Thus,  we may assume that $f'=f$ and $b'|_Q=b|_Q$. Since $V=V^{\perp_b}\oplus Q$ by the definition of $\MaVd$, we have $f(e)=1$ for a suitable $e\in V^{\perp_b}$. Define $g\in V^*$ by the formula $g(v)=b'(e,v)$. Then $b'=b+f\wedge g$.\qed 

\section{\repl{The normal shape} of the $2$nd-type \repl{symplectic} forms and the contact forms}\label{5}

\begin{theorem}\label{5.1}{\normalfont ($K$ is perfect of characteristic $p>2$).}
Any $2$nd-type \repl{symplectic} form \repl{corresponding} to the flag $\Fcal$ is conjugate under $G(\Fcal)$ to a form expressed as
\begin{equation}
\label{(5.1)}
d\bigl(\exp(x_{i_0})(x_{i_1}dx_{i_1'}+\ldots+x_{i_s}dx_{i_s'})\bigr),
\end{equation}
where $i_0\in \{1,\ldots,n\}$ is an index, $\,\{i_1i_1'\},\ldots,\{i_si_s'\}$ a partition of the set $\{1,\ldots,n\}$ into pairs, $s=n/2$.

In order that the form \eqref{(5.1)} be conjugate to the form
\begin{equation}
\label{(5.2)}
\lambda\,d\bigl(\exp(x_{j_0})(x_{j_1}dx_{j_1''}+\ldots+x_{j_s}dx_{j_s''})\bigr),
\end{equation}
where $j_0,\,\{j_1j_1''\},\ldots,\{j_sj_s''\}$ have the same meaning and $\lambda\in K^\times,$ it is necessary and sufficient that there exist a permutation of indices which sends $i_0$ to $j_0,$ every pair in the partition $\{i_1i_1'\},\ldots,\{i_si_s'\}$ to a pair occurring in the second partition $\{j_1j_1''\},\ldots,\{j_sj_s''\},$ and maps each subset of indices corresponding to indeterminates\footnote{By indeterminates we mean the chosen basis elements $x_1,\ldots,x_n\in E$.} of the same height into itself.
\end{theorem}

\noindent
\textbf{Proof.} Set 
$$
V=E^*,\qquad V_q=\{v\in V\mid\langle v,E_q\rangle=0\},\qquad\Gcal=(V_q)_{0\leqslant q\leqslant m}.
$$
Then $\,\St_\Fcal\cong\St_\Gcal$, $\,E_{k-1}/E_k\cong(V_k/V_{k-1})^*$, and the set in the right-hand side of \eqref{(2.11)} is nothing else but
$$
\coprod_{k\geqslant 1}\,\Bigl(\LaVdnd\times\bigl((V_k/V_{k-1})^*\,\diagdown\{0\}\bigr)\Bigr)/\St_\Gcal.
$$
By Lemma~\ref{4.5} this set and, thanks to bijection \eqref{(2.11)}, the set $\Ham^2(\Fcal)/G(\Fcal)$ are parameterized by the collection of invariants $(k,\ell,n_{qt})$, where $1\leqslant k\leqslant m$, and $(\ell,n_{qt})$ satisfies the relations \eqref{(4.5)} for $r=m$. (We assume that $E_m=0$.)

Let us compute these invariants for the form \eqref{(5.1)}. The mapping \eqref{(2.10)} sends this form to
$$
(x_{i_1}\wedge x_{i_1'}+\ldots+x_{i_s}\wedge x_{i_s'},\,\overline x_{i_0})\in\LaEnd\times\bigl(E_{m_{i_0}-1}/E_{m_{i_0}}\bigr),
$$
where $\overline x_{i_0}$ is the class of $x_{i_0}$ in $E_{m_{i_0}-1}/E_{m_{i_0}}$. In particular, $k=m_{i_0}$. Let $v_1,\ldots,v_n$ be the basis of $V$ dual to the basis $x_1,\ldots,x_n$. The partition $\{i_1i_1'\},\ldots,\{i_si_s'\}$ determines an involutive permutation $i\mapsto i'$ on the set $\{1,\ldots,n\}$. The elements $\{v_i\mid m_i\leqslant q\}$ form a basis of $V_q$ and the elements $\{v_i\mid m_{i'}\geqslant t\}$ form a basis of $V_{t-1}^\perp$. Thus, $\ell=m_{i_0'}$. The classes of elements $\{v_i\mid m_i=q,\ m_{i'}=t\}$ form a basis of
$$
(V_q\cap V_{t-1}^\perp)/(V_q\cap V_t^\perp+V_{q-1}\cap V_{t-1}^\perp).
$$
This implies that $n_{qt}=\card\{i\mid m_i=q,\ m_{i'}=t\}$. Observe that the same set of invariants corresponds to a form differing from the form \eqref{(5.1)} by a non-zero scalar factor.    

Every given collection $k,\ell,n_{qt}$ with necessary properties corresponds to a form of shape \eqref{(5.1)} which can be constructed as follows. Split $\{1,\ldots,n\}$ into a disjoint union of subsets $I_{qt}$, where $I_{qt}\subseteq\{i\mid m_i=q\}$ and $\card I_{qt}=n_{qt}$. For all pairs $(q,t)$ with $q\ne t$ choose bijections of $I_{qt}$ onto $I_{tq}$ so that the bijections $I_{qt}\to I_{tq}$ and $I_{tq}\to I_{qt}$ be inverse to each other. Also choose involutive fixed point free permutations on the sets $I_{qq}$.

Combining these bijections we get a permutation $i\mapsto i'$ of the set $\{1,\ldots,n\}$ which determines a partition into pairs: $\{i_1i_1'\},\ldots,\{i_si_s'\}$. Using the fact that $n_{k\ell}\ne 0$, pick any $i_0\in I_{k\ell}$. A form of shape \eqref{(5.1)} has been thus constructed.

Forms \eqref{(5.1)} and \eqref{(5.2)} have identical collections of invariants in the case where $\,m_{i_0}=m_{j_0}$, $\,m_{i_0'}=m_{j_0''}\,$, and the sets $I_{qt}=\{i\mid m_i=q,\ m_{i'}=t\}$ and $J_{qt}=\{i\mid m_i=q,\ m_{i''}=t\}$ are of the same cardinality for all $q,t$ (here $i\mapsto i''$ is the involutive permutation of the set $\{1,\ldots,n\}$ corresponding to the partition into pairs $\{j_1j_1''\},\ldots,\{j_sj_s''\}$).  If these conditions are satisfied, then there exists a permutation of indices $1,\ldots,n$ which maps $I_{qt}$ onto $J_{qt}$, sends $i_0$ to $j_0$, and, moreover, whenever some $i$ goes to $j$, the paired index $i'$ goes to $j''$, that is, the pairs $\{j_1j_1''\},\ldots,\{j_sj_s''\}$ in a certain order correspond to the pairs $\{i_1i_1'\},\ldots,\{i_si_s'\}$.\qed 

\begin{theorem}\label{5.2}{\normalfont($p>2$).}
Any contact form \repl{corresponding} to the flag $\Fcal$ is conjugate under $G(\Fcal)$ to a form expressed as
\begin{equation}
\label{(5.3)}
f\mskip1mu(dx_{i_0}+x_{i_1} dx_{i_1'}+\ldots+x_{i_s}dx_{i_s'}),
\end{equation}
where $f\in\Ocal(\Fcal)^\times$, $i_0\in \{1,\ldots,n\}$ is an index, $\{i_1i_1'\},\ldots,\{i_si_s'\}$ are a partition of the set $\{1,\ldots,n\}\diagdown\{i_0\}$ into pairs, $s=(n-1)/2$. In order that the form \eqref{(5.3)} be conjugate for some $g\in\Ocal(\Fcal)^\times$ to the form
\begin{equation}
\label{(5.4)}
g\mskip1mu(dx_{j_0}+x_{j_1} dx_{j_1''}+\ldots+x_{j_s} dx_{j_s''}),
\end{equation}
where $j_0\mskip1mu,\{j_1j_1''\},\ldots,\{j_sj_s''\}$ have the same meaning\footnote{I.e., $j_0\in \{1,\ldots,n\}$, and the collection $\{j_1j_1''\},\ldots,\{j_sj_s''\}$ is a partition of the set $\{1,\ldots,n\}\diagdown\{j_0\}$ into pairs.}, it is necessary and sufficient that there exist a permutation of indices which sends $i_0$ to $j_0,$ every pair in the partition $\{i_1i_1'\},\ldots,\{i_si_s'\}$ to a pair occurring in the second partition $\{j_1j_1''\},\ldots,\{j_sj_s''\},$ and maps each subset of indices corresponding to indeterminates of the same height into itself.
\end{theorem}

\noindent
\textbf{Proof.}
Let us pass to $V=E^*$ in the same way as in the proof of Theorem~\ref{5.1}. By doing this the set in the right-hand side of \eqref{(3.6)} is interpreted as the set described in Lemma~\ref{4.7}. Since~\eqref{(3.6)} is a bijection, the orbits of the contact forms under the action of the group $\OFrtmGF$ are in a one-to-one correspondence with the collections of invariants $k,n_{qt}$ subject to constraints \eqref{(4.8)} with $r=m$.

Let us compute them for the form \eqref{(5.3)}. Let $i\mapsto i'$ be the involutive permutation of the set $\{1,\ldots,n\}\diagdown\{i_0\}$ defined by the partition into pairs $\{i_1i_1'\},\ldots,\{i_si_s'\}$. The mapping \eqref{(3.4)} sends the form \eqref{(5.3)} to 
$$
(x_{i_0},x_{i_1}\wedge x_{i_1'}+\ldots+x_{i_s}\wedge x_{i_s'})\in\MaE.
$$
The vectors $\{v_i\mid i\ne i_0\}$ constitute a basis of the space $Q=\Ker x_{i_0}\subseteq V$. The vectors $\{v_i\mid m_i\leqslant q\}$ constitute a basis of $V_q$ and the vectors $\{v_i\mid i\ne i_0,\ m_{i'}\geqslant t\}$ constitute a basis of $Q\cap (Q\cap V_{t-1})^\perp$. The classes of vectors $\{v_i\mid i\ne i_0,\ m_i=q,\ m_{i'}=t\}$ constitute a basis of
$$
Q\cap V_q\cap(Q\cap V_{t-1}^\perp)/\bigl(Q\cap V_q\cap(Q\cap V_t)^\perp+Q\cap V_{q-1}\cap(Q\cap V_{t-1})^\perp\bigr).
$$
Thanks to formulas \eqref{(4.9)}, \eqref{(4.10)} we have
$$
k=m_{i_0},\quad n_{qt}=\card\{i\ne i_0\mid m_i=q,\ m_{i'}=t\}.
$$
The remaining part of the proof is similar  to that of Theorem 5.1.\qed  

\vfil\eject

\section{Preparation to the classification of 1st-type \repl{symplectic} forms}\label{6}

Let the field $K$   be perfect. In a vector space $V$ consider flags of subspaces $V(k)$ of the form
\begin{gather}
\label{(6.1)}
0=V(0)\subseteq V(1)\subseteq\ldots,\qquad{\textstyle\bigcup\limits_{k=0}^\infty}\,V(k)=V;\\
\label{(6.2)}
0=V(0)\subseteq V(1)\subseteq\ldots\subseteq V(\infty)={\textstyle\bigcup\limits_{k=0}^\infty}V(k)\subseteq V(\infty+1)=V;\\
\label{(6.3)}
0=V(\infty)\subseteq V(\infty+1)=V;\\
\label{(6.4)}
V=V(0)\supseteq V(1)\supseteq\ldots,\qquad{\textstyle\bigcap\limits_{k=0}^\infty}V(k)=0;\\
\label{(6.5)}
V=V(0)\supseteq V(1)\supseteq\ldots\supseteq V(\infty)={\textstyle\bigcap\limits_{k=0}^\infty V(k)}\supseteq V(\infty+1)=0;\\
\label{(6.6)}
V=V(\infty)\supseteq V(\infty+1)=0.
\end{gather}
We denote the sets of all flags on $V$ of the form \eqref{(6.1)}--\eqref{(6.6)}, respectively, by
$$
\FNup(V),\quad\FNinftyup(V),\quad\Finftyup(V),\quad\FNdown(V),\quad\FNinftydown (V),\quad\Finftydown (V).
$$
Let $\Gcal$ be one of the flags \eqref{(6.1)}--\eqref{(6.6)}. Its $k$-th factor defined as $V(k+1)/V(k)$ for increasing flags and $V(k)/V(k+1)$ for decreasing flags, where $k\in\{0\}\cup\Nee$ for the flags \eqref{(6.1)} and \eqref{(6.4)}, $\,k\in\{0,\infty\}\cup\Nee$ for the flags \eqref{(6.2)} and \eqref{(6.5)}, and $k=\infty$ for the flags \eqref{(6.3)} and \eqref{(6.6)}, will be denoted by $\Phi_k\Gcal$.

Define a category $\Acal=\Acal_0$. Its objects are tuples
\begin{equation}
\label{(6.7)}
A=(V,\Gcal,b,W,\Hcal,c,(\nu^k)_{k=0,1,2,\ldots}),
\end{equation}
where $V$, $W$ are finite-dimensional vector spaces, $\Gcal\in\FNup(V)$, $\Hcal\in\FNup(W)$, $b:V\times V\to K$ and $c:W\times W\to K$ are bilinear \repl{antisymmetric} forms, where $b$ is non-degenerate, and $\nu^k\!:\Phi_k\Gcal\to\Phi_k\Hcal$ are some $p^{k+1}$-semilinear isomorphisms. We will define morphisms in~$\Acal$ a bit later.

Our goal is to determine objects of the category $\Acal$ up to an equivalence. For this we have to consider more cumbrous objects.

Denote by $I_n'$ (resp. $I_n''$) the set of sequences $\iov=(i_1i_2,\ldots,i_n)$,  where $i_s\in\{0,\infty\}\cup\Nee$ and, if among the indices we encounter $\infty$ with $i_t$ being the first such index, then $t$ is even (resp. odd), and all the subsequent indices are also equal to $\infty$.

The set $I_0'=I_0''$ is supposed to consist of the single empty sequence of indices, while $I_\infty'$ and $I_\infty''$ are the sets of infinite sequences with the property just described. For any $\overline i\in I_n'\cup I_n''$, we always assume that $i_1$ is the first index of the sequence, $i_2$ is the second one, etc. If $\overline i\in I_n'\cup I_n''$ and $j\in\{0,\infty\}\cup\Nee$, then $\overline ij$ (resp. $j\overline i$) denotes the sequence
$$
(i_1,\ldots,i_nj)\quad\text{(resp. $(ji_1,\ldots,i_n)$)}.
$$

The objects of the category $\Acal_n$ for $n>0$ are, by definition, collections
\begin{equation}
\label{(6.8)}
A=\left((V\klpair,\,\Gcal\klpair,\,b\klpair)_{\klov\in I_n'},\,(W\klpair,\,\Hcal\klpair)_{\klov\in I_n'',\,k_1\ne\infty},\,(c\klpair)_{\klov\in I_n'',\,k_1,\ell_1\ne\infty},\,(\nu\klpair)_{\kov\in I_{n+1}',\,\lov\in I_n'}\right)
\end{equation}
where $V\klpair$ and $W\klpair$ are finite-dimensional vector spaces among which there are finitely many non-zero ones, $\Gcal\klpair$ (resp. $\Hcal\klpair$) is a flag of subspaces $V\klpair(q)$ (resp. $W\klpair(q)$) in $V\klpair$ (resp. $W\klpair$) such that
\begin{gather} 
\label{(6.9)}
\Gcal\klpair\in\begin{cases}
\FNup(V\klpair)&\text{for }n\text{ even, }k_n\ne\infty;\\
\noalign{\smallskip}
\Finftyup(V\klpair)&\text{for }n\text{ even, }k_n=\infty;\\
\noalign{\smallskip}
\FNinftydown(V\klpair)&\text{for }n\text{ odd, }k_n\ne\infty;\\
\noalign{\smallskip}
\Finftydown(V\klpair)&\text{for }n\text{ odd, }k_n=\infty;
\end{cases}\\
\label{(6.10)}
\Hcal\klpair\in\begin{cases}
\FNinftyup(W\klpair)&\text{for }n\text{ even, }k_n\ne\infty;\\
\noalign{\smallskip}
\Finftyup(W\klpair)&\text{for }n\text{ even, }k_n=\infty;\\
\noalign{\smallskip}
\FNdown(W\klpair)&\text{for }n\text{ odd, }k_n\ne\infty;\\
\noalign{\smallskip}
\Finftydown(W\klpair)&\text{for }n\text{ odd, }k_n=\infty;
\end{cases}\\
\label{(6.11)}
b\klpair:V\klpair\times V\lkpair\to K\quad\text{and}\quad c\klpair:W\klpair\times W\lkpair\to K
\end{gather}
are non-degenerate bilinear pairings such that
\begin{equation}
\label{(6.12)}
b\klpair(v,v')=-b\lkpair(v',v)\quad\text{and}\quad c\klpair(w,w')=-c\lkpair(w',w)
\end{equation}
for all $v\in V\klpair$, $v'\in V\lkpair\mskip1mu$, $w\in W\klpair$, $w'\in W\lkpair\mskip1mu$,
\begin{equation}
\label{(6.13)}
b\kkpair(v,v)=0\quad\text{and}\quad c\kkpair(w,w)=0\quad\text{for all }\,v\in V\kkpair\mskip1mu,\ \, w\in W\kkpair\mskip1mu,
\end{equation}
and
\begin{equation}
\label{(6.14)} 
\nu\klpair:\Phi_{k_{n+1}}\Gcal_{(\ell_1\ldots\ell_n)}^{(k_1\ldots k_n)}\to\Phi_{\ell_n}\Hcal_{(k_2\ldots k_nk_{n+1})}^{(k_1\ell_1\ldots\ell_{n-1})}
\end{equation}
are $p^{k_1+1}$-semilinear isomorphisms.

Observe slight differences in the definition of objects of the category $\Acal_0$ in conditions \eqref{(6.10)}, \eqref{(6.11)}.

A morphism of object \eqref{(6.8)} to an object 
\begin{equation}
\label{(6.15)}
\Tilde A=\bigl(\Tilde V\klpair,\,\Tilde\Gcal\klpair,\,\Tilde b\klpair,\,\Tilde W\klpair,\,\Tilde\Hcal\klpair,\,\Tilde c\klpair,\,\Tilde\nu\klpair\,\bigr)\in\Acal_n,
\end{equation}
for $n\geqslant 0$, is a collection
\begin{equation}
\label{(6.16)}
\bigl((\varphi\klpair)_{\klov\in I_n'},(\psi\klpair)_{\klov\in I_n'',\,k_1\ne\infty}\bigr)
\end{equation}
of isomorphisms of vector spaces $\varphi\klpair:V\klpair\to\Tilde V\klpair$ and $\psi\klpair:W\klpair\to\Tilde W\klpair$ such that
\begin{gather}
\label{(6.17)}
\varphi\klpair\,\Gcal\klpair=\Tilde\Gcal\klpair\quad\text{and}\quad\psi\klpair\,\Hcal\klpair=\Tilde\Hcal\klpair,\\
\label{(6.18)}
\Tilde b\klpair(\varphi\klpair v,\varphi\lkpair v')=b\klpair(v,v')\quad\text{and}\quad\Tilde c\klpair(\psi\klpair w,\psi\lkpair w')=c\klpair(w,w')
\end{gather}
for all $v\in V\klpair$, $v'\in V\lkpair$, $w\in W\klpair$, $w'\in W\lkpair$, and the diagrams 
\begin{equation}
\label{(6.19)}
\begin{CD}
\Phi_{k_{n+1}}\Gcal_{(\ell_1\ldots\ell_{n-1}\ell_n)}^{(k_1k_2\ldots k_n)}@>\nu\klpair >>\Phi_{\ell_n}\Hcal_{(k_2\ldots k_nk_{n+1})}^{(k_1\ell_1\ldots\ell_{n-1})}\\
@V\varphi_{(\ell_1\ldots\ell_{n-1}\ell_n)}^{(k_1k_2\ldots k_n)}VV@V\psi_{(k_2\ldots k_nk_{n+1})}^{(k_1\ell_1\ldots\ell_{n-1})}VV\\
\Phi_{k_{n+1}}\Tilde\Gcal{}_{(\ell_1\ldots\ell_{n-1}\ell_n)}^{(k_1k_2\ldots k_n)}@>\Tilde\nu {}\klpair >>\Phi_{\ell_n}\Tilde\Hcal{}_{(k_2\ldots k_nk_{n+1})}^{(k_1\ell_1\ldots\ell_{n-1})}
\end{CD}
\end{equation}
where the vertical arrows are induced by the indicated mappings, are commutative.  

Observe that all morphisms of the categories considered in this section are equivalences. 

The \textit{orthogonal sum of objects} \eqref{(6.15)} and
\begin{equation}
\label{(6.20)}
\Tilde{\Tilde A}=\bigl(\Tilde{\Tilde V}\klpair,\,\Tilde{\Tilde\Gcal}\klpair,\,\Tilde{\Tilde b}\klpair,\,\Tilde{\Tilde W}\klpair,\Tilde{\Tilde\Hcal}\klpair,\,\Tilde{\Tilde c}\klpair,\,\Tilde{\Tilde\nu}\klpair\,\bigr)\in\Acal_n
\end{equation}
is an object \eqref{(6.8)} in which
$$
V\klpair=\Tilde V\klpair\oplus\Tilde{\Tilde V}\klpair,\quad W\klpair=\Tilde W\klpair\oplus\Tilde{\Tilde W}\klpair,
$$
flags $\Gcal\klpair$ and $\Hcal\klpair$ are defined by formulas $V\klpair(q)=\Tilde V{}\klpair(q)\oplus\Tilde{\Tilde V}\klpair(q)$ and  $W\klpair(q)=\Tilde W{}\klpair(q)\oplus\Tilde{\Tilde W}\klpair(q)$, respectively, the pairings $b\klpair$ and $c\klpair$ by the formulas 
$$
b\klpair(\Tilde v+\Tilde{\Tilde v},\Tilde v'+\Tilde{\Tilde v}')=\Tilde b\klpair(\Tilde v,\Tilde v')+\Tilde{\Tilde b}\klpair(\Tilde{\Tilde v},\Tilde{\Tilde v}')\quad\text{and}\quad c\klpair(\Tilde w+\Tilde{\Tilde w},\Tilde w'+\Tilde{\Tilde w}')=\Tilde c\klpair(\Tilde w,\Tilde w')+\Tilde{\Tilde c}\klpair(\Tilde{\Tilde w},\Tilde{\Tilde w}')
$$
for all $\Tilde v\in\Tilde V\klpair$, $\Tilde{\Tilde v}\in\Tilde{\Tilde V}\klpair$, $\Tilde v'\in\Tilde V\lkpair$, $\Tilde{\Tilde v}'\in\Tilde{\Tilde V}\lkpair$, $\Tilde w\in\Tilde W\klpair$, $\Tilde{\Tilde w}\in\Tilde{\Tilde W}\klpair$, $\Tilde w'\in\Tilde W\lkpair$, $\Tilde{\Tilde w}'\in\Tilde{\Tilde W}\lkpair$,
and the isomorphisms $\nu\klpair$ are defined from the commutative diagrams
$$
\begin{array}{r@{\mskip\medmuskip}c@{\mskip\medmuskip}l}
\Phi_{k_{n+1}}\Gcal_{(\ell_1\ldots\ell_{n-1}\ell_n)}^{(k_1k_2\ldots k_n)}&\xrightarrow{\nu\klpair }&\Phi_{\ell_n}\Hcal_{(k_2\ldots k_nk_{n+1})}^{(k_1\ell_1\ldots\ell_{n-1})}\\
\rotatebox{45}{\raisebox{-10pt}[0pt][0pt]{$\cong$}}\qquad\mbox{}&&\qquad\rotatebox{-45}{$\cong$}\\
\Phi_{k_{n+1}}\Tilde\Gcal_{(\ell_1\ldots\ell_{n-1}\ell_n)}^{(k_1k_2\ldots k_n)}\oplus\Phi_{k_{n+1}}\Tilde{\Tilde\Gcal}{}_{(\ell_1\ldots\ell_{n-1}\ell_n)}^{(k_1k_2\ldots k_n)}&\xrightarrow{\Tilde\nu\klpair\,\oplus\,\Tilde{\Tilde\nu}\klpair}&\Phi_{\ell_n}\Tilde\Hcal{}_{(k_2\ldots k_nk_{n+1})}^{(k_1\ell_1\ldots\ell_{n-1})}\oplus \Phi_{\ell_n}\Tilde{\Tilde\Hcal}{}_{(k_2\ldots k_nk_{n+1})}^{(k_1\ell_1\ldots\ell_{n-1})}
\end{array}
$$
An object will be called \textit{decomposable} if it is equivalent to the orthogonal sum of some objects $\Tilde A$ and $\Tilde{\Tilde A}$ in which among the spaces $\Tilde V\klpair$ and $\Tilde W\klpair$, and also among the spaces $\Tilde{\Tilde V}\klpair$ and $\Tilde{\Tilde W}\klpair$, there exist non-zero ones. Otherwise the object will be called \textit{indecomposable}.

Define categories $\Bcal_m$, where $m=0,1,\ldots$, whose objects are collections
\begin{multline}
\label{(6.21)}
B=\Bigl(\bigl(R\ijpair,\,\Kcal\ijpair\bigr)_{\iov\in I_{m+1}',\,\jov\in I_m'},\,\bigl(d\ijpair\bigr)_{\ijov\in I_{m+1}'},\,\bigl(S\ijpair,\,\Lcal\ijpair\bigr)_{\iov\in I_{m+1}'',\,i_1\ne\infty,\,\jov\in I_m''},\\
\bigl(e\ijpair\bigr)_{\ijov\in I_{m+1}'',\,i_1,j_1\ne\infty},\,\bigl(\mu\ijpair\bigr)_{\iov\in I_{m+1}',\,\jov\in I_m'}\Bigr),
\end{multline}
where $R\ijpair$ and $S\ijpair$ are finite-dimensional vector spaces among 
which there are finitely many non-zero ones, $\Kcal\ijpair$ (resp. 
$\Lcal\ijpair$) are flags of subspaces $R\ijpair(q)$ (resp. $S\ijpair(q)$) 
in $R\ijpair$ (resp. $S\ijpair$) such that
\begin{gather}
\label{(6.22)}
\Kcal\ijpair\in\begin{cases}
\FNdown(R\ijpair)&\text{ for }m\text{ even, }j_m\ne\infty;\\
\noalign{\smallskip}
\Finftydown(R\ijpair)&\text{ for }m\text{ even, }j_m=\infty;\\
\noalign{\smallskip}
\FNinftyup(R\ijpair)&\text{ for }m\text{ odd, }j_m\ne\infty;\\
\noalign{\smallskip}
\Finftyup(R\ijpair)&\text{ for }m\text{ odd, }j_m=\infty;\\
\end{cases}\\
\label{(6.23)}
\Lcal\ijpair\in\begin{cases}
\FNinftydown(S\ijpair)&\text{ for }m\text{ even, }j_m\ne\infty;\\
\noalign{\smallskip}
\Finftydown(S\ijpair)&\text{ for }m\text{ even, }j_m=\infty;\\
\noalign{\smallskip}
\FNup(S\ijpair)&\text{ for }m\text{ odd, }j_m\ne\infty;\\
\noalign{\smallskip}
\Finftyup(S\ijpair)&\text{ for }m\text{ odd, }j_m=\infty;\\
\end{cases}
\end{gather}
\begin{equation}
\label{(6.24)}
\begin{aligned}
d\ijpair:{}&\Phi_{j_{m+1}}\Kcal_{(j_1\ldots j_m)}^{(i_1\ldots i_mi_{m+1})}\times \Phi_{i_{m+1}}\Kcal_{(i_1\ldots i_m)}^{(j_1\ldots j_mj_{m+1})}\to K\quad\text{and}\\
\noalign{\smallskip}
e\ijpair:{}&\Phi_{j_{m+1}}\Lcal_{(j_1\ldots j_m)}^{(i_1\ldots i_mi_{m+1})}\times \Phi_{i_{m+1}}\Lcal_{(i_1\ldots i_m)}^{(j_1\ldots j_mj_{m+1})}\to K
\end{aligned}
\end{equation}
non-degenerate bilinear pairings such that 
\begin{equation}
\label{(6.25)}
\arraycolsep1em
\left.\begin{array}{cc}
d\ijpair(y,y')=-d\jipair(y',y),&e\ijpair(z,z')=-e\jipair(z',z),\\
d\iipair(y,y)=0,&e\iipair(z,z)=0,
\end{array}
\right\}
\end{equation}
and
\begin{equation}
\label{(6.26)}
\mu\ijpair:R_{(j_1\ldots j_m)}^{(i_1\ldots i_mi_{m+1})}\to S_{(i_2\ldots i_{m+1})}^{(i_1j_1\ldots j_m)}
\end{equation}
are $p^{i_1+1}$-semilinear isomorphisms. A \textit{morphism} of the object \eqref{(6.21)} to another object
\begin{equation}
\label{(6.27)}
\Tilde B=\bigl(\Tilde R\ijpair,\,\Tilde{\Kcal}\ijpair,\,\Tilde d\mskip1mu\ijpair,\,\Tilde S\ijpair,\,\Tilde{\Lcal}\ijpair,\,\Tilde e\ijpair,\,\Tilde\mu\ijpair\,)\in\Bcal_m
\end{equation}
is a collection 
\begin{equation}
\label{(6.28)}
\bigl((f\ijpair)_{\iov\in I_{m+1}',\,\jov\in I_m'},\,(g\ijpair)_{\iov\in I_{m+1}'',\,\jov\in I_m'',\,i_1\ne\infty}\bigr)
\end{equation}
of isomorphisms of vector spaces $f\ijpair:R\ijpair\to\Tilde R\ijpair$ and $g\ijpair:S\ijpair\to\Tilde S\ijpair$ such that 
\begin{gather}
\label{(6.29)}
f\ijpair\,\Kcal\ijpair=\Tilde{\Kcal}\ijpair,\qquad g\ijpair\,\Lcal\ijpair=\Tilde{\Lcal}\ijpair,\\
\label{(6.30)}
\Tilde d\mskip1mu\ijpair(\varphi\ijpair y,\varphi\jipair y')=d\ijpair(y,y'),\qquad\Tilde e\ijpair(\psi\ijpair z,\psi\jipair z')=e\ijpair(z,z'),
\end{gather}
where the mapping $\varphi\ijpair:\Phi_{j_{m+1}}\Kcal_{(j_1\ldots j_m)}^{(i_1\ldots i_{m+1})}\to\Phi_{j_{m+1}}\Tilde{\Kcal}{}_{(j_1\ldots j_n)}^{(i_1\ldots i_{m+1})}$ for $\ijov\in I_{m+1}'$ is induced by $f_{(j_1\ldots j_m)}^{(i_1\ldots i_mi_{m+1})}$, the mapping $\psi\ijpair:\Phi_{j_{m+1}}\Lcal_{(j_1\ldots j_m)}^{(i_1\ldots i_mi_{m+1})}\to\Phi_{j_{m+1}}\Tilde{\Lcal}{}_{(j_1\ldots j_m)}^{(i_1\ldots i_mi_{m+1})}$ for $\ijov\in I_{m+1}''$ with $i_1,j_1\ne\infty$ is induced by $g_{(j_1\ldots j_m)}^{(i_1\ldots i_mi_{m+1})}$. The last condition in the definition of a morphism is that the diagrams 
\begin{equation}
\label{(6.31)}
\begin{CD}
R_{(j_1\ldots j_m)}^{(i_1i_2\ldots i_{m+1})}@>\mu\ijpair >>S_{(i_2\ldots i_{m+1})}^{(i_1j_1\ldots j_m)}\\
@Vf\ijpair VV@VVg_{(i_2\ldots i_{m+1})}^{(i_1j_1\ldots j_m)}V\\
\Tilde R{}_{(j_1\ldots j_m)}^{(i_1i_2\ldots i_{m+1})}@>\Tilde \mu {}\ijpair >>\Tilde S{}_{(i_2\ldots i_{m+1})}^{(i_1j_1\ldots j_m)}
\end{CD}
\end{equation}
must be commutative.

The notions of orthogonal sums of objects, decomposable and indecomposable objects are defined by analogy with the categories $\Acal_n$.

\medskip
We will define the \textit{grinding functors}
$$
\Gfrak_{\Bcal_n}^{\Acal_n}:\Acal_n\leadsto\Bcal_n\quad\text{and}\quad\Gfrak_{\Acal_{m+1}}^{\Bcal_m}:\Bcal_m \leadsto\Acal_{m+1}.
$$

First we will make several preliminary remarks. Let $V$ be a finite-dimensional vector space, $\Gcal=\bigl(V(k)\bigr)$ and $\Gcal'=\bigl(V[\ell]\mskip1mu\bigr)$ flags of its subspaces. In each quotient $\Phi_k\Gcal$ define a flag denoted by $\Gcal'_{\Phi_k\Gcal}$ whose $\ell$-th space coincides with the image of $V[\ell]\cap V(k+1)$ in $\Phi_k\Gcal$ if the flag $\Gcal$ is increasing and with the image of $V[\ell]\cap V(k)$ if the flag $\Gcal$ is decreasing. 

Observe that if $\Gcal'\in\FNup(V)$, resp. $\Finftyup(V),\FNinftyup(V),\ldots\ $, then $\Gcal'_{\Phi_k\Gcal}\in\FNup(\Phi_k\Gcal)$, resp. $\Finftyup(\Phi_k\Gcal)$, $\FNinftyup(\Phi_k\Gcal),\ldots\ $. 

Suppose there is given another finite-dimensional vector space $U$, a flag $\Fcal$ in it, and a bilinear pairing $b:V\times U\to K$. Consider the flag $\Fcal^{\perp_b}$ whose spaces coincide with the orthogonals of the spaces of $\Fcal$ with the respective indices.

We will say that the flag $(\Fcal^{\perp_b})_{\Phi_k\Gcal}$ is \textit{obtained by transferring the flag $\Fcal$ by means of the pairing $b$ to the factor} $\Phi_k\Gcal$. The flags $(\Gcal^{\perp_b})_{\Phi_\ell\Fcal}$ are similarly defined.

The pairing $b$ induces non-degenerate pairings
\begin{equation}
\label{(6.32)}
\Phi_\ell((\Fcal^{\perp_b})_{\Phi_k\Gcal})\times\Phi_k((\Gcal^{\perp_b})_{\Phi_{\ell\Fcal}})\to K 
\end{equation}
by the method described in \ref4, see \eqref{(4.6)}. If the pairing $b$ is non-degenerate and $\Fcal\in\FNup(U)$ or $\FNinftyup(U)$, etc., then $\Fcal^{\perp_b}\in\FNdown(V),\FNinftydown(V)$, etc., respectively, and the flag $(\Fcal^{\perp_b})_{\Phi_k\Gcal}$ is of the same type.

Suppose we are given not a pairing, but a $p^t$-semilinear isomorphism $\mu:V\to U$, where $t$ is an integer. If $\Fcal$ is of one of the types \eqref{(6.1)}--\eqref{(6.6)}, then the flags $\mu^{-1}\Fcal$ and $(\mu^{-1}\Fcal)_{\Phi_k\Gcal}$ have the same type. The same is true for the flags $\mu\Gcal$ and $(\mu\Gcal)_{\Phi_\ell\Fcal}$. We will say that the flag $(\mu^{-1}\Fcal)_{\Phi_k\Gcal}$ (resp. $(\mu\Gcal)_{\Phi_\ell\Fcal}$) \textit{is obtained by transferring the flag} $\Fcal$ (resp. $\Gcal$) \textit{by means of the isomorphism $\mu$ to the factor} $\Phi_k\Gcal$ (resp. $\Phi_\ell\Fcal$). The mapping $\mu$ induces a $p^t$-semilinear isomorphism
\begin{equation}
\label{(6.33)}
\Phi_\ell\bigl((\mu^{-1}\Fcal)_{\Phi_k\Gcal}\bigr)\to\Phi_k\bigl((\mu\Gcal)_{\Phi_\ell\Fcal}\bigr).
\end{equation}
If, for example, the flags $\Fcal=\bigl(U(\ell)\bigr)$ and $\Gcal=\bigl(V(k)\bigr)$ are increasing, then this isomorphism is expressed explicitly as the through mapping
\begin{multline*}
\bigl(\mu^{-1}U(\ell+1)\cap V(k+1)+V(k)\bigr)/\bigl(\mu^{-1}U(\ell)\cap V(k+1)+V(k)\bigr)\to\\
\to\bigl(U(\ell+1)\cap\mu V(k+1)+\mu V(k)\bigr)/\bigl(U(\ell)\cap\mu V(k+1)+\mu V(k)\bigr)\to\\
\to\bigl(U(\ell+1)\cap\mu V(k+1)+U(\ell)\bigr)/\bigl(U(\ell+1)\cap\mu V(k)+U(\ell)\bigr),
\end{multline*}
where the first arrow is induced by $\mu$ and the second one is the canonical isomorphism.

\medskip
Let us construct the functor $\Gfrak_{\Bcal_0}^{\Acal_0}$. Given an object \eqref{(6.7)}, we set $\Gfrak_{\Bcal_0}^{\Acal_0}A=B$, where $B$ is the object \eqref{(6.21)} for $m=0$, where 
$$
R^i=\Phi_i\Gcal=V(i+1)/V(i),\quad S^i=\Phi_i\Hcal=W(i+1)/W(i),\quad i=0,1,\ldots\ .
$$
The flag $\Kcal^i$ is obtained by transferring the flag $\Gcal$ by means of the form $b$ to the factor $R^i$, i.e.,
$$
R^i(k)=\bigl(V(i+1)\cap V(k)^{\perp_b}+V(i)\bigr)/V(i).
$$
Since $b$ is non-degenerate, we have $\Kcal^i\in\FNdown(R^i)$. The flag $\Lcal^i$ is obtained by transferring the flag $\Hcal$ by means of the form $c$ and adding the spaces $S^i(\infty)$, $S^i(\infty+1)$, i.e., 
\begin{gather*}
S^i(k)=\bigl(W(i+1)\cap W(k)^{\perp_c}+W(i)\bigr)/W(i)\quad\text{for}\quad k=0,1,\ldots\ ,\\
S^i(\infty)=\bigl(W(i+1)\cap W^{\perp_c}+W(i)\bigr)/W(i),\quad S^i(\infty+1)=0.
\end{gather*}
Thus,  $\Lcal^i\in\FNinftydown(S^i)$. The pairings 
$$
d_j^{\mskip2mu i}:\Phi_j\Kcal^i\times\Phi_i\Kcal^j\to K,\quad\text{where}\quad i,j\in\{0\}\cup\Nee,
$$
are induced by $b$ and the pairings $e_j^{\mskip1mu i}:\Phi_j\Lcal^i\times\Phi_i\Lcal^j\to K$ are induced by $c$ according to the construction 
\eqref{(6.32)}. Finally, $\,\mu^i=\nu^{\mskip2mu i}$.

\medskip
The functors $\Gfrak_{\Bcal_n}^{\Acal_n}$ for $n>0$ are similarly defined. For an object \eqref{(6.8)} we set $\Gfrak_{\Bcal_n}^{\Acal_n}A=B$, where $B$ is the object \eqref{(6.21)} for $m=n$, where
$$
R_{(j_1\ldots j_n)}^{(i_1\ldots i_ni_{n+1})}=\Phi_{i_{n+1}}\Gcal_{(j_1\ldots j_n)}^{(i_1\ldots i_n)}\,,\qquad S_{(j_1\ldots j_n)}^{(i_1\ldots i_ni_{n+1})}=\Phi_{i_{n+1}}\Hcal_{(j_1\ldots j_n)}^{(i_1\ldots i_n)}\,.
$$
The flag $\Kcal_{(j_1\ldots j_n)}^{(i_1\ldots i_ni_{n+1})}$ is obtained by transferring the flag $\Gcal_{(i_1\ldots i_n)}^{(j_1\ldots j_n)}$ by means of the pairing $b_{(j_1\ldots j_n)}^{(i_1\ldots i_n)}$ to the factor $R_{(j_1\ldots j_n)}^{(i_1\ldots i_ni_{n+1})}$. The flag $\Lcal_{(j_1\ldots j_n)}^{(i_1\ldots i_ni_{n+1})}$ for $j_1\ne\infty$ is obtained by transferring the flag $\Hcal_{(i_1\ldots i_n)}^{(j_1\ldots j_n)}$ by means of the pairing $c_{(j_1\ldots j_n)}^{(i_1\ldots i_n)}$ to the factor $S_{(j_1\ldots j_n)}^{(i_1\ldots i_ni_{n+1})}$, while for $j_1=\infty$ it is determined by its membership in $\Finftyup(S\ijpair)$ when $n$ is odd and in $\Finftydown(S\ijpair)$ when $n$ is even (conditions \eqref{(6.9)} and \eqref{(6.10)} imply \eqref{(6.22)} and \eqref{(6.23)}, respectively). The pairings $d_{(j_1\ldots j_nj_{n+1})}^{(i_1\ldots i_ni_{n+1})}$ and $e_{(j_1\ldots j_nj_{n+1})}^{(i_1\ldots i_ni_{n+1})}$ are obtained by applying construction \eqref{(6.32)} (conditions \eqref{(6.12)}, \eqref{(6.13)} imply \eqref{(6.25)}). Finally, $\,\mu_{(j_1\ldots j_n)}^{(i_1\ldots i_ni_{n+1})}=\nu_{(j_1\ldots j_n)}^{(i_1\ldots i_ni_{n+1})}$.

If \eqref{(6.16)} is a morphism of the object \eqref{(6.8)} to the object \eqref{(6.15)}, then $\varphi\klpair$ (resp. $\psi\klpair$) induces isomorphisms of the factors of the flag $\Gcal\klpair$ onto the factors of the flag $\Tilde\Gcal\klpair$ (resp. isomorphisms of the factors of the flag $\Hcal\klpair$ onto the factors of the flag $\Tilde\Hcal{}\klpair$), and the collection of these isomorphisms is a morphism $\Gfrak_{\Bcal_n}^{\Acal_n}A\to\Gfrak_{\Bcal_n}^{\Acal_n}\Tilde A$. This defines the functor $\Gfrak_{\Bcal_n}^{\Acal_n}$ on morphisms. 

\medskip
Define the functor $\Gfrak_{\Acal_{m+1}}^{\Bcal_m}$ setting $\Gfrak_{\Acal_{m+1}}^{\Bcal_m}B=A$ for the object \eqref{(6.21)}, where $A$ is defined as follows (see \eqref{(6.8)} for $n=m+1$):
$$
V_{(\ell_1\ldots\ell_m\ell_{m+1})}^{(k_1\ldots k_mk_{m+1})}=\Phi_{\ell_{m+1}}\Kcal_{(\ell_1\ldots\ell_m)}^{(k_1\ldots k_mk_{m+1})},\qquad W_{(\ell_1\ldots\ell_m)}^{(k_1\ldots k_mk_{m+1})}=\Phi_{\ell_{m+1}}\Lcal_{(\ell_1\ldots\ell_m)}^{(k_1\ldots k_mk_{m+1})}.
$$
The flag $\Gcal_{(\ell_1\ldots\ell_m\ell_{m+1})}^{(k_1\ldots k_mk_{m+1})}$ is obtained by transferring the flag $\Lcal_{(k_2\ldots k_{m+1})}^{(k_1\ell_1\ldots\ell_m)}$ by means of the isomorphism $\mu_{(\ell_1\ldots\ell_m)}^{(k_1\ldots k_mk_{m+1})}$ to the factor $V_{(\ell_1\ldots\ell_m\ell_{m+1})}^{(k_1\ldots k_mk_{m+1})}$. The flag $\Hcal_{(\ell_1\ldots\ell_m\ell_{m+1})}^{(k_1\ldots k_mk_{m+1})}$ is obtained by transferring the flag $\Kcal_{(k_2\ldots k_{m+1})}^{(k_1\ell_1\ldots\ell_m)}$ by means of the isomorphism $\mu_{(k_2\ldots k_{m+1})}^{(k_1\ell_1\ldots\ell_m)}$ to the factor $W_{(\ell_1\ldots\ell_m\ell_{m+1})}^{(k_1\ldots k_mk_{m+1})}$ (from \eqref{(6.22)}, \eqref{(6.23)} it follows that \eqref{(6.9)}, \eqref{(6.10)} hold). The isomorphisms $\nu_{(\ell_1\ldots\ell_{m+1})}^{(k_1\ldots k_{m+2})}$ are determined by the construction \eqref{(6.33)}. Finally,
$$
b_{(\ell_1\ldots\ell_{m+1})}^{(k_1\ldots k_{m+1})}=d_{(\ell_1\ldots\ell_{m+1})}^{(k_1\ldots k_{m+1})},\qquad c_{(\ell_1\ldots\ell_{m+1})}^{(k_1\ldots k_{m+1})}=e_{(\ell_1\ldots\ell_{m+1})}^{(k_1\ldots k_{m+1})}
$$
(from \eqref{(6.25)} it follows that \eqref{(6.12)}, \eqref{(6.13)} hold).

\begin{lemma}\label{6.1}
Objects $A,\Tilde A\in\Acal_n$ are equivalent if and only if so are\/ $\Gfrak_{\Bcal_n}^{\Acal_n}A$ and\/ $\Gfrak_{\Bcal_n}^{\Acal_n}\Tilde A$.
\end{lemma}

\noindent
\textbf{Proof.}
Let $A$ be the object \eqref{(6.8)}, let $\Tilde A$ be the object \eqref{(6.15)}, and let \eqref{(6.28)} be an equivalence $\Gfrak_{\Bcal_n}^{\Acal_n}A\to\Gfrak_{\Bcal_n}^{\Acal_n}\Tilde A$. For each pair $\klov\in I_n'$ choose arbitrarily isomorphisms $\sigma\klpair:V\klpair\to\Tilde V\klpair$ so that $\sigma\klpair$ sends $\Gcal\klpair$ to $\Tilde\Gcal\klpair$ and the induced mappings of factors $\Phi_i\Gcal\klpair\to\Phi_i\Tilde\Gcal\klpair$ coincide with $f_{\lov}^{\kov i}$. Similarly, for $\klov\in I_n''$, where $k_1\ne\infty$, choose isomorphisms $\tau\klpair:W\klpair\to\Tilde W\klpair$ so that $\tau\klpair$ transforms $\Hcal\klpair$ to $\Tilde\Hcal{}\klpair$ and the induced mappings of factors $\Phi_i\Hcal\klpair\to\Phi_i\Tilde\Hcal\klpair$ coincide with $g_{\lov}^{\kov i}$.

Let us transfer the structure of the object $\Tilde A$ by means of isomorphisms $\sigma\klpair$ and $\tau\klpair$ to the spaces $V\klpair$ and $W\klpair$. We obtain an object, say \eqref{(6.20)}, in which 
$$
\Tilde{\Tilde V}\klpair=V\klpair,\quad\Tilde{\Tilde W}\klpair=W\klpair,\quad\Tilde{\Tilde\Gcal}\klpair=\Gcal\klpair,\quad\Tilde{\Tilde\Hcal}\klpair=\Hcal\klpair.
$$
By construction, the collection $(\sigma\klpair,\tau\klpair)$ defines an equivalence $\Tilde{\Tilde A}\to\Tilde A$ and the induced equivalence $\Gfrak_{\Bcal_n}^{\Acal_n}\Tilde{\Tilde A}\to\Gfrak_{\Bcal_n}^{\Acal_n}\Tilde A$ is given by the same collection $(f\ijpair,g\ijpair)$ as the equivalence $\Gfrak_{\Bcal_n}^{\Acal_n}A\to\Gfrak_{\Bcal_n}^{\Acal_n}\Tilde A$. Therefore, the collection of identity mappings gives an equivalence $\Gfrak_{\Bcal_n}^{\Acal_n}A\to\Gfrak_{\Bcal_n}^{\Acal_n}\Tilde{\Tilde A}$.

This implies, first, thanks to condition \eqref{(6.31)} and the definition of equivalence in the category~$\Bcal_n$ that $\nu\ijpair=\Tilde{\Tilde\nu}\ijpair$ for all $\iov\in I_{n+1}'$, $\jov\in I_n'$, i.e., $\Tilde{\Tilde A}$ differs from $A$ only by a collection of pairings.

Second, thanks to \eqref{(6.29)}, the flags $\Kcal_{\lov}^{\kov i}$ and $\Tilde{\Tilde{\Kcal}}_{\lov}^{\kov i}$ in $\Phi_i\Gcal\klpair$ that enter the structure of the objects $\Gfrak_{\Bcal_n}^{\Acal_n}A$ and $\Gfrak_{\Bcal_n}^{\Acal_n}\Tilde{\Tilde A}$, respectively, coincide, and the same is true for the flags $\Lcal_{\lov}^{\kov i}$ and $\Tilde{\Tilde{\Lcal}}{}_{\lov}^{\kov i}$ in $\Phi_i\Hcal\klpair$.

Third, by \eqref{(6.30)} the pairings $b\klpair$ and $\Tilde{\Tilde b}\klpair$ (resp. $c\klpair$ and $\Tilde{\Tilde c}\klpair$) induce the same pairings 
$$
\Phi_j\Kcal_{\lov}^{\kov i}\times\Phi_i\Kcal_{\lov}^{\kov j}\to K\quad\text{(resp. }\Phi_j\Lcal_{\lov}^{\kov i}\times\Phi_i\Lcal_{\lov}^{\kov j}\to K\text{)}.
$$
By the definition of the grinding functors this means that for any $\klov\in I_n'$ such that $\kov\ne\lov$ (resp. $\klov\in I_n''$ such that $k_1,\ell_1\ne\infty$ and $\kov\ne\lov$) the mapping \eqref{(4.7)} sends $b\klpair,\,\Tilde{\Tilde b}\mskip1mu\klpair\in B(V\klpair,V\lkpair)$ (resp. $c\klpair,\,\Tilde{\Tilde c}\mskip1mu\klpair\in B(W\klpair,W\lkpair)$) to the same point. By Proposition \ref{4.6} there exist $\varphi\klpair=N(\Gcal\klpair)$ and $\varphi\lkpair\in N(\Gcal\lkpair)$ (resp. $\psi\klpair\in N(\Hcal\klpair)$ and $\psi\lkpair\in N(\Hcal\lkpair)$) such that 
\begin{gather*}
\Tilde{\Tilde b}\mskip1mu\klpair(\varphi\klpair v,\varphi\lkpair v')=B\klpair(v,v'),\qquad\Tilde{\Tilde c}\mskip1mu\klpair(\psi\klpair w,\mskip1mu\psi\lkpair w')=c\klpair(w,w')\\
\text{for all}\quad v\in V\klpair,\quad v'\in V\lkpair,\quad w\in W\klpair,\quad w'\in W\lkpair.
\end{gather*}
For $\kov=\lov$ using Proposition \ref{4.2} instead of Proposition~\ref{4.6} we find $\varphi\kkpair$ and $\psi\kkpair$ with similar properties. The collection of mappings $(\varphi\klpair,\psi\klpair)$ determines an equivalence $A\to\Tilde{\Tilde A}$.

Indeed, conditions \eqref{(6.17)} and \eqref{(6.18)} in the definition of equivalence in $\Acal_n$ are satisfied by construction. Condition \eqref{(6.19)} holds since $\varphi\klpair$ (resp. $\psi\klpair$) induce the identity transformations of the factors $\Phi_i\Gcal\klpair$ (resp. $\Phi_i\Hcal\klpair$), Therefore, the objects $A$ and $\Tilde A$ are also equivalent. 

\begin{lemma}\label{6.2}Let $V,\Tilde V,W,\Tilde W$ be finite-dimensional vector spaces, $\mu:V\to W$ and $\Tilde\mu:\Tilde V\to\Tilde W$ semilinear isomorphisms,
$$
\Gcal=\bigl(V(k)\bigr),\quad\Tilde\Gcal=\bigl(\Tilde V(k)\bigr),\quad\Hcal=\bigl(W[\mskip1mu\ell\mskip2mu]\mskip1mu\bigr),\quad\Tilde\Hcal=\bigl(\Tilde W[\mskip1mu\ell\mskip2mu]\mskip1mu\bigr)
$$
flags of types \eqref{(6.1)}--\eqref{(6.6)} such that $\Gcal$ and $\Tilde\Gcal$ are of the same type, $\Hcal$ and $\Tilde\Hcal$ are also of the same type. Suppose that we are given collections of linear isomorphisms $\varphi_k:\Phi_k\Gcal\to\Phi_k\Tilde\Gcal$ and $\psi_\ell:\Phi_\ell\Hcal\to\Phi_\ell\Tilde\Hcal$ of all factors that make sense.

In order for isomorphisms $f:V\to\Tilde V$ and $g:W\to\Tilde W$ which transform, respectively, $\Gcal$ to $\Tilde\Gcal$ and $\Hcal$ to $\Tilde\Hcal$  to exist, induce isomorphisms $\varphi_k$ and $\psi_\ell$ of the corresponding factors, and make the diagram
$$
\begin{CD}
V@>\mu >>W\\
@VfVV@VVgV\\
\Tilde V@>\Tilde\mu>>\Tilde W
\end{CD}
$$
commutative, it is necessary and sufficient that every $\varphi_k$ transform $(\mu^{-1}\Hcal)_{\Phi_k\Gcal}$ to $(\Tilde\mu^{-1}\Tilde\Hcal)_{\Phi_k\Tilde\Gcal}\,,$ every $\psi_\ell$ transform $(\mu\Gcal)_{\Phi_\ell\Hcal}$ to $(\Tilde\mu\Tilde\Gcal)_{\Phi_\ell\Tilde\Hcal}\,,$ and all the diagrams
$$
\begin{CD}
\Phi_\ell\bigl((\mu^{-1}\Hcal)_{\Phi_k\Gcal}\bigr)@>>>\Phi_k\bigl((\mu\Gcal)_{\Phi_\ell\Hcal}\bigr)\\
@VVV@VVV\\
\Phi_\ell\bigl((\Tilde\mu^{-1}\Tilde\Hcal)_{\Phi_k\Tilde\Gcal}\bigr)@>>>\Phi_k\bigl((\Tilde\mu\Tilde\Gcal)_{\Phi_\ell\Tilde\Hcal}\bigr)
\end{CD}
$$
in which horizontal arrows are the mappings \eqref{(6.33)}, vertical arrows are induced, respectively, by $\varphi_k$ and $\psi_\ell$ be commutative.  
\end{lemma}

\noindent
\textbf{Proof.}
Everything boils down to the case where $W=V$, $\mu=id_V$, $\Tilde W=\Tilde V$, $\Tilde\mu=id_{\Tilde V}$. If an isomorphism $f:V\to\Tilde V$ sending $\,\Gcal$ to $\Tilde\Gcal$, $\,\Hcal$ to $\Tilde\Hcal\,$ and inducing isomorphisms $\varphi_k$ and $\psi_\ell$ exists, then, as is easy to see, $\varphi_k$ sends $\,\Hcal_{\Phi_k\Gcal}$ to $\Tilde\Hcal_{\Phi_k\Tilde\Gcal}\,$, $\,\psi_\ell$ sends $\,\Gcal_{\Phi_\ell\Hcal}$ to $\Tilde\Gcal_{\Phi_\ell\Tilde\Hcal}\,$ and the diagrams
\begin{equation}
\label{(6.34)}
\begin{array}{ccc}
\Phi_\ell(\Hcal_{\Phi_k\Gcal})&\cong &\Phi_k(\Gcal_{\Phi_\ell\Hcal})\\
\noalign{\smallskip}
\llap{$\varphi_k$}\Big\downarrow&&\Big\downarrow\rlap{$\psi_\ell$}\\
\noalign{\smallskip\smallskip}
\Phi_\ell(\Tilde\Hcal_{\Phi_k\Tilde\Gcal})&\cong &\Phi_k(\Tilde\Gcal_{\Phi_\ell\Tilde\Hcal})
\end{array}
\end{equation}
commute.

Conversely, let these properties of the collection $(\varphi_k,\psi_\ell)$ be satisfied. Choose a basis of $V$ compatible with both the flag $\Gcal$ and the flag $\Hcal$. Let $e$ be a basis element. Without loss of generality we may assume all the flags to be increasing. Suppose that 
$$
e\in V(k+1)\cap V[\ell+1],\quad e\notin V(k),\quad e\notin V[\mskip1mu\ell\mskip2mu].
$$
Let $\overline e$ be the image of $e$ in $\Phi_k\Gcal$ and $\overline{\overline e}$ the image of $e$ in $\Phi_\ell\Hcal$. Then
\begin{align*}
\overline u&=\varphi_k\mskip1mu\overline e\in\bigl(\Tilde V(k+1)\cap\Tilde V[\ell+1]+\Tilde V(k)\bigr)/\Tilde V(k),\\
\overline{\overline u}&=\psi_\ell\mskip2mu\overline{\overline e}\in\bigl(\Tilde V(k+1)\cap\Tilde V[\ell+1]+\Tilde V[\mskip1mu\ell\mskip2mu]\bigr)/\Tilde V[\mskip1mu\ell\mskip2mu],
\end{align*}
where thanks to \eqref{(6.34)} the image of $\overline u$ in $\Phi_\ell(\Tilde\Hcal_{\Phi_k\Tilde\Gcal})$ and the image of $\overline{\overline u}$ in $\Phi_k(\Tilde\Gcal_{\Phi_\ell\Tilde\Hcal})$ correspond to each other. Therefore, there exists $u\in\Tilde V(k+1)\cap\Tilde V[\ell+1]$ whose image in $\Phi_k\Gcal$ coincides with $\overline u$ and image in $\Phi_\ell\Hcal$ coincides with $\overline{\overline u}$. Assign this vector $u$ to the vector $e$. The linear mapping $f:V\to\Tilde V$ defined on the basis elements in this way has the desired properties.\qed 

\begin{lemma}\label{6.3}
Objects $B,\Tilde B\in\Bcal_m$ are equivalent if and only if so are\/ $\Gfrak_{\Acal_{m+1}}^{\Bcal_m}B$ and\/ $\Gfrak_{\Acal_{m+1}}^{\Bcal_m}\Tilde B$.
\end{lemma}

\noindent
\textbf{Proof.}
Let $B$ be object \eqref{(6.21)} and $\Tilde B$ object \eqref{(6.27)}; let \eqref{(6.16)} be an equivalence 
$$
\Gfrak_{\Acal_{m+1}}^{\Bcal_m}B\to\Gfrak_{\Acal_{m+1}}^{\Bcal_m}\Tilde B.
$$
By conditions \eqref{(6.17)}, \eqref{(6.19)} and the definition of an equivalence in $\Acal_{m+1}$ it follows that for every $\iov\in I_{m+1}'$ and $\jov\in I_m'$ we can apply Lemma~\ref{6.2} to the spaces $R\ijpair\,$, $\Tilde R\ijpair\,$, $S_{(i_2\ldots i_{m+1})}^{(i_1j_1\ldots j_m)}\,$, $\Tilde S_{(i_2\ldots i_{m+1})}^{(i_1j_1\ldots j_m)}\,$, the isomorphisms $\mu\ijpair\,$, $\Tilde\mu\ijpair\,$, the flags $\Kcal\ijpair\,$, $\Tilde{\Kcal}\ijpair\,$, $\Lcal_{(i_2\ldots i_{m+1})}^{(i_1j_1\ldots j_m)}\,$, $\Tilde{\Lcal}_{(i_2\ldots i_{m+1})}^{(i_1j_1\ldots j_m)}\,$, and the collection of isomorphisms of factors
$$
\varphi_{\jov k}^{\iov}:\Phi_k\Kcal\ijpair\to\Phi_k\Tilde{\Kcal}\ijpair\,,\qquad\psi_{(i_2\ldots i_{m+1}\ell)}^{(i_1j_1\ldots j_m)}:\Phi_\ell\Kcal_{(i_2\ldots i_{m+1}\ell)}^{i_1j_1\ldots j_m)}\to\Phi_\ell\Tilde{\Lcal}_{(i_2\ldots i_{m+1}\ell)}^{(i_1j_1\ldots j_m)}\,.
$$
Therefore, there exist isomorphisms
$$
f\ijpair:R\ijpair\to\Tilde R{}\ijpair,\qquad g_{(i_2\ldots i_{m+1})}^{(i_1j_1\ldots j_m)}:S_{(i_2\ldots i_{m+1})}^{(i_1j_1\ldots j_m)}\to\Tilde S_{(i_2\ldots i_{m+1})}^{(i_1j_1\ldots j_m)}
$$
which transform the flag $\Kcal\ijpair$ to $\Tilde{\Kcal}{}\ijpair$, the flag $\Lcal_{(i_2\ldots i_{m+1})}^{(i_1j_1\ldots j_m)}$ to $\Tilde{\Lcal}_{(i_2\ldots i_{m+1})}^{(i_1j_1\ldots j_m)}$, induce isomorphisms $\varphi_{\jov k}^{\iov}$ and $\psi_{(i_2\ldots i_{m+1}k)}^{(i_1j_1\ldots j_m)}$, and make the diagram \eqref{(6.31)} commutative. The obtained collection \eqref{(6.28)} gives an equivalence $B\to\Tilde B$ since condition \eqref{(6.30)} coincides precisely with condition \eqref{(6.18)} applied to the initial equivalence $\Gfrak_{\Acal_{m+1}}^{\Bcal_m}B\to\Gfrak_{\Acal_{m+1}}^{\Bcal_m}\Tilde B$.\qed 

\medskip
An object of $\Acal_n$ for $n>0$ (similarly, an object of $\Bcal_m$) will be called \textit{primitive} if all the flags entering its structure do not contain proper non-zero subspaces. We denote by $\Prim\Acal_n$ (resp. by $\Prim\Bcal_m$) the full subcategory of all primitive objects in $\Acal_n$ (resp. in $\Bcal_m$).

The grinding functors preserve primitive objects. If \eqref{(6.8)} is a primitive object, then its \textit{support} is the set 
\begin{gather*}
\supp A=\supp A_1\cup\supp_2A,\quad\text{where}\quad\supp_1A=\{(\kov,\overline l)\in I_n'\times I_n'\mid V\klpair\ne 0\},\\
\supp_2A=\{(\klov)\in I_n''\times I_n''\mid k_1\ne\infty,\,W\klpair\ne 0\}.
\end{gather*}
Similarly, the \textit{support} of a primitive object \eqref{(6.21)} is 
\begin{gather*}
\supp B=\supp_1B\cup\supp_2B,\quad\text{where}\quad\supp_1B=\{(\ijov)\in I_{m+1}'\times I_m'\mid R\ijpair\ne 0\},\\
\supp_2B=\{(\ijov)\in I''_{m+1}\times I_m''\mid i_1\ne\infty,\,S\ijpair\ne 0\}.
\end{gather*}
To study the primitive objects we introduce a new category $\Prim$. Its objects are collections 
\begin{equation}
\label{(6.35)}
C=\Bigl(\bigl(V\ijpair,b\ijpair\mskip1mu\bigr)_{\ijov\in I_\infty'},\,\bigl(W\ijpair\mskip1mu\bigr)_{\ijov\in I_\infty'',\,i_1\ne\infty},\,\bigl(c\ijpair\mskip1mu\bigr)_{\ijov\in I_\infty'',\,i_1,j_1\ne\infty},\,\bigl(\nu\ijpair\mskip1mu\bigr)_{\ijov\in I_\infty'}\Bigr),
\end{equation}
where $\,V\ijpair\,$, $\,W\ijpair\,$ are finite-dimensional vector spaces among which there are finitely many non-zero ones,
$$
b\ijpair:V\ijpair\times V\jipair\to K\quad\text{and }\quad c\ijpair:W\ijpair\times W\jipair\to K
$$
are non-degenerate bilinear pairings satisfying identities \eqref{(6.12)} and \eqref{(6.13)}, and
$$
\nu\ijpair:V_{(j_1j_ 2\ldots)}^{(i_1i_2i_3\ldots)}\to W_{(i_2i_3\ldots)}^{(i_1j_1j_2\ldots)}
$$
are $p^{i_1+1}$-semilinear isomorphisms. The \textit{support} of object \eqref{(6.35)} is the set 
\begin{gather*}
\supp C=\supp_1 C\cup\supp_2C,\quad\text{where}\quad\supp_1C=\{(\ijov)\in I_\infty'\times I_\infty'\mid V\ijpair\ne 0\},\\
\supp_2C=\{(\ijov)\in I''_\infty\times I''_\infty\mid i_1\ne\infty,\,W\ijpair\ne 0\}.
\end{gather*}
The morphisms (equivalences), as well as orthogonal sums, decomposable and indecomposable objects are defined in an obvious way.

\medskip
Now we will define functors $\,\Prim\Acal_n\leadsto\Prim$. For $m\geqslant n$, denote by $\Gfrak_{\Acal_m}^{\Acal_n}:\Acal_n\leadsto\Acal_m$ and $\Gfrak_{\Bcal_m}^{\Acal_n}:\Acal_n\leadsto\Bcal_m$ the composition of functors $\Gfrak_{\Bcal_r}^{\Acal_r}$ and $\Gfrak_{\Acal_{s+1}}^{\Bcal_s}$ ordered as needed. Suppose that \eqref{(6.8)} is a primitive object in $\Acal_n$. Consider the sequence of objects for $m\geqslant n$:
\begin{multline*}
A_m=\Gfrak_{\Acal_m}^{\Acal_n}A=\Bigl(\bigl(V\klpair,\,\Gcal\klpair,\,b\klpair\mskip1mu\bigr)_{\klov\in I_m'},\, \bigl(W\klpair,\,\Hcal\klpair\mskip1mu\bigr)_{\klov\in I_m'',\,k_1\ne\infty},\\
\bigl(c\klpair\mskip1mu\bigr)_{\klov\in I_m'',\,k_1,\ell_1\ne\infty},\,\bigl(\nu\klpair\mskip1mu\bigr)_{\kov\in I_{m+1}',\,\lov\in I_m'}\Bigr)\in\Acal_m\,,
\end{multline*}
\begin{multline*} 
B_m=\Gfrak_{\Bcal_m}^{\Acal_n}A=\Bigl(\bigl(V\klpair,\,\Gcal\klpair\mskip1mu\bigr)_{\kov\in I_{m+1}',\,\lov\in I_m'},\,\bigl(b\klpair\mskip1mu\bigr)_{\klov\in I_{m+1}'},\,\bigl(W\klpair,\,\Hcal\klpair\mskip1mu\bigr)_{\kov\in I_{m+1}'',\,k_1\ne\infty,\,\lov\in I_m''},\\
\bigl(c\klpair\mskip1mu\bigr)_{\klov\in I_{m+1}'',\,k_1,\ell_1\ne\infty},\bigl(\nu\klpair\mskip1mu\bigr)_{\kov\in I_{m+1}',\,\lov\in I_m'}\Bigr)\in\Bcal_m\,.
\end{multline*}
Since the object $A_m$ is primitive, it follows that any space of the flag $\Gcal\klpair$ is either zero, or coincides with $V\klpair$. Therefore, for any $(\klov)\in\supp_1A_m$, there exists a unique index $q$ such that $V_{\lov}^{\kov q}=\Phi_q\Gcal\klpair\ne 0$, i.e., $(\kov q,\lov )\in\supp_1B_m$, and for this index $V_{\lov}^{\kov q}\cong V\klpair$. Therefore, the mapping 
$$
((k_1\ldots k_mk_{m+1}),(\ell_1\ldots\ell_m))\mapsto((k_1\ldots k_m),(\ell_1\ldots\ell_m))
$$
gives a bijection $\supp_1B_m\to\supp_1A_m$ and, similarly, a bijection $\supp_2B_m\to\supp_2A_m$. Moreover, there are isomorphisms
\begin{gather*}
V_{(\ell_1\ldots\ell_m)}^{(k_1\ldots k_mk_{m+1})}\cong V_{(\ell_1\ldots\ell_m)}^{(k_1\ldots k_m)}\quad\text{for all }\,((k_1\ldots k_mk_{m+1}),(\ell_1\ldots\ell_m))\in\supp_1B_m,\\
\noalign{\smallskip}
W_{(\ell_1\ldots\ell_m)}^{(k_1\ldots k_mk_{m+1})}\cong W_{(\ell_1\ldots\ell_m)}^{(k_1\ldots k_m)}\quad\text{for all }\,((k_1\ldots k_mk_{m+1}),(\ell_1\ldots\ell_m))\in\supp_2B_m.
\end{gather*}
Similarly, the mapping
$$
((k_1\ldots k_mk_{m+1}),(\ell_1\ldots\ell_m\ell_{m+1}))\mapsto((k_1\ldots k_mk_{m+1}),(\ell_1\ldots\ell_m))
$$ 
gives bijections $\supp_1A_{m+1}\to\supp_1B_m$ and $\supp_2A_{m+1}\to\supp_2B_m\mskip1mu$, and there are canonical isomorphisms
\begin{gather*}
V_{(\ell_1\ldots\ell_m\ell_{m+1})}^{(k_1\ldots k_mk_{m+1})}\cong V_{(\ell_1\ldots\ell_m)}^{(k_1\ldots k_mk_{m+1})}\quad\text{for all }\,((k_1\ldots k_mk_{m+1}),(\ell_1\ldots\ell_m\ell_{m+1}))\in\supp_1A_{m+1},\\
\noalign{\smallskip}
W_{(\ell_1\ldots\ell_m\ell_{m+1})}^{(k_1\ldots k_mk_{m+1})}\cong W_{(\ell_1\ldots\ell_m)}^{(k_1\ldots k_mk_{m+1})}\quad\text{for all }\,((k_1\ldots k_mk_{m+1}),(\ell_1\ldots\ell_m\ell_{m+1}))\in\supp_2A_{m+1}. 
\end{gather*}
The obtained identifications of $V_{(\ell_1\ldots\ell_m)}^{(k_1\ldots k_m)}$ with $V_{(\ell_1\ldots\ell_m\ell_{m+1})}^{(k_1\ldots k_mk_{m+1})}$ and of $W_{(\ell_1\ldots\ell_m)}^{(k_1\ldots k_m)}$ with $W_{(\ell_1\ldots\ell_m\ell_{m+1})}^{(k_1\ldots k_mk_{m+1})}$ are compatible with the pairings, while the identifications of $V_{(\ell_1\ldots\ell_m)}^{(k_1k_2\ldots k_{m+1})}$ with $V_{(\ell_1\ldots\ell_m\ell_{m+1})}^{(k_1k_2\ldots k_{m+1}k_{m+2})}$ and of $W_{(k_2\ldots k_{m+1})}^{(k_1\ell_1\ldots\ell_m)}$ with $W_{(k_2\ldots k_{m+1}k_{m+2})}^{(k_1\ell_1\ldots\ell_m\ell_{m+1})}$ are compatible with the semilinear isomorphisms. Set 
\begin{align*}
Z_1&=\{(\ijov)\in I_\infty'\times I'_\infty\mid((i_1\ldots i_m),(j_1\ldots j_m))\in\supp_1A_m\quad\text{for all }\,m\geqslant n\},\\
\noalign{\smallskip}
Z_2&=\{(\ijov)\in I''_\infty\times I''_\infty\mid((i_1\ldots i_m),(j_1\ldots j_m))\in\supp_2A_m\quad\text{for all }\,m\geqslant n\}.
\end{align*}
To object \eqref{(6.8)} we assign object \eqref{(6.35)} in which 
\begin{equation}
\label{(6.36)}
V\ijpair=V_{(j_1\ldots j_n)}^{(i_1\ldots i_n)}\,,\qquad W\ijpair=W_{(j_1\ldots j_n)}^{(i_1\ldots i_n)}
\end{equation}
if $(\ijov)\in Z_1$ (resp. $(\ijov)\in Z_2$), and we set $V\ijpair=0$ (resp. $W\ijpair=0$) otherwise. If $(\ijov)\in Z_1$, then
$$
((j_1\ldots j_m),(i_1\ldots i_m))\in\supp_1A_m\quad\text{for all }\,m\geqslant n
$$ 
by the non-degeneracy of pairings $b_{(j_1\ldots j_m)}^{(i_1\ldots i_m)}$, whence $(\jov,\iov)\in Z_1$. Similarly, $(\ijov)\in Z_2$ and $j_1\ne\infty$ imply $(\jov,\iov)\in Z_2$. For each $(\ijov)\in Z_1$ (resp. $(\ijov)\in Z_2$ with $j_1\ne\infty$) we set 
$$
b\ijpair=b_{(j_1\ldots j_n)}^{(i_1\ldots i_n)}\,,\qquad c\ijpair=c_{(j_1\ldots j_n)}^{(i_1\ldots i_n)}\,.
$$
If $(\ijov)\in Z_1$, then 
$((i_1i_2\ldots i_{m+1}),(j_1\ldots j_m))\in\supp_1B_m\mskip1mu$, 
for all $m\geqslant n$, whence
$$
((i_1j_1\ldots j_m),(i_2\ldots i_{m+1}))\in\supp_2B_m
$$
by bijectivity of the mappings $\nu_{(j_1\ldots j_m)}^{(i_1i_2\ldots i_{m+1})}$. This implies $((i_1j_1\ldots j_{m-1}),(i_2i_3\ldots i_{m+1}))\in\supp_2A_m\mskip1mu$, i.e., $\,((i_1j_1j_2\ldots),(i_2i_3\ldots))\in Z_2\,$. Set
$$
\nu\ijpair=\nu_{(j_1\ldots j_n)}^{(i_1i_2\ldots i_{n+1})}.
$$
Thus,  we have bijections $\supp_1C\to\supp_1A$ and $\supp_2C\to\supp_2A$, as well as equalities \eqref{(6.36)} for all $(\ijov)\in\supp_1C$ (resp. $(\ijov)\in\supp_2C$). This entails immediately

\begin{lemma}\label{6.4}
Let \eqref{(6.8)} and \eqref{(6.15)} be two primitive objects in $\Acal_n$ to which objects $C$, $\Tilde C$ of\/ $\Prim$ correspond. Equivalences $A\to\Tilde A$ are in a bijective correspondence with equivalences $C\to\Tilde C$.
\end{lemma}

\medskip
We also have commutative diagrams
\begin{equation}
\label{(6.37)}
\arraycolsep2pt
\begin{array}{ccc}
\Prim\Acal_n&&\\
\llap{$\displaystyle\Gfrak_{\Acal_{n+1}}^{\Acal_n}\bigg\downarrow$}&\vcenter to40pt{\hbox{\rotatebox{-30}{$\longrightarrow$}}\vfil\hbox{\rotatebox{30}{$\longrightarrow$}}}&\Prim\\
\Prim\Acal_{n+1}&&
\end{array}
\end{equation}

Let $A$ be an object of $\Acal$. If $\Gfrak_{\Acal_m}^{\Acal}A$ is not primitive, then one of the spaces entering the structure of the object $\Gfrak_{\Acal_m}^{\Acal}A$ yields under grinding a factor of dimension strictly less than the dimension of the space itself. Therefore, for $n$ sufficiently large the object $\Gfrak_{\Acal_n}^{\Acal}A$ will be primitive.

Denote by $\Gfrak A$ the image of $\Gfrak_{\Acal_n}^{\Acal}A$ in $\Prim$. Since the diagram \eqref{(6.37)} is commutative, this image does not depend on $n$. Thus,  we have defined a functor
$$
\Gfrak:\Acal\leadsto\Prim.
$$

\begin{lemma}\label{6.5}
Two objects $A$, $\Tilde A\in\Acal$ are equivalent if and only if so are $\Gfrak A$ and $\Gfrak\Tilde A$.
\end{lemma}

This follows from Lemmas~\ref{6.1}, \ref{6.3}, \ref{6.4}.

\begin{lemma}\label{6.6}
Any object of the category\/ $\Prim$ is equivalent to $\Gfrak A$ for a suitable $A\in\Acal$.
\end{lemma}

\noindent
\textbf{Proof.}
Consider an object \eqref{(6.35)}. For $n\geqslant 0$, set
$$
V_{(\ell_1\ldots\ell_n)}^{(k_1\ldots k_m)}=\bigoplus_{{\scriptstyle\ijov\in I'_\infty\,\mid\,i_1=k_1,\ldots,\,i_m=k_m,\atop\scriptstyle\hphantom{\ijov\in I'_\infty}j_1=\ell_1,\ldots,\,j_n=\ell_n}}V\ijpair\,,\qquad W_{(\ell_1\ldots\ell_n)}^{(k_1\ldots k_m)}=\bigoplus_{{\scriptstyle\ijov\in I''_\infty\,\mid\,i_1=k_1,\ldots,\,i_m=k_m,\atop\scriptstyle\hphantom{\ijov\in I'_\infty}j_1=\ell_1,\ldots,\,j_n=\ell_n}}W\ijpair\,,
$$
where $m=n$ or $n+1$. Define the pairings
$$
b_{(\ell_1\ldots\ell_n)}^{(k_1\ldots k_n)}:V_{(\ell_1\ldots\ell_n)}^{(k_1\ldots k_n)}\times V_{(k_1\ldots k_n)}^{(\ell_1\ldots\ell_n)}\to K,\qquad c_{(\ell_1\ldots\ell_n)}^{(k_1\ldots k_n)}:W_{(\ell_1\ldots\ell_n)}^{(k_1\ldots k_n)}\times W_{(k_1\ldots k_n)}^{(\ell_1\ldots\ell_n)}\to K
$$
by making them on $V\ijpair\times V_{\overline{j'}}^{\overline{i'}}$ (resp. on $W\ijpair\times W_{\overline{j'}}^{\overline{i'}}$), where $\iov,\overline{i'},\jov,\overline{j'}\in I_\infty'$ (resp. $\iov,\overline{i'},\jov,\overline{j'}\in I_\infty''$ with $i_1,i_1'\ne\infty$) equal to $b\ijpair$ (resp. $c\ijpair\mskip1mu$) when $\,\overline{i'}=\jov$, $\,\overline{j'}=\iov\,$ and zero otherwise. 

Define $p^{k_1+1}$-semilinear isomorphisms 
$$
\nu_{(\ell_1\ldots\ell_n)}^{(k_1k_2\ldots k_{n+1})}:V_{(\ell_1\ldots\ell_n)}^{(k_1k_2\ldots k_{n+1})}\to W_{(k_2\ldots k_{n+1})}^{(k_1\ell_1\ldots\ell_n)}
$$
assuming them to be equal to $\nu\ijpair$ on $V\ijpair$ for $\ijov\in I_\infty'$.

Define the flags $\Gcal_{(\ell_1\ldots\ell_n)}^{(k_1\ldots k_n)}$ and $\Gcal_{(\ell_1\ldots\ell_n)}^{(k_1\ldots k_nk_{n+1})}$ as follows. The space $V_{(\ell_1\ldots\ell_n)}^{(k_1\ldots k_n)}(q)$ is the sum of all spaces $V\ijpair$ with $\ijov\in I_\infty'$ such that 
\begin{gather*}
i_1=k_1,\ldots,i_n=k_n,\qquad j_1=\ell_1,\ldots,j_n=\ell_n,\\
i_{n+1}<q\quad\text{for $n$ even},\qquad i_{n+1}\geqslant q\quad\text{for $n$ odd}.
\end{gather*}
The space $V_{(\ell_1\ldots\ell_n)}^{(k_1\ldots k_nk_{n+1})}(q)$ is the sum of all $V\ijpair$ with $\ijov\in I_\infty'$ such that
\begin{gather*}
i_1=k_1,\ldots,i_{n+1}=k_{n+1},\qquad j_1=\ell_1,\ldots,j_n=\ell_n,\\
j_{n+1}\geqslant q\quad\text{for $n$ even},\qquad j_{n+1}<q\quad\text{for $n$ odd}.
\end{gather*}
The values of $q$ are determined in each case by condition \eqref{(6.9)} or \eqref{(6.22)}.

The flags $\Hcal_{(\ell_1\ldots\ell_n)}^{(k_1\ldots k_n)}$ and $\Hcal_{(\ell_1\ldots\ell_n)}^{(k_1\ldots k_nk_{n+1})}$ are similarly defined. Then for $\ijov\in I_\infty'$ (resp. $\ijov\in I_\infty''$ with $i_1\ne\infty$) we have
\begin{alignat*}{2}
\Phi_{i_{n+1}}\Gcal_{(j_1\ldots j_n)}^{(i_1\ldots i_n)}&\cong 
V_{(j_1\ldots j_n)}^{(i_1\ldots i_ni_{n+1})},&\qquad\Phi_{j_{n+1}}\Gcal_{(j_1\ldots j_n)}^{(i_1\ldots i_ni_{n+1})}&\cong V_{(j_1\ldots j_nj_{n+1})}^{(i_1\ldots i_ni_{n+1})},\\
\noalign{\smallskip}
\Phi_{i_{n+1}}\Hcal_{(j_1\ldots j_n)}^{(i_1\ldots i_n)}&\cong 
W_{(j_1\ldots j_n)}^{(i_1\ldots i_ni_{n+1})},&\qquad\Phi_{j_{n+1}}\Hcal_{(j_1\ldots j_n)}^{(i_1\ldots i_ni_{n+1})}&\cong W_{(j_1\ldots j_nj_{n+1})}^{(i_1\ldots i_ni_{n+1})}.
\end{alignat*}
From the constructed ingredients we build up $A_n\in\Acal_n$ and $B_n\in\Bcal_n$ for each $n\geqslant 0$. Then $B_n\cong\Gfrak_{\Bcal_n}^{\Acal_n}A_n$ and $A_{n+1}\cong\Gfrak_{\Acal_{n+1}}^{\Bcal_n}B_n$. For $n$ sufficiently large, the objects $A_n$ are primitive and the corresponding objects of $\Prim$ are equivalent to $C$. Therefore, the objects $\Gfrak A_0$ and $C$ are equivalent.

\begin{lemma}\label{6.7}An object $A\in\Acal$ is decomposable if and only if\/ $\Gfrak A$ is decomposable.
\end{lemma}

\noindent
\textbf{Proof.}
If $A$ is decomposable, then all objects $\Gfrak_{\Acal_n}^{\Acal}A$ are decomposable, which implies decomposability of $\Gfrak A$.

Conversely, suppose $\Gfrak A$ is equivalent to the orthogonal sum $\Tilde C\oplus\Tilde{\Tilde C}$,  where $\Tilde C,\,\Tilde{\Tilde C}\ne 0$. By Lemma~\ref{6.6} $\,\Tilde C$ and $\Tilde{\Tilde C}$ are equivalent, respectively, to $\Gfrak\Tilde A$ and $\Gfrak\Tilde{\Tilde A}$ for some $\Tilde A,\,\Tilde{\Tilde A}\in\Acal$. Then $\Gfrak A$ is equivalent to $\Gfrak(\Tilde A\oplus\Tilde{\Tilde A})$. By Lemma~\ref{6.5} $A$ is equivalent to $\Tilde A\oplus\Tilde{\Tilde A}$, i.e., $A$ is decomposable.\qed 

\medskip
Set 
\begin{alignat*}{2}
\Qfrakov:&=I'_\infty\times I'_\infty,
&\qquad\Rfrakov:
&=\{(\ijov)\in I''_\infty\times I''_\infty\mid i_1\ne\infty\},\\
\Qfrakov_\sigma:&=\{(\ijov)\in\Qfrakov\mid i_2\ne\infty\},
&\qquad\Rfrakov_\tau:
&=\{(\ijov)\in\Rfrakov\mid j_1\ne\infty\}.
\end{alignat*}
Define\footnote{Note that $\tau$ and $\tau'$ are involutive permutations, $\rho$ is a bijection, and $\sigma$ an injection.}
$$
\tau:\Qfrakov\to\Qfrakov,\qquad\tau':\Rfrakov_\tau\to\Rfrakov_\tau,\qquad\rho:\Qfrakov\to\Rfrakov,\qquad\sigma:\Qfrakov_\sigma\to\Qfrakov
$$
by the formulas 
\begin{align*}
\tau\bigl((i_1i_2\ldots),(j_1j_2\ldots)\bigr)&=\bigl((j_1j_2\ldots),(i_1i_2\ldots)\bigr),\\
\tau'\bigl((i_1i_2\ldots),(j_1j_2\ldots)\bigr)&=\bigl((j_1j_2\ldots),(i_1i_2\ldots)\bigr),\\
\rho\bigl((i_1i_2i_3\ldots),(j_1j_2\ldots)\bigr)&=\bigl((i_1j_1j_2\ldots),(i_2i_3\ldots)\bigr),\\
\sigma\bigl((i_1i_2i_3\ldots),(j_1j_2j_3\ldots)\bigr)&=\bigl((i_3i_4i_5\ldots),(i_2i_1j_1j_2\ldots)\bigr).
\end{align*}
Thus, $\Qfrakov_\sigma=\rho^{-1}(\Rfrakov_\tau)$ and $\sigma=\tau\rho^{-1}\tau'\rho$.  Since $\tau^2=\id$ and $\tau^{\prime 2}=\id$, we have
\begin{equation}
\label{(6.38)}
\sigma\tau\sigma|_{\,\Qfrakov_\sigma}=\tau|_{\,\Qfrakov_\sigma}.
\end{equation}
Consider $\Qfrakov$ as a quiver with an arrow going from point $P\in\Qfrakov$ to point $Q\in\Qfrakov$ if and only if $P\in\Qfrakov_\sigma$ and $Q=\sigma P$. Formula \eqref{(6.38)} implies that $\tau$ is an anti-automorphism of the quiver $\Qfrakov$, i.e., $\tau$ reverses arrows. The connected components of the quiver $\Qfrakov$ are linear or cyclic directed chains without branching. Define $\alpha:\Qfrakov\to\Nee$ by the formula
$$
\alpha((i_1i_2\ldots\quad),(j_1j_2\ldots\quad))=i_1+1.
$$

Consider a category whose objects will be called \textit{symplectic representations} or, for short, simply \textit{representations of the quiver} $\Qfrakov$. They are described as collections
\begin{equation}
\label{(6.39)}
\bigl((V_P)_{P\in\Qfrakov},\,(b_P)_{P\in\Qfrakov},\,(h_P)_{P\in\Qfrakov_\sigma}\bigr),
\end{equation}
where the $V_P$ are finite-dimensional vector spaces among which there are finitely many non-zero ones,
$$
b_P:V_P\times V_{\tauP}\to K
$$
are non-degenerate pairings,
$$
h_P:V_P\to V_{\siP}
$$
are $p^{\alpha(P)-\alpha(\tau\sigma P)}$-semilinear isomorphisms such that 
\begin{gather}
\label{(6.40)}
b_P(u,v)=-b_{\tauP}(v,u)\quad\text{for all }\,u\in V_P,\ \,v\in V_{\tauP},\\
\noalign{\smallskip}
\label{(6.41)}
b_P(u,u)=0\quad\text{for all }\,u\in V_P,\quad\text{if }\,\tau P=P,\\
\noalign{\smallskip}
\label{(6.42)}
b_P(u,h_{\tausiP}\mskip1mu v)^{\,p^{\alpha(P)}}\!=b_{\siP}(h_P\mskip1mu u,v)^{\,p^{\alpha(\tau\sigma P)}}\quad\text{for all }\,P\in\Qfrakov_\sigma,\ u\in V_P,\ v\in V_{\tausiP},\\
\noalign{\smallskip}
\label{(6.43)}b_P(u,h_P\mskip1mu u)=0\quad\text{for all }\,u\in V_P,\quad
\text{if}\ \,P\in\Qfrakov_\sigma,\ \,\sigma P=\tau P.
\end{gather}
The \textit{support of the representation} \eqref{(6.39)} is the set $\{P\in\Qfrakov\mid V_P\ne 0\}$. The notions of equivalence of representations, orthogonal sums of representations, decomposable and indecomposable representations are defined in a natural way. 

\begin{lemma}\label{6.8}The category $\Prim$ is equivalent to the category of symplectic representations of the quiver $\Qfrakov$. The orthogonal sum decompositions in the two categories  correspond to each other.
\end{lemma}

\noindent
\textbf{Proof.} Given an object \eqref{(6.35)}, we write
\begin{gather*}
V_P=V\ijpair,\quad b_P=b\ijpair,\quad\nu_P=\nu\ijpair\quad\text{for }\,P=(\ijov)\in\Qfrakov,\\
W_Q=W\ijpair\quad\text{for }\,Q=(\ijov)\in\Rfrakov,\qquad\text{and}\qquad c_Q=c\ijpair\quad\text{for }\,Q=(\ijov)\in\Rfrakov_\tau.
\end{gather*}
For $P\in\Qfrakov_\sigma$, define $h_P$ as the composite of isomorphisms
$$
h_P:\,V_P\stackrel{\nu_P}{\displaystyle
\longrightarrow}W_{\rho P}\longrightarrow(W_{\tau'\!\rho P})^*=(W_{\rho\mskip1mu\tausiP})^*\longrightarrow(V_{\tausiP})^*\longrightarrow V_{\siP}.
$$
The second mapping here comes from the pairing $c_{\rho P}$. The penultimate mapping assigns to a linear form $\psi\in(W_{\rho\mskip1mu\tausiP})^*$ the linear form $\varphi\in(V_{\tausiP})^*$ given by the formula
$$
\varphi(v)=\psi\bigl(\nu_{\tausiP}(v)\bigr)^{\,p^{-\alpha(\tau\sigma P)}}.
$$
It is therefore $p^{-\alpha(\tau\sigma P)}$-semilinear. The last mapping comes from the pairing $b_{\siP}$. As a result, $h_P$ is $p^{\alpha(P)-\alpha(\tau\sigma P)}$-semilinear. The definition of $h_P$ amounts to the identity 
\begin{equation}
\label{(6.44)}
b_{\sigma P}(h_P\mskip1mu u,v)=c_{\rho P}(\nu_P\mskip1mu u,\nu_{\tausiP}\mskip1mu v)^{\,p^{-\alpha(\tau\sigma P)}}\quad\text{for all }\,u\in V_P,\ \,v\in V_{\tausiP}.
\end{equation}
Let us replace here $P$ with $\tau\sigma P$. Taking into account equality \eqref{(6.38)} we get
$$
b_{\tauP}(h_{\tausiP}\mskip1mu v,u)=c_{\tau'\!\rho P}(\nu_{\tausiP}\mskip1mu v,\nu_P\mskip1mu u)^{\,p^{-\alpha(P)}},\\
$$
which can be rewritten as
$$
b_P(u,h_{\tausiP}\mskip1mu v)=c_{\rho P}(\nu_P\mskip1mu u,\nu_{\tausiP}\mskip1mu v)^{\,p^{-\alpha(P)}}.
$$
Comparing this with \eqref{(6.44)} we get \eqref{(6.42)}. If $\sigma P=\tau P$, then \eqref{(6.44)} with $v=u$ yields \eqref{(6.43)}. Thus,  a representation \eqref{(6.39)} is constructed.

On the other hand, from the representation \eqref{(6.39)} we can recover the object \eqref{(6.35)} as follows. For $P\in\Qfrakov$, select arbitrarily a space $W_{\rho P}$ of the same dimension as $V_P$ and a $p^{\alpha(P)}$-semilinear isomorphism $\nu_P:V_P\to W_{\rho P}$. Now define pairings $c_{\rho P}$ by formula \eqref{(6.44)}.\qed 

\medskip
Denote by $\Qfrak$ the subset of pairs $(\ijov)\in\Qfrakov$ such that $\ijov$ simultaneously either contain $\infty$, or do not contain it and in this case are periodic.

\begin{lemma}\label{6.9}
The support of any representation \eqref{(6.39)} is contained in $\Qfrak$.
\end{lemma}

\noindent
\textbf{Proof.}
Since the set $\{P\in\Qfrakov\mid V_P\ne 0\}$ is finite, for each $Q=(\ijov)$ taken from this set there exists $n>0$ such that either $\sigma^nQ$ and $\sigma^{-n}Q$ are not defined, or $\sigma^nQ=Q$. By the definition of $\sigma$, in the first case $\iov$ and $\jov $ contain $\infty$, while in the second case they do not contain it and are periodic.\qed

\smallskip
A set of the form $X\cup\tau X$, where $X$ is a connected component, will 
be called a \textit{$\tau$-connected component of the quiver} $\Qfrakov$. 
Denote by $\Pfrak$ the set of pairs $(\klov)$, where $\klov$ are sequences of 
positive integers both either finite of the same length, or infinite periodic. 
Consider an equivalence relation on $\Pfrak$ generated by relations 
\begin{gather}
\label{(6.45)}
((k_1\ldots k_n),(\ell_1\ldots\ell_n))\sim((\ell_n\ldots\ell_1),(k_n\ldots k_1)),\\
\label{(6.46)}
((k_1k_2\ldots\quad),(\ell_1\ell_2\ldots\quad))\sim((k_2k_3\ldots\quad),(\ell_2\ell_3\ldots\quad)),\\
\label{(6.47)}
((k_1k_2\ldots\quad),(\ell_1\ell_2\ldots\quad))\sim((\ell_t\ell_{t-1}\ldots\ell_1\ell_t\ell_{t-1}\ldots\quad),(k_tk_{t-1}\ldots k_1k_t\ldots)),
\end{gather}
where $t$ is the common period of the two sequences.

\begin{lemma}\label{6.10}
The $\tau$-connected components of\/ $\Qfrakov$ contained in $\Qfrak$ are in a one-to-one correspondence with the points of the set $\Pfrak/\sim$.
\end{lemma}

\noindent
\textbf{Proof.} The components of interest to us are of the form
\begin{gather}
\notag
\mbox{}\\
\label{(6.48)}
\lefteqn{\raisebox{-220pt}[0pt][0pt]{\hspace*{112pt}\includegraphics{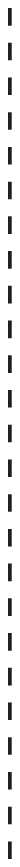}}}\raisebox{-10pt}{\includegraphics{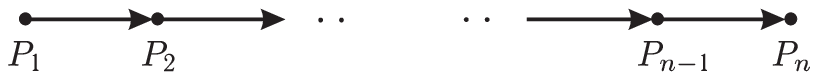}}\\[15pt]
\label{(6.49)}
\raisebox{-10pt}{\includegraphics{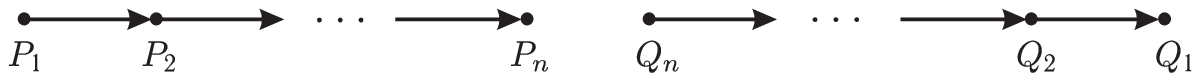}}\\[15pt]
\label{(6.50)}
\raisebox{-30pt}{\hspace*{-16pt}\includegraphics{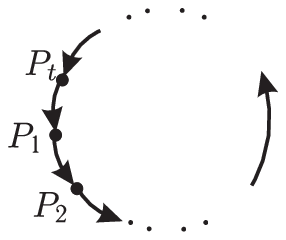}}\\[15pt]
\label{(6.51)}
\raisebox{-30pt}{\includegraphics{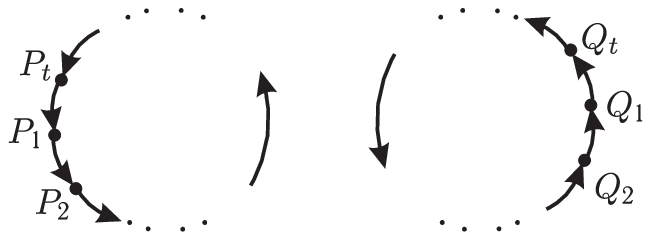}}\\
\notag
\mbox{}
\end{gather}
where the reflection in the dashed line corresponds to the action of the anti-automorphism $\tau$.

To components of the form \eqref{(6.48)}, \eqref{(6.49)} we assign the class of the pair
$$
\bigl((\alpha(P_1),\ldots,\alpha(P_n)),\ (\alpha(\tau P_1),\ldots,\alpha(\tau P_n))\bigr)
$$
in $\Pfrak/\sim$, and to components of the form \eqref{(6.50)}, \eqref{(6.51)} the class of
$$
\bigl((\alpha(P_1),\ldots,\alpha(P_t),\alpha(P_1),\ldots,\alpha(P_t),\ldots),\ (\alpha(\tau P_1),\ldots,\alpha(\tau P_t),\alpha(\tau P_1),\ldots,\alpha(\tau P_t),\ldots)\bigr).
$$
If 
$$
P_1=((i_1i_2\ldots\quad),(j_1j_2\ldots\quad)),
$$
then
\begin{align*}
\sigma^rP_1&=((i_{2r+1}i_{2r+2}\ldots\quad),(i_{2r}\ldots i_1j_1j_2\ldots\quad)),\\
\sigma^{-r}P_1&=((j_{2r}\ldots j_1i_1i_2\ldots\quad),(j_{2r+1}j_{2r+2}\ldots\quad)).
\end{align*}
This implies that the class in $\Pfrak/\sim$ represented by $((k_1\ldots k_n),(\ell_1\ldots\ell_n))\in\Pfrak$ corresponds to the unique $\tau$-connected component containing the point\footnote{Subtracting 1 from $k_1,\ell_1,\ldots$ is required to conform with the previously used indexing scheme,  where we have $\alpha(\ijov)=i_1+1$. This corrects the formulas for $P_1$ in the original manuscript.}
$$
P_1=\bigl((k_1-1,\ell_2-1,k_2-1,\ldots,\ell_n-1,k_n-1,\infty\ldots\quad),
\ (\ell_1-1,\infty\ldots\quad)\bigr),
$$
and the class represented by $((k_1k_2\ldots\quad),(\ell_1\ell_2\ldots\quad))\in\Pfrak$ corresponds to the component containing the point
\begin{multline*}
P_1=\bigl((k_1-1,\ell_2-1,k_2-1,\ldots,\ell_t-1,k_t-1,\ell_1-1,k_1-1,\ell_2-1,k_2-1,\ldots\quad),\\
(\ell_1-1,k_t-1,\ell_t-1,\ldots k_2-1,\ell_2-1,k_1-1,\ell_1-1,k_t-1,\ell_t-1,\ldots)\bigr),
\end{multline*}
where $t$ is the common period of the sequences $\klov$.\qed

\medskip
Clearly, any representation \eqref{(6.39)} decomposes uniquely as the orthogonal sum of representations with the support of each summand being a $\tau$-connected component. Let us examine one by one representations whose support is a component of some type.

First, consider the component \eqref{(6.49)}. Here a representation is determined, up to an equivalence, by the space $V_{P_1}$, and any decomposition of the space into a direct sum leads to the decomposition of the representation into an orthogonal sum. In this case, there exists a unique indecomposable representation corresponding to a one-dimensional space $V_{P_1}$.

Consider representations with support \eqref{(6.48)}. Define a non-degenerate bilinear form on $V_{P_1}$ setting
$$
(u,v)=b_{P_1}\bigl(u,h_{P_{n-1}}\circ\ldots\circ h_{P_1}(v)\bigr).
$$
This form is \repl{antisymmetric} since
$$
(u,u)=b_{P_s}\bigl(h_{P_{s-1}}\cdots h_{P_1}(u),h_{P_s}h_{P_{s-1}}\cdots h_{P_1}(u)\bigr)^{\,p^{\alpha(P_{n-1})+\ldots+\alpha(P_{s+1})-\alpha(P_1)-\ldots-\alpha(P_{s-1})}}=0
$$
by \eqref{(6.42)}, \eqref{(6.43)}, when $n=2s$ is even, and 
$$
(u,u)=b_{P_s}\bigl(h_{P_{s-1}}\cdots h_{P_1}(u),h_{P_{s-1}}\cdots h_{P_1}(u)\bigr)^{\,p^{\alpha(P_{n-1})+\ldots+\alpha(P_s)-\alpha(P_1)-\ldots-\alpha(P_{s-1})}}=0
$$
by \eqref{(6.42)}, \eqref{(6.41)}, when $n=2s-1$ is odd.

The non-degenerate symplectic space $V_{P_1}$ determines the representation up to an equivalence, and any decomposition of $V_{P_1}$ into a direct sum of orthogonal subspaces leads to the decomposition of the representation. We have again a unique indecomposable representation corresponding to a two-dimensional symplectic space.

Consider representations with support \eqref{(6.51)}. Let
$$
k_1=\alpha(P_1),\ldots,k_t=\alpha(P_t),k_{t+1}=\alpha(P_1),\ldots,\ell_1=\alpha(Q_1),\ldots,\ell_t=\alpha(Q_t),\ell_{t+1}=\alpha(Q_1),\ldots 
$$
So the component \eqref{(6.51)} corresponds to the class of $(\klov)$ in $\Pfrak/\sim$.

The composite mapping $h=h_{P_{\mskip 1mu t}}\circ h_{P_{\mskip 1mu t-1}}\circ\ldots\circ h_{P_{\mskip 1mu 1}}$ is a non-singular $p^{k_1+\ldots+k_t-\ell_1-\ldots-\ell_t}$-semilinear endomorphism of the space $V_{P_1}$.

The orthogonal decompositions of the representation correspond to the decompositions of $V_{P_1}$ into direct sums of $h$-invariant subspaces. Two representations are equivalent if and only if the corresponding pairs $(V_{P_1},h)$ and $(\Tilde V_{P_1},\Tilde h)$ are equivalent in the sense that there exists a linear isomorphism $\varphi:V_{P_1}\to\Tilde V_{P_1}$ such that $\Tilde h=\varphi h\varphi^{-1}$.

In \cite[Ch.~3]{3}, it is demonstrated that the theory of semilinear endomorphisms reduces to the representation theory of an associative ring. From the Krull-Remak-Schmidt theorem\footnote{This is one of basic results in the theory of rings and modules. More generally, the theorem can be formulated for groups with a set of endomorphisms (e.g., \cite[Ch.~1, Th.~11]{3}). Application of this theorem in our case is justified by the fact that only finite-dimensional representations are considered.} it follows therefore that the decomposition of the representation into an orthogonal sum of indecomposable representations is unique up to a permutation of summands and their equivalence.

Finally, let us study representations with support \eqref{(6.50)}. Let $\tau P_1=P_s$. On the space $V_{P_1}$ define a bilinear form 
$$
(u,v)=b_{P_1}\bigl(u,\ h_{P_{s-1}}\circ h_{P_{s-2}}\circ\ldots\circ h_{P_1}(v)\bigr)
$$
and a non-singular linear endomorphism $h=h_{P_t}\circ h_{P_{t-1}}\circ\ldots\circ h_{P_1}$. The form is non-degenerate thanks to non-degeneracy of the pairing $b_{P_1}$. In view of \eqref{(6.41)}, \eqref{(6.42)}, \eqref{(6.43)} it is also \repl{antisymmetric} and 
\begin{equation}
\label{(*)}
(u,hu)=0\quad\text{for all }\,u\in V_{P_1}.
\end{equation}
In this case everything is reduced to the study of the pairs $(V_{P_1},h)$, where $V_{P_1}$ is a non-degenerate symplectic space and $h$ is a non-singular linear endomorphism satisfying condition~\eqref{(*)}.

\begin{lemma}\label{6.11}
{\rm1)} Let $V$ be a non-degenerate symplectic space, $f$ a non-singular linear endomorphism of $V$ such that $(u,fu)=0$ for all $u\in V$. Then $V$ can be represented as a direct sum $U\oplus U'$ of two isotropic $f$-invariant subspaces.

\smallskip
{\rm2)} If $(\Tilde V,\Tilde f)$ is another such pair and $\Tilde V=\Tilde U\oplus\Tilde U'$ is a similar decomposition, then the pairs $(V,f)$ and $(\Tilde V,\Tilde f)$ are equivalent if and only if the endomorphisms $f|_U$ and $\Tilde f|_{\Tilde U}$ are equivalent {\normalfont(cf. Exercise 15b) in \S~5 of Chapter 9 in \cite{1}).}
\end{lemma}

\noindent
\textbf{Proof.} 1) We will prove by induction on $\dim V$. Let $U_1$ be an invariant indecomposable with respect to $f$ subspace of maximal dimension. Since $U_1$ is indecomposable, it is spanned by $\{v,fv,f^2 v,\ldots\}$ for a certain $v\in U_1$. In particular, $U_1$ is isotropic.

Further, there exists a smallest non-zero $f$-invariant subspace $U_1^{\min}\subseteq U_1$. Let $U_1'$ be an arbitrary invariant indecomposable with respect to $f$ subspace such that $(U_1^{\min},U_1')\ne 0$. Since $U_1^{\min}\not\subseteq U_1\cap(U_1')^\perp$, we have $U_1\cap U_1^{\prime\perp}=0$.

The pairing $U_1\times U_1'\to K$ has no left kernel, and since $\dim U_1'\leqslant\dim U_1$, it has no right kernel either. Thanks to indecomposability $U_1'$ is isotropic. Therefore, the restriction of the form to $U_1\oplus U_1'$ is non-degenerate. It remains to apply the induction hypothesis to the orthogonal complement of the subspace $U_1\oplus U_1'$.

2) Let $V=U\oplus U'$ be any decomposition with the required properties. Let us represent $U$ as a direct sum of indecomposable with respect to $f$ subspaces $U_1,\ldots,U_n$. Then $U'$ can be expressed as a direct sum of indecomposable with respect to $f$ subspaces $U_1',\ldots,U_n'$ so that $(U_i,U_j')=0$ for $i\ne j$, and the restriction of the form to $U_i\times U_i'$ is a non-degenerate pairing for each $i$. The minimal polynomials of the endomorphisms $f|_{U_i}$ and $f|_{U_i'}$ coincide thanks to the identity $(u,fv)=(fu,v)$. Therefore, $f|_{U_i}$ and $f|_{U_i'}$ are equivalent for each $i$. Hence, so are the endomorphisms $f|_{U}$ and $f|_{U'}$. 

If the pairs $(V,f)$ and $(\Tilde V,\Tilde f)$ are equivalent, then the endomorphisms $f$ and $\Tilde f$ are equivalent. By the above, the endomorphisms $f|_{U}$ and $\Tilde f|_{\Tilde{U}}$ are also equivalent. Conversely, if $f|_{U}$ and $\Tilde f|_{\Tilde{U}}$ are equivalent, then the identity $(u,fv)=(fu,v)$ implies the equivalence of $(V,f)$ and $(\Tilde V,\Tilde f)$.\qed 

\medskip
Thus, the equivalence classes of representations of the quiver $\Qfrakov$ with support \eqref{(6.50)} correspond to the equivalence classes of non-degenerate linear endomorphisms. In particular, the decomposition of such a representation into an orthogonal sum of indecomposables is unique.

\medskip
Summarizing what we have said above, we conclude that an arbitrary representation \eqref{(6.39)} decomposes into an orthogonal sum of indecomposables uniquely, up to a permutation of summands and their equivalence.

By Lemmas~\ref{6.5},~\ref{6.8} the equivalence classes of objects of $\Acal$ are in a one-to-one correspondence with the equivalence classes of representations of the quiver $\Qfrakov$. Moreover, it follows from Lemmas~\ref{6.7} and~\ref{6.8} that the orthogonal sum decompositions are respected under this correspondence. This implies uniqueness of the decomposition of any object of the category $\Acal$ as an orthogonal sum of indecomposable objects, up to a permutation of summands and their equivalence.

\medskip
Let us find indecomposable objects of the category $\Acal$ explicitly.

For any object \eqref{(6.7)}, there exist bases $e_1,\ldots,e_N$ in $V$ and $e_1',\ldots,e_N'$ in $W$ compatible with the flags $\Gcal$ and $\Hcal$, respectively, such that 
$$
m_i=\min\{q\mid e_i\in V(q)\}=\min\{q\mid e_i'\in W(q)\}
$$
and the mapping $\nu^{m_i-1}$ transforms the class $\overline e_i$ of the element $e_i$ in $\Phi_{m_i-1}\Gcal$ into the class $\overline{e_i'}$ of the element $e'_i$ in $\Phi_{m_i-1}\Hcal$. Such bases will be called \textit{compatible}.

The collection $m_1,\ldots,m_N$ carries complete information about the flags $\Gcal$ and $\Hcal$ and the isomorphisms $\nu^k$:
\begin{gather*}
V(q)=\langle e_i\mid m_i\leqslant q\rangle_K,\qquad W(q)=\langle e_i'\mid m_i\leqslant q\rangle_K,\\
\nu^k\overline e_i=\overline{e_i'}\quad\text{for all $i$ such that $\,m_i-1=k$}.
\end{gather*}
Thus,  any object \eqref{(6.7)} can be described by the collection $m_1,\ldots,m_N$ and the matrices of the forms $b$, $c$ in compatible bases. 

Let us show how to select compatible bases for the object corresponding to a representation \eqref{(6.39)}. We will use explicit constructions given in the proofs of Lemmas~\ref{6.6} and~\ref{6.8}. We have
$$ 
V={\textstyle\bigoplus\limits_{P\in\Qfrak}}\,V_P,\qquad W={\textstyle\bigoplus\limits_{P\in\Qfrak}}\,W_{\rho P}.  
$$
Choose bases in all spaces $V_P$ with $P\in\Qfrak$ so that for $\tau P\ne P$ the bases in $V_P$ and in $V_{\tau P}$ are dual to each other with respect to $b_P$, whereas for $\tau P=P$ the matrix of the form $b_P$ in the chosen basis of $V_P$ is in the normal form.

In $V$ consider the basis $e_1,\ldots,e_N$ obtained as the union of the bases of the spaces $V_P$. Then $m_i=\alpha(P)$ if $e_i\in V_P$. Let us order the basis of $V$ so that the matrix of the form $b$ were equal to (here $E$ is the unit matrix of the corresponding size) 
\begin{equation} 
\label{(6.52)}
\left(\begin{array}{c|c}
0&E\\ \hline -E&0\vbox to12pt{}
\end{array}
\right).
\end{equation}
In each space $W_{\rho P}$ with $P\in\Qfrak$ take the basis consisting of the vectors $\{e_i'=\nu_P\mskip1mu e_i\mid e_i\in V_P\}$. The union of these bases gives a basis $e_1',\ldots,e_N'$ of $W$ compatible with the basis of $V$. The matrix of the form $c$ is obtained from \eqref{(6.44)}.

Let us specify the choice of a basis in $V$ for the case of an indecomposable representation. Assuming a basis of $W$ to be constructed as above, let us write down the quantities that determine the object \eqref{(6.7)}. 

Consider the representation with support \eqref{(6.49)} corresponding to a 1-dimensional space $V_{P_1}$, or with support \eqref{(6.48)} corresponding to a 2-dimensional symplectic space $V_{P_1}$. Select $e_1\in V_{P_1}$ and $e_{n+1}\in V_{\tau P_1}$ so that $b_{P_1}(e_1,e_{n+1})=1$. Set
$$ 
e_2:=h_{P_1}e_1,\ldots,\quad e_n:=h_{P_{n-1}}e_{n-1},\quad e_{n+2}:=h_{\tau P_2}^{-1}e_{n+1},\ldots,\quad e_{2n}:=h_{\tau P_n}^{-1}e_{2n-1}.
$$
Then
$$
m_1=\alpha(P_1),\ldots,m_n=\alpha(P_n),\ m_{n+1}=\alpha(\tau P_1),\ldots,m_{2n}=\alpha(\tau P_n)
$$
and the matrix of the form $c$ acquires the shape\footnote{The superscript $\top$ is used to denote transposed matrices.} 
\begin{equation}\label{(6.53)}
\left(\begin{array}{c|c} 
0&\Jcal_n(0)\vtop to7pt{}\\
\hline
-\Jcal_n(0)^{\top}&0\vbox to14pt{}
\end{array}
\right),
\end{equation}
where $\Jcal_n(0)$ is the Jordan cell with eigenvalue 0.  

Consider a representation with support \eqref{(6.51)} or \eqref{(6.50)}. Set $h:=h_{P_t}\circ h_{P_{t-1}}\circ\ldots \circ  h_{P_{1}}$. In the case of support \eqref{(6.50)} we use the decomposition $V_{P_1}=U_{P_1}\oplus U_{P_1}'$ with the properties described in Lemma~\ref{6.11}. Set
$$ 
U_{P_2}:=h_{P_1}U_{P_1},\quad U_{P_2}':=h_{P_1}U_{P-1}',\ldots,\quad U_{P_t}:=h_{P_{t-1}}U_{P_{t-1}},\quad U_{P_t}':=h_{P_{t-1}}U_{P_{t-1}}'.
$$
Then $\,b_{P_i}(U_{P_i},U_{\tau P_i})=0$, $\,b_{P_i}(U_{P_i}',U_{\tau P_i}')=0$, and the restriction of $b_{P_i}$ to $U_{P_i}\times U'_{\tau P_i}$ gives a non-degenerate pairing. Take any basis $e_1,\ldots,e_n$ of the space $V_{P_1}$ (resp. $U_{P_1}'$) in the case of component \eqref{(6.51)} (resp.~\eqref{(6.50)}) and the dual with respect to $b_{P_1}$ basis $e_{tn+1},\ldots,e_{tn+n}$ of the space $V_{\tau P_1}$ (resp. $U'_{\tau P_1}$). Set
\begin{gather*}
e_{n+1}=h_{P_1}e_1,\ldots,e_{2n}=h_{P_1}e_n,\ldots,e_{tn-n+1}=h_{P_{t-1}}e_{tn-2n+1},\ldots,e_{tn}=h_{P_{t-1}}e_{tn-n},\\
e_{tn+n+1}=h_{\tau P_2}^{-1}e_{tn+1},\ldots,e_{tn+2n}=h_{\tau P_2}^{-1}e_{tn+n},\ldots,e_{2tn-n+1}=h_{\tau P_t}^{-1}e_{2tn-2n+1},\ldots,e_{2tn}=h_{\tau P_t}^{-1}e_{2tn-n}.
\end{gather*}
Then
\begin{gather*}
m_1=\ldots=m_n=\alpha(P_1),\quad\ldots,\quad m_{tn-n+1}=\ldots=m_{tn}=\alpha(P_t),\\
m_{tn+1}=\ldots=m_{tn+n}=\alpha(\tau P_1),\quad\ldots,\quad m_{2tn-n+1}=\ldots=m_{2tn}=\alpha(\tau P_t).
\end{gather*}
Let $(\xi_{ij})$ be the matrix of the endomorphism $h$ (resp. $h|_{U_{P_1}}$) in the basis $e_1,\ldots,e_n$ in the case of representation with support \eqref{(6.51)} (resp. \eqref{(6.50)}). In both cases the matrix of the form $c$ in the basis 
$e_1',\ldots,e_N'$, where $N=2tn$, of the space $W$ obtained from the basis $e_1,\ldots,e_N$ of $V$ as explained earlier is
\begin{equation}
\label{(6.54)}
\left(\begin{array}{cccc|cccc}
&&&&&E\\
\multicolumn{4}{c|}{\raisebox{-4pt}{\Large0}}&&&\ddots\\
&&&&&&&E\\
&&&&\Efrak\vtop to6pt{}\\
\hline
&&&-\Efrak^{\mskip1mu\top}\vbox to14pt{}&\\
-E&&&&\\
&\ddots&&&\multicolumn{4}{c}{\raisebox{6pt}{\Large0}}\\
&&-E&&
\end{array}\right)
\end{equation}
where
$$
\Efrak=\left(\begin{array}{c}
\xi_{11}^{p^{\ell_1}}\ldots \xi_{t1}^{p^{\ell_1}}\\
\hdotsfor{1}\\
\xi_{1t}^{p^{\ell_1}}\ldots \xi_{tt}^{p^{\ell_1}}
\end{array}
\right),\quad\ell_1=\alpha(\tau P_1).
$$
The obtained object of the category $\Acal$ is indecomposable if the matrix 
$(\xi_{ij})$ defines an indecomposable $p^{\alpha(P_1)+\ldots+\alpha(P_t)-\alpha(\tau P_1)-\ldots\alpha(\tau P_t)}$-semilinear endomorphism.

In the case of an algebraically closed field we can specify further. The indecomposable linear endomorphisms are precisely the Jordan blocks. Therefore, in \eqref{(6.54)} we should take
$$
\Efrak=\Jcal_n(\lambda)\ \,\text{with }\,\lambda\ne 0\quad\text{when }\,\alpha(P_1)+\ldots+\alpha(P_t)=\alpha(\tau P_1)+\ldots\alpha(\tau P_t).
$$
On the other hand, for $r\ne 0$ any non-degenerate $p^r$-semilinear endomorphism of a finite-dimensional vector space acts as the identity transformation on a suitable basis. For $r=1$, this is proved in~\cite[Ch.~5, Section~8]{4}. In the general case the proof goes through without changes. In particular, the space of indecomposable such an endomorphism is 1-dimensional. Therefore, we should replace \eqref{(6.54)} with (1s and $-1$s occupy main diagonals of respective blocks)
\begin{equation}
\label{(6.55)}
\left(\begin{array}{cccc|cccc}
&&&&&1\\
\multicolumn{4}{c|}{\raise-4pt\hbox{\Large0}}&&&\ddots \\
&&&&&&&1\\
&&&&1\vtop to 5pt{}\\
\hline
&&&-1&\vbox to 12pt{}\\
-1&&&&\\
&\ddots&&&
\multicolumn{4}{c}{\raise6pt\hbox{\Large0}}\\
&&-1&&
\end{array}
\right)
\end{equation}
when $\ \alpha(P_1)+\ldots+\alpha(P_t)\ne\alpha(\tau P_1)+\ldots+\alpha(\tau P_t)$.

\section{\repl{Normal shape} of the the 1st-type \repl{symplectic} forms}\label{7}

Consider all possible triples $(E,\Fcal,\om)$, where $E$ is a space, $\Fcal$ a flag, $\om\in\Om^2(E,\Fcal)$ a \repl{symplectic} form. Triples $(E,\Fcal,\om)$ and $(E',\Fcal',\om')$ will be called \textit{equivalent} if there exists an isomorphism ${\sigma:\Ocal(E,\Fcal)\to\Ocal(E',\Fcal')}$ which preserves the divided powers and transforms $\om$ to $\om'$.  Preservation of divided powers under $\sigma$ means that for each $f\in\fm(E,\Fcal)$ we have $f^{(r)}\in\Ocal(E,\Fcal)$ if and only if $(\sigma f)^{(r)}\in\Ocal(E',\Fcal')$, and if $f^{(r)}\in\Ocal(E,\Fcal)$, then $\sigma f^{(r)}=(\sigma f)^{(r)}$. (Note that such an isomorphism extends to an isomorphism of complexes $\Om(E,\Fcal)\to\Om(E',\Fcal')$.)

Suppose that $E=E'\oplus E''$ is a decomposition such that 
$$
E_k=(E'\cap E_k)+(E''\cap E_k)\quad\text{for all }\,k.
$$
Define flags $\Fcal'$ (resp. $\Fcal''$) in $E'$ (resp. $E''$) by setting $E_k'=E'\cap E_k$ and $E''_k=E''\cap E_k$. There are canonical embeddings $\Om(E',\Fcal')\to\Om(E,\Fcal)$ and $\Om(E'',\Fcal'')\to\Om(E,\Fcal)$. 

A form $\om\in\Ham^1\Fcal$ will be called \textit{decomposable} if it is conjugate under $G(\Fcal)$ to a form which can be written as $\om'+\om''$ for some $\om'\in\Ham^1\Fcal'$ and $\om''\in\Ham^2\Fcal''$.

To every triple $(E,\Fcal,\om)$ we assign an object \eqref{(6.7)} as follows. 
Set
$$
V=E^*\qquad\text{and}\qquad W=H^1(\Om(\Fcal))^*.
$$
Define the flags $\Gcal$ and $\Hcal$ by the formulas
\begin{align*}
V(k)&=\{v\in V\mid\langle v,E_k\rangle=0\},\\
W(k)&=\{w\in W\mid\langle w,H^1(\Om(\Fcal))_k\rangle=0\} 
\end{align*}
(the flag in $H^1(\Om(\Fcal))$ has been defined in \ref{1}). The~isomorphisms of Proposition~\ref{1.11} induce $p^{k+1}$-semilinear isomorphisms 
$$
V(k+1)/V(k)\cong(E_k/E_{k+1})^*\to\bigl(H^1(\Om(\Fcal))_k/H^1(\Om(\Fcal))_{k+1}\bigr)^*\cong W(k+1)/W(k)
$$
which we take for $\nu^k$. The mapping \eqref{(2.5)} assigns to $\om$ bivectors in $\LaEnd$ and $\bigwedge^2\!H^1(\Om(\Fcal))$ which we identify with forms $b$ and $c$, respectively.

\begin{proposition}\label{7.1}
\repl{Symplectic} forms $\om\in\Om^2(E,\Fcal)$ and $\Tilde\om\in\Om^2(\Tilde E,\Tilde{\Fcal})$ are equivalent if and only if so are the objects $A$ and $\Tilde A$ of the category $\Acal$ corresponding to them.

The form $\om$ is indecomposable if and only if the corresponding object $A$ is indecomposable.
\end{proposition}

\noindent
\textbf{Proof.} Suppose that $A$ and $\Tilde A$ are equivalent. Then the identically indexed factors of the flags $\Gcal$ and $\Tilde\Gcal$, and therefore also such factors of the flags $\Fcal$ and $\Tilde{\Fcal}$, have the same dimension. Hence, we may assume that $\Tilde E=E$ and $\Tilde{\Fcal}=\Fcal$. Then $A$ being equivalent to $\Tilde A$ means precisely that the images of $\om$ and $\Tilde\om$ in the set $\LaEnd\times\bigwedge^2\!H^1(\Om(\Fcal))$ are conjugate under the group $T$ (see the definition at the end of~\ref{2}). Since \eqref{(2.7)} is a bijection, the forms $\om$ and $\Tilde\om$ are conjugate under $G(\Fcal)$.

Suppose that $A\cong A'\oplus A''$ is a decomposable object. The objects $A'$ and $A''$ correspond to some \repl{symplectic} forms $\om'\in\Om^2(E',\Fcal')$ and $\om''\in\Om^2(E'',\Fcal'')$. On $E'\oplus E''$ define a flag $\Fcal'\oplus\Fcal''$ whose spaces are $E_k'\oplus E''_k$. Then $\om'+\om''\in\Om^2(E'\oplus E'',\Fcal'\oplus\Fcal'')$ is a decomposable \repl{symplectic} form, and the corresponding to it object in $\Acal$ is equivalent to $A'\oplus A''$, i.e., to $A$. By the preceding part $\om$ is equivalent to $\om'+\om''$, i.e., is decomposable.\qed 

\medskip
Using the classification of objects in the category $\Acal$ given in~\ref{6} let us formulate our final results. Recall that $\Pfrak$ stands for the set of pairs $(\klov)$,  where $\kov$ and $\lov$ are sequences of positive integers, both of which are either finite, or infinite periodic. The equivalence relation considered on $\Pfrak$ is generated by \eqref{(6.45)}--\eqref{(6.47)}.

\begin{theorem}[$K$ is perfect]\label{7.2}
For each \repl{symplectic} form of the 1st type the equivalent form represented as the sum of indecomposable \repl{symplectic} forms in disjoint groups of indeterminates\footnote{By indeterminates $x_1,\ldots,x_n$ we mean here elements of an $\Fcal$-compatible basis of the space $E$. As explained in \ref{1}, the flag $\Fcal$ determines the \textit{height} $m_i$ of each $x_i$. Since $\Fcal$ is not fixed in this section, a precise description of a \repl{symplectic} form should also include the values of the heights $m_1,\ldots,m_n$, from which the flag can be reconstructed.} is unique up to a permutation of summands and their equivalence.
\end{theorem}

\begin{theorem}[$K$ is perfect]\footnote{This is a new formulation which should be, hopefully, better understandable than the original formulation.}
\label{7.3}
The indecomposable \repl{symplectic} forms of the 1st type split up into disjoint classes which are in a bijective correspondence with the elements of the set\/ $\Pfrak/\mathord\sim\,$.

\smallskip
Any indecomposable 1st-type \repl{symplectic} form is equivalent to a form
$$
\sum_{i<j}\left(a_{ij}+b_{ij}x_i^{(p^{m_i}-1)}x_j^{(p^{m_j}-1)}dx_i\wedge dx_j\right),
$$
where $a_{ij},b_{ij}\in K$, the matrix $(a_{ij})$ is as presented in \eqref{(6.52)}, whereas the shape of the matrix $(b_{ij})$ is determined by the class containing the given form.

\smallskip
Let $(\klov)\in\Pfrak$ be a representative of an element of the set\/ $\Pfrak/\mathord\sim\,$. The corresponding class consists of equivalent \repl{symplectic} forms with the matrix $(b_{ij})$ in their normal shape given by \eqref{(6.53)}, if the sequences 
$$
\kov=(k_1,\ldots,k_n)\quad\text{and}\quad\lov =(\ell_1,\ldots,\ell_n)
$$
are finite. Moreover, the heights of the indeterminates $x_1,\ldots,x_{2n}$ in this case are
$$
m_1=k_1,\ldots,m_n=k_n,\ m_{n+1}=\ell_1,\ldots,m_{2n}=\ell_n.
$$

\smallskip
If the sequences $\kov$, $\lov$ are infinite and $t$ is their least common period, then the matrices $(b_{ij})$ of the forms in the class corresponding to the pair $(\klov)$ are as in \eqref{(6.54)} with the blocks arising from a splitting of the indeterminates $x_1,x_2,\ldots$ into $2t$ groups of equal cardinality. 

All indeterminates in the first of these groups have height $k_1,\ldots\,,$ indeterminates in the $t$-th group have height $k_t$, indeterminates in the $(t+1)$st group have height $\ell_1,\ldots\,,$ and indeterminates  in the last $2t$-th group have height $\ell_t$. Each matrix of shape \eqref{(6.54)} depends on a matrix $(\xi_{ij})$ of an indecomposable $p^{k_1+\ldots+k_t-\ell_1-\ldots-\ell_t}$-semilinear endomorphism. 

Two \repl{symplectic} forms in the same class are equivalent if and only if so are the corresponding semilinear endomorphisms.

\smallskip
If $K$ is algebraically closed, then in the second case, i.e., if $\klov$ are infinite periodic, the matrix $(b_{ij})$ is the matrix \eqref{(6.54)},  where $\Efrak$ is a Jordan cell with nonzero eigenvalue if $k_1+\ldots+k_t=\ell_1+\ldots+\ell_t$, or the matrix \eqref{(6.55)} of size $2t\times 2t$ if $k_1+\ldots+k_t\ne\ell_1+\ldots+\ell_t$.
\end{theorem}

\end{document}